\documentclass[11pt]{amsart}

\usepackage[T1]{fontenc}
\usepackage{lmodern}
\usepackage{amssymb,amsmath}
\allowdisplaybreaks
\usepackage{calc}
\usepackage[inline]{enumitem}
\usepackage{hyphenat}
\usepackage{xcolor}
\usepackage{mathtools}
\usepackage{tikz-cd}
\usetikzlibrary{arrows.meta,decorations.pathreplacing,calc,fadings}

\usepackage[labelformat=empty]{caption}
\usepackage[left=3cm,right=3cm,top=2.8cm,bottom=2.8cm]{geometry}
\usepackage[colorlinks,linkcolor=blue,citecolor=blue,urlcolor=teal]{hyperref}
\usepackage{aliascnt}
\usepackage[nameinlink,noabbrev]{cleveref}

\crefformat{section}{\S#2#1#3}
\Crefformat{section}{\S#2#1#3}
\crefformat{subsection}{\S#2#1#3}
\Crefformat{subsection}{\S#2#1#3}

\newcommand{\Q}{\mathbb{Q}}
\newcommand{\Z}{\mathbb{Z}}
\newcommand{\C}{\mathbb{C}}
\newcommand{\cE}{\mathcal{E}}
\newcommand{\cF}{\mathcal{F}}
\newcommand{\cO}{\mathcal{O}}
\newcommand{\cI}{\mathcal{I}}
\DeclareMathOperator{\Hom}{Hom}
\DeclareMathOperator{\Spec}{Spec}
\DeclareMathOperator{\Spf}{Spf}
\DeclareMathOperator{\Frac}{Frac}
\DeclareMathOperator{\rk}{rk}
\newcommand{\GM}{\mathrm{GM}}
\newcommand{\conj}{\mathrm{conj}}

\theoremstyle{plain}
\newtheorem{prop}{Proposition}[section]
\crefname{prop}{Proposition}{Propositions}
\Crefname{prop}{Proposition}{Propositions}

\newaliascnt{proposition}{prop}
\newtheorem{proposition}[proposition]{Proposition}
\aliascntresetthe{proposition}
\crefname{proposition}{Proposition}{Propositions}
\Crefname{proposition}{Proposition}{Propositions}

\newaliascnt{lemma}{prop}
\newtheorem{lemma}[lemma]{Lemma}
\aliascntresetthe{lemma}
\crefname{lemma}{Lemma}{Lemmas}
\Crefname{lemma}{Lemma}{Lemmas}

\newaliascnt{theorem}{prop}
\newtheorem{theorem}[theorem]{Theorem}
\aliascntresetthe{theorem}
\crefname{theorem}{Theorem}{Theorems}
\Crefname{theorem}{Theorem}{Theorems}

\newaliascnt{corollary}{prop}
\newtheorem{corollary}[corollary]{Corollary}
\aliascntresetthe{corollary}
\crefname{corollary}{Corollary}{Corollaries}
\Crefname{corollary}{Corollary}{Corollaries}

\newaliascnt{conjecture}{prop}

\aliascntresetthe{conjecture}
\crefname{conjecture}{Conjecture}{Conjectures}
\Crefname{conjecture}{Conjecture}{Conjectures}

\theoremstyle{definition}
\newaliascnt{definition}{prop}
\newtheorem{definition}[definition]{Definition}
\aliascntresetthe{definition}
\crefname{definition}{Definition}{Definitions}
\Crefname{definition}{Definition}{Definitions}

\newaliascnt{remark}{prop}
\newtheorem{remark}[remark]{Remark}
\aliascntresetthe{remark}
\crefname{remark}{Remark}{Remarks}
\Crefname{remark}{Remark}{Remarks}

\newaliascnt{example}{prop}
\newtheorem{example}[example]{Example}
\aliascntresetthe{example}
\crefname{example}{Example}{Examples}
\Crefname{example}{Example}{Examples}

\newaliascnt{para}{prop}
\newtheorem{para}[para]{}
\aliascntresetthe{para}
\crefformat{para}{#2#1#3}
\Crefformat{para}{#2#1#3}

\theoremstyle{plain}
\newtheorem{claim}{Claim}[prop]
\crefname{claim}{Claim}{Claims}
\Crefname{claim}{Claim}{Claims}

\renewcommand{\theequation}{\theprop.\arabic{equation}}
\makeatletter
\@addtoreset{equation}{prop}
\makeatother
\setlist[enumerate]{label=\textup{(\arabic*)},leftmargin=*}
\setlist[enumerate,2]{label=\textup{(\alph*)},leftmargin=*}

\title[$K3$ Picard--Fuchs equations with initial conditions in the Hodge line]{$K3$ Picard--Fuchs equations with initial conditions\\
in the Hodge line}
\author{Fei Ren}
\address{Bergische Universit\"at Wuppertal, Gau\ss str.\ 20, D-42119 Wuppertal, Germany}
\email{renfei@uni-wuppertal.de}

\providecommand{\eq}[2]{\begin{equation}\label{#1}#2 \end{equation}}

\providecommand{\beq}{\begin{equation}}
\providecommand{\eeq}{\end{equation}}

\providecommand{\F}{\mathbb{F}}

\providecommand{\W}{\mathbf{W}}

\providecommand{\Ker}{\operatorname{Ker}}

\providecommand{\dR}{\operatorname{dR}}

\providecommand{\leftrarrows}{\mathrel{\raise.75ex\hbox{\oalign{
  $\scriptstyle\leftarrow$\cr
  \vrule width0pt height.5ex$\hfil\scriptstyle\relbar$\cr}}}}
\providecommand{\lrightarrows}{\mathrel{\raise.75ex\hbox{\oalign{
  $\scriptstyle\relbar$\hfil\cr
  $\scriptstyle\vrule width0pt height.5ex\smash\rightarrow$\cr}}}}
\providecommand{\Rrelbar}{\mathrel{\raise.75ex\hbox{\oalign{
  $\scriptstyle\relbar$\cr
  \vrule width0pt height.5ex$\scriptstyle\relbar$}}}}

\makeatletter
\providecommand{\leftrightarrowsfill@}{\arrowfill@\leftrarrows\Rrelbar\lrightarrows}
\providecommand{\xleftrightarrows}[2][]{\ext@arrow 3399\leftrightarrowsfill@{#1}{#2}}
\makeatother

\providecommand{\bthm}{\begin{theorem}}
\providecommand{\ethm}{\end{theorem}}
\providecommand{\bcor}{\begin{corollary}}
\providecommand{\ecor}{\end{corollary}}
\providecommand{\bprop}{\begin{proposition}}
\providecommand{\eprop}{\end{proposition}}
\providecommand{\ble}{\begin{lemma}}
\providecommand{\ele}{\end{lemma}}
\providecommand{\bex}{\begin{example}}
\providecommand{\eex}{\end{example}}
\providecommand{\bexc}{\begin{exercise}}
\providecommand{\eexc}{\end{exercise}}
\providecommand{\brmk}{\begin{remark}}
\providecommand{\ermk}{\end{remark}}
\providecommand{\bdefn}{\begin{definition}}
\providecommand{\edefn}{\end{definition}}
\providecommand{\bpf}{\begin{proof}}
\providecommand{\epf}{\end{proof}}
\providecommand{\benu}{\begin{enumerate}}
\providecommand{\eenu}{\end{enumerate}}
\providecommand{\bit}{\begin{itemize}}
\providecommand{\eit}{\end{itemize}}
\providecommand{\bqn}{\begin{qn}}
\providecommand{\eqn}{\end{qn}}
\providecommand{\D}{\mathbb{D}}

\renewcommand{\hat}{\widehat}
\renewcommand{\tilde}{\widetilde}
\begin{document}

\begin{abstract}
Lam and Litt conjectured that a formal horizontal section of a Gauss--Manin connection is algebraic if and only if it is integral, if and only if it is $\omega(p)$-integral.
We prove it for families of $K3$ surfaces over a finitely generated $\Z$-algebra in a number field, when the initial condition lies in the Hodge line.
In this case, all three conditions are equivalent to the vanishing of the Higgs field on the Hodge line.
Along the way, we show that a smooth projective family $X\to S$ with $H^0(X_s,T_{X_s})=0$ is constant on an \'{e}tale neighborhood of $s$ if and only if its Kodaira--Spencer class vanishes near $s$.
\end{abstract}

\maketitle

\tableofcontents

\section{Introduction}
\label{sec:intro}

Let $R\subset \C$ be a finitely generated $\Z$-algebra, with fraction field $K$, and let $S$ be a smooth $R$-scheme.
Let $\pi\colon X\to S$ be a smooth projective morphism, and write
\[
(\cE,\nabla)=\bigl(R^{m}\pi_*\Omega^\bullet_{X/S},\nabla_{\GM}\bigr)
\]
for a Gauss--Manin system of $\pi$ (a Picard--Fuchs equation in the sense of \cite[\S1.2.1]{LL25}).
Fix $s\in S(R)$ and a morphism $v\colon\Spec R\to E$ over $s$, where $E$ is the total space of $\cE$.
There is a unique formal horizontal section of $(\cE,\nabla)_K$ through $v_K$.
Lam--Litt, in \cite[Conjecture 2.3.1]{LL25}, conjecture that the formal horizontal section through $v_K$ is algebraic if and only if it is integral, if and only if it is $\omega(p)$-integral in the sense of \cite[Definition 2.2.5]{LL25}.
The linear version of this conjecture, stated in \textit{loc.\ cit.}\ Conjecture 1.3.1, is strictly stronger than the Grothendieck--Katz $p$-curvature conjecture: see \textit{loc.\ cit.}\ \S1.2.6.
They prove that an $\omega(p)$-integral formal flat section of a Picard--Fuchs system is algebraic if the initial condition lies in the image of the cycle class map \cite[Theorem 1.3.3]{LL25}.
They also prove the linear conjecture for $m=1$ for every nonzero initial condition, when $\pi$ is a non-isotrivial family of elliptic curves \cite[Theorem 4.1.1]{LL25}.

The elliptic argument uses that, for a smooth projective family $\pi\colon X\to S$ of elliptic curves in characteristic $p$, the conjugate filtration on $R^1\pi_*\Omega^\bullet_{X/S}$ has two graded pieces $\mathrm{gr}^0_{\conj}$ and $\mathrm{gr}^1_{\conj}$, both of rank one.
A class whose reduction at $s$ is nonzero in $\mathrm{gr}^0_{\conj}$ therefore generates that line near $s$, so vanishing to order $\omega(p)$ of the graded $p$-curvature $\mathrm{gr}^0_{\conj}\psi_p$ on the class is its vanishing to that order on the whole graded piece.
The identification, up to sign, of graded conjugate $p$-curvature with the Frobenius pullback of the Higgs field \cite[Theorem 4.2.1]{LL25} (citing \cite[Theorem 3.2]{Kat72}) then makes $\mathrm{gr}^1_{\mathrm{Hdg}}\nabla_{\GM}$ of the Gauss--Manin connection $\nabla_{\GM}$ vanish to order $\omega(p)/p$ at a Zariski-dense set of primes, hence identically, so the Hodge line $F^1_{\mathrm{Hdg}}$ is preserved by $\nabla_{\GM}$, i.e., it is horizontal and underlies a sub-local system.
The monodromy of an elliptic family preserves the cup product pairing $\begin{pmatrix}0&1\\-1&0\end{pmatrix}$ on $R^1\pi^{\mathrm{an}}_*\mathbb{Z}$, so it lands in $\mathrm{SL}_2$.
By \cite[Corollaire (4.2.9)]{Del71} the identity component of its Zariski closure is semisimple.
The only connected semisimple algebraic subgroups of $\mathrm{SL}_2$ are the trivial group and $\mathrm{SL}_2$.
For a non-isotrivial family the monodromy is infinite, so that identity component can only be $\mathrm{SL}_2$.
The group $\mathrm{SL}_2$ has no invariant line in $\C^2$.
If $F^1_{\mathrm{Hdg}}$ were horizontal, its fiber at a base point would be a monodromy-invariant line, hence invariant under the Zariski closure of the monodromy and its identity component $\mathrm{SL}_2$; so $F^1_{\mathrm{Hdg}}$ is not horizontal.

The same rank-one phenomenon occurs for $R^2\pi_*\Omega^\bullet$ of a $K3$ surface, but only in the Hodge line $F_{\mathrm{Hdg}}^2$.
A $K3$ surface $Y$ over a perfect field $k$ and its Frobenius twist $Y'$ satisfy $\dim_k H^0(Y,\Omega^2_{Y/k})=\dim_k H^0(Y',\Omega^2_{Y'/k})=1$, so $\mathrm{gr}^0_{\conj}H^2_{\dR}(Y/k)\simeq \Gamma(Y',\Omega^2_{Y'/k})$ via the Cartier isomorphism is one-dimensional, while $\mathrm{gr}^1_{\conj}H^2_{\dR}(Y/k)$ has dimension $20$.
A single vector in $\mathrm{gr}^1_{F_{\conj}}H^2_{\dR}(X_{k(\mathfrak{p})}/S_{k(\mathfrak{p})})$ does not force the $p$-curvature of that piece to vanish.
The elliptic argument, applied to $R^2\pi_*\Omega^\bullet$, therefore shows that $\theta=0$ for an initial condition $v\in F_{\mathrm{Hdg}}^2$ satisfying the hypothesis of \cite[Conjecture 1.3.1]{LL25}, i.e., such that the formal horizontal section of $(\cE,\nabla_{\GM})$ through $v$ is $\omega(p)$-integral, provided the reduction of $v$ modulo $\mathfrak{p}$ generates $\mathrm{gr}^0_{F_{\conj}}H^2_{\dR}(X_{k(\mathfrak{p})}/S_{k(\mathfrak{p})})$ near $s$ for a Zariski-dense set of maximal ideals $\mathfrak{p}$ of $R$.
For a nonzero $v\in H^1_{\dR}(X_s/R)$, Lam--Litt show that the reduction of $v$ modulo $\mathfrak{p}$ does not lie in $F^1_{\conj}H^1_{\dR}(X_{s,k(\mathfrak{p})}/k(\mathfrak{p}))$ for a Zariski-dense set of maximal ideals $\mathfrak{p}$ of $R$ \cite[Lemma 4.2.2]{LL25}, using the Grothendieck--Katz conjecture for extensions of $(\cO,d)$ by itself, which is known.
For a general class in $R^2\pi_*\Omega^\bullet$ the analogous lemma seems to be the crucial difficulty, and we do not prove that lemma here.
Instead, for a class in the Hodge line there is a different input: on an ordinary $K3$ the Hodge and conjugate flags on $R^2\pi_*\Omega^\bullet$ are complementary \cite[IV, Th\'eor\`eme (4.13) (i)(vi)]{IR83}, and a $K3$ over a number field has ordinary good reduction at a density one set of places after a finite extension \cite[Theorem 0.1]{BZ09}.

In this article we carry out the strategy described above, and thereby obtain a new class of initial conditions for which \cite[Conjecture 1.3.1]{LL25} holds for non-isotrivial families of $K3$ surfaces.

\begin{theorem}[{\Cref{thm:main}}]
\label{thm:main-intro}
Let $K$ be a number field and $R\subset K$ a finitely generated $\Z$-algebra.
Let $S$ be a smooth integral $R$-scheme of finite type, and let $\pi\colon X\to S$ be a smooth projective morphism whose fibers are $K3$ surfaces.
Let
\[
(\cE,\nabla)=\bigl(R^2\pi_*\Omega^\bullet_{X/S},\nabla_{\GM}\bigr),
\]
and write
\[
\theta\colon F_{\mathrm{Hdg}}^2\cE\longrightarrow \bigl(F_{\mathrm{Hdg}}^1\cE/F_{\mathrm{Hdg}}^2\cE\bigr)\otimes_{\cO_S}\Omega^1_{S/R}
\]
for the Higgs field of the Hodge line, i.e., the associated graded of $\nabla_{\GM}$ in filtration degree $2$.
Assume:
\begin{enumerate}[label=\textup{(\alph*)}]
\item $\theta$ is not identically zero.
\item $s\in S(R)$ and $v\colon\Spec R\to E$ is a morphism over $s$ that factors through $F_{\mathrm{Hdg}}^2E$ and is not the zero section.
\end{enumerate}
Then the formal horizontal section of $(\cE,\nabla)_K$ through $v_K$ is not $\omega(p)$-integral.
\end{theorem}

We assume that $R$ is contained in the number field $K$ because the proof of \Cref{thm:main} uses \cite[Theorem 0.1]{BZ09}: that theorem is a density one statement for places of a number field, not a Zariski density statement on $\mathrm{Spec}(R)$ for a general finitely generated $\Z$-algebra $R$.

Using deformation theory, we upgrade \Cref{thm:main-intro} to an equivalence, and along the way we prove the following characterization of local isotriviality, which might be of independent interest.

\begin{theorem}[{\Cref{cor:ks-iff}}]
\label{thm:ks-iff-intro}
Let $K$ be a field of characteristic $0$, let $S$ be a smooth $K$-scheme of finite type, let $s\in S(K)$, and let $\pi\colon X\to S$ be a smooth projective morphism with $H^0(X_s,T_{X_s})=0$.
Then the following are equivalent.
\begin{enumerate}
\item The Kodaira--Spencer class $\mathrm{KS}_{\pi}$ vanishes on an open neighborhood of $s$.
\item There exist an \'{e}tale morphism $S'\to S$ with $S'$ connected and of finite type over $K$, and a $K$-point $s'\in S'(K)$ over $s$, such that, over $S'$,
\[
X\times_S S'\simeq X_s\times_K S'.
\]
\end{enumerate}
\end{theorem}

The following theorem is the main result of this article.
In particular, we give a proof of (1)$\Rightarrow$(2) that does not rely on variations of Hodge structures.

\begin{theorem}[{\Cref{thm:hodge-line-iff}}]
\label{thm:hodge-line-iff-intro}
Let $K$, $R$, $S$, $\pi$, $s$, and $v$ be as in \Cref{thm:main-intro}, but without assumption (a).
Write $\tilde{v}_K$ for the formal horizontal section of $(\cE,\nabla)_K$ through $v_K$.
The following are equivalent.
\begin{enumerate}
\item
The Higgs field $\theta$ is identically zero.
\item
The section $\tilde{v}_K$ is algebraic.
\item
The section $\tilde{v}_K$ is integral.
\item
The section $\tilde{v}_K$ is $\omega(p)$-integral.
\end{enumerate}
\end{theorem}

The $K3$-specific inputs of this article are used exclusively in \Cref{sec:theorem} and \Cref{sec:dropping}, together with the standard facts on $K3$ surfaces recalled in \Cref{sec:ordinary}.
In particular, they are as follows.
\begin{enumerate}[label=\textup{(\alph*)}]
\item
\label{item:k3-main}
For \Cref{thm:main-intro} (equivalently \Cref{thm:main}):
\begin{enumerate}[label=\textup{(\roman*)}]
\item the following two results of \cite{BZ09} for $K3$ surfaces, used in Step~1 of the proof of \Cref{thm:main}:
\begin{enumerate}[label=\textup{(\arabic*)}]
\item density one ordinary good reduction after a finite extension \cite[Theorem~0.1]{BZ09} (\Cref{thm:bz}), applied to $X_{s_K}$;
\item ordinarity of a $K3$ surface over a finite field depends only on the characteristic polynomial of geometric Frobenius \cite[Lemma~1.1 and \S3.1]{BZ09}, used in \Cref{claim:good-primes} to pass from the model of \cite{BZ09} to $X_{s,k(\mathfrak p)}$;
\end{enumerate}
\item $\dim_{\kappa(s)}H^0(X_s,\Omega^2_{X_s/\kappa(s)})=1$ for every $s\in S$, so that $F_{\mathrm{Hdg}}^2$ and $\mathrm{gr}^0_{F_{\conj}}$ are lines (\Cref{ex:K3-fzip}\ref{ex:K3-family}), used in Steps~0 and~2 of the proof of \Cref{thm:main};
\item Hodge--de Rham degeneration at $E_1$ for a $K3$ surface over a perfect field of characteristic $p>0$ (\Cref{ex:K3-fzip}), applied in Step~2 of the proof of \Cref{thm:main} to the reductions $X_{s,k(\mathfrak p)}$;
\item for an ordinary $K3$ reduction, \eqref{eq:opposite-n2-i2}, the case $m=2$ of \Cref{lem:opposite} for ordinary $F$-zips, applies.
\end{enumerate}
The number $b_2=22$ is not used.
$\mathrm{gr}^1_{F_{\conj}}$, of rank $20$, is likewise unused for the main theorem: only that $\mathrm{gr}^0_{F_{\conj}}$ is a line enters the argument.
The rank $20$ of $F^1_{\mathrm{Hdg}}\cE/F^2_{\mathrm{Hdg}}\cE$ enters the proof of \Cref{thm:main} only as the length of a basis.
\item
\label{item:k3-isotrivial}
For (1)$\Rightarrow$(2) in \Cref{thm:hodge-line-iff} (and \Cref{lem:isotrivial-hodge-line} directly before it):
\begin{enumerate}[label=\textup{(\roman*)}]
\item the isomorphism $T_{X/S}\xrightarrow{\sim}\Omega^1_{X/S}$ given by a nowhere vanishing $2$-form, which exists locally on $S$ because $\Omega^2_{X_t/\kappa(t)}\simeq\cO_{X_t}$ for every $t\in S$.
Combined with \cite[(1.4.1.7.2)]{Kat72}, it shows that Kodaira--Spencer vanishes exactly when the Higgs field of the Hodge line does;
\item vanishing of $H^0(X_t,T_{X_t})$ for a fiber $X_t$ of $\pi$ \cite[Theorem 9.5.1]{Huy16}, which gives the hypothesis $H^0(Y,T_Y)=0$ of \Cref{thm:ks-isotrivial} for $Y=X_{s_K}$;
\item $\dim_KH^0(Y,\Omega^2_{Y/K})=1$, so that $F^2_{\mathrm{Hdg}}$ of the constant family in \Cref{lem:isotrivial-hodge-line}\ref{item:isohl-3} is a line bundle, and $(F_{\mathrm{Hdg}}^2\cE,\nabla_{\GM})_K$ becomes isomorphic to $(\cO_{S'},d)$ in the proof of \Cref{thm:hodge-line-iff}.
\end{enumerate}
\end{enumerate}

\medskip
\noindent \textbf{AI disclosure.} Following what has become common practice, the author includes the AI disclosure below, although he does not think such a rule is going to work.
This project started with an observation (now \Cref{thm:main}) during a seminar in BUW on \cite{LL25}.
Noticing the comments of Litt, the author decided to explore it further based on this observation.
The author had intense discussions with AI chatbots in this process (about 14\,800 messages, including discussions towards some substantial extensions of this paper); he enjoyed them a lot, and they were pleasant and fruitful.
In particular, the chatbots greatly helped in navigating the deformation theory literature, which led to the exposition in \Cref{sec:deformation-theory}.
They identified the subtle difference between the author's definition of $R[[x_1,\dots,x_d]]^g$ (see \Cref{para:g-integral}; this, especially the structural identity \eqref{eq:g-series}, was a priori a misread) and that in \cite[\S2.2]{LL25}, and supplied a comparison lemma (\Cref{lem:omega-comparison} and \Cref{rmk:omega-comparison}).
They also helped the author clear up \Cref{rmk:local-freeness} and pointed out that it does not depend on \cite{DI87}.
They also drew the nice pictures in \Crefrange{fig:reduction}{fig:schemes-61}.
In addition, they made useful comments on many other minor points in proofreading earlier versions of this article: English, tightening specific arguments, and identification of missing or redundant assumptions.

\section{Recollections}
\label{sec:recollections}

Throughout, $R$ is an integral domain of characteristic $0$, finitely generated as a $\Z$-algebra, and $K=\Frac(R)$.

\subsection{Algebraicity and integrality of formal horizontal sections}
\label{subsec:integrality}

Let $\cE$ be a quasi-coherent $\cO_S$-module on a smooth $R$-scheme $S$.
Write $T_S:=(\Omega^1_{S/R})^\vee$.

\begin{definition}[Connection]
\label{def:connection}
\begin{enumerate}
\item
A \emph{connection} on $\cE$ is an $R$-linear map of sheaves
\[
\nabla\colon\cE\to\Omega^1_{S/R}\otimes_{\cO_S}\cE
\]
satisfying the Leibniz rule $\nabla(fs)=f\nabla(s)+df\otimes s$ for local sections $f$ of $\cO_S$ and $s$ of $\cE$.
\item
For an open $U\subset S$ and $\partial\in T_S(U)$, write
\[
\nabla_\partial:=(\partial\otimes\mathrm{id}_{\cE})\circ\nabla\colon\cE|_U\to\cE|_U.
\]
\end{enumerate}
\end{definition}

\begin{definition}[Horizontal section]
\label{def:horizontal}
Let $\nabla$ be a connection on $\cE$.
A local section $s\in\cE(U)$ is \emph{horizontal} if $\nabla s=0$ in $(\Omega^1_{S/R}\otimes_{\cO_S}\cE)(U)$.
\end{definition}

\begin{definition}[Curvature, $p$-curvature]
\label{def:curvature}
A connection $\nabla\colon\cE\to\Omega^1_{S/R}\otimes_{\cO_S}\cE$ extends to $\nabla\colon\Omega^1_{S/R}\otimes_{\cO_S}\cE\to\Omega^2_{S/R}\otimes_{\cO_S}\cE$ by $\nabla(\omega\otimes s)=d\omega\otimes s-\omega\wedge\nabla s$.
\begin{enumerate}
\item
The \emph{curvature} of $\nabla$ is the $\cO_S$-linear map
\[
\nabla^2:=\nabla\circ\nabla\colon\cE\to\Omega^2_{S/R}\otimes_{\cO_S}\cE.
\]
A connection is \emph{integrable}, or \emph{flat}, when $\nabla^2=0$.
\item
Let the base have characteristic $p>0$, let $\nabla$ be integrable, and let $F_{\mathrm{abs}}\colon S\to S$ be the absolute Frobenius.
The \emph{$p$-curvature} of $\nabla$ is the $\cO_S$-linear map
\[
\psi_p\colon\cE\to F_{\mathrm{abs}}^*\Omega^1_{S/R}\otimes_{\cO_S}\cE
\]
of \cite[(3.0.3)]{Kat72}.
Equivalently, by transposition, it is the map
\[
\psi_p\colon T_S\to\operatorname{End}_{\cO_S}(\cE),\qquad
\psi_p(\partial)(s)=\nabla_\partial^p(s)-\nabla_{\partial^p}(s),
\]
where $\partial^p$ is the $p$-th power of the derivation $\partial$.
It is $p$-linear in the vector field, $\psi_p(f\partial)=f^p\psi_p(\partial)$.
\end{enumerate}
\end{definition}

From now on $\cE$ is a vector bundle, i.e., a locally free $\cO_S$-module of finite rank.
Write
\[
E:=\underline{\Spec}_S\operatorname{Sym}(\cE^\vee)\xrightarrow{\tau} S
\]
for its total space with structure morphism.
For an open $U\subset S$, morphisms $\sigma\colon U\to E$ with $\tau\circ\sigma=\mathrm{id}_U$ correspond bijectively to elements of $\cE(U)$, and we call both of them \emph{local sections} of $\cE$ over $U$.
The formation of $E$ commutes with base change, so $E_K$ is the total space of $\cE_K$.
For a subbundle $\cF\subset\cE$, i.e., a locally direct summand, the surjection $\cE^\vee\to\cF^\vee$ makes the total space of $\cF$ a closed subscheme of $E$.
We write $F_{\mathrm{Hdg}}^iE\subset E$ for the total space of $F_{\mathrm{Hdg}}^i\cE$.
Write $T_{E}:=(\Omega^1_{E/R})^\vee$.

The following computation is standard; we include it for completeness.

\begin{lemma}[Horizontal lift of vector fields]
\label{lem:connection-splitting}
Let $\nabla$ be a connection on $\cE$, let $U\subset S$ be an open subset over which $\cE|_U$ is free with basis $e_1,\dots,e_r$, and let $t_1,\dots,t_r\in\cE^\vee(U)\subset\cO_E(\tau^{-1}U)$ be the dual basis, so that $\cO_E(\tau^{-1}U)=\cO_S(U)[t_1,\dots,t_r]$.
Write
\[
\nabla e_j=\sum_i\omega_{ij}\otimes e_i\quad\text{with }\omega_{ij}\in\Omega^1_{S/R}(U).
\]
For $\partial\in T_S(U)$, let $\partial^h\in T_E(\tau^{-1}U)$ be the vector field
\begin{equation}\label{eq:horizontal-lift}
\partial^h:=\partial-\sum_{i,j}\omega_{ij}(\partial)\,t_j\,\frac{\partial}{\partial t_i},
\end{equation}
where $\partial$ acts on $\cO_E(\tau^{-1}U)=\cO_S(U)[t_1,\dots,t_r]$ coefficientwise.
Then the following hold.
\begin{enumerate}
\item The vector field $\partial^h$ does not depend on the basis.
\item For every $f\in\cO_E(\tau^{-1}U)$ and every $\sigma\in\cE(U)$, viewed as a morphism $\sigma\colon U\to E$, such that $\nabla_\partial\sigma=0$, we have
\[
\sigma^*(\partial^h f)=\partial(\sigma^*f)
\]
in $\cO_S(U)$.
The same holds for $\sigma\in\cE(\hat U_Z)$, where $\hat U_Z$ is the formal completion of $U$ along a closed subscheme $Z\subset U$, with the equality in $\Gamma(\hat U_Z,\cO_{\hat U_Z})$.
\end{enumerate}
\end{lemma}

\begin{proof}
(1) If $e'_j=\sum_kg_{kj}e_k$ is another basis, write $\nabla e'_j=\sum_i\omega'_{ij}\otimes e'_i$, let $t'_1,\dots,t'_r$ be the coordinates dual to $e'_1,\dots,e'_r$, and set $\omega=(\omega_{ij})$, $\omega'=(\omega'_{ij})$, $g=(g_{kj})$, $t=(t_1,\dots,t_r)^{\mathrm T}$ and $t'=(t'_1,\dots,t'_r)^{\mathrm T}$.
Then
\[\omega'=g^{-1}\omega g+g^{-1}dg\quad\text{and}\quad t=gt',\]
and substituting into the formula for $\partial^h$ shows that it is unchanged.

(2) Write $\sigma=\sum_i\sigma_ie_i$.
By definition of $\nabla_\partial$, we have $\nabla_\partial e_j=\sum_i\omega_{ij}(\partial)e_i$, so by the Leibniz rule,
\[
0=\nabla_\partial\sigma=\sum_j\bigl(\partial(\sigma_j)e_j+\sigma_j\nabla_\partial e_j\bigr)=\sum_i\Bigl(\partial(\sigma_i)+\sum_j\omega_{ij}(\partial)\sigma_j\Bigr)e_i,
\]
so $\partial(\sigma_i)=-\sum_j\omega_{ij}(\partial)\sigma_j$ for all $i$.
By the chain rule,
\[
\partial(\sigma^*f)=\sigma^*(\partial f)+\sum_i\sigma^*\Bigl(\frac{\partial f}{\partial t_i}\Bigr)\partial(\sigma_i)=\sigma^*(\partial^h f).
\]
\end{proof}

\begin{definition}[Formal section]
\label{def:leaf-notions}
Let $\nabla$ be an integrable connection on $\cE$, let $s\in S(R)$, and let $v\colon\Spec R\to E$ be a morphism with $\tau\circ v=s$, and let $v_K\colon\Spec K\to E_K$ be its generic fiber.
Write
\[
\hat S_{s_K}:=\Spf\hat{\cO}_{S_K,s_K},
\]
and let $d$ be the relative dimension of $S$ over $R$ at $s$.
A section of a smooth morphism is a regular immersion \cite[\href{https://stacks.math.columbia.edu/tag/067R}{Tag 067R}]{SP}.
So the ideal sheaf of $s$ is generated by a regular sequence on a Zariski open neighborhood of each point of $s(\Spec R)$.
After replacing $R$ by $R[1/a]$ for some nonzero $a\in R$, we may fix a regular sequence $x_1,\dots,x_d$ generating it, and a basis $e_1,\dots,e_r$ of $\cE$, on a single affine open neighborhood $U_s$ of $s(\Spec R)$. 
The completion of $\cO(U_s)$ with respect to the ideal $(x_1,\dots,x_d)$ of $s$ is $R[[x_1,\dots,x_d]]$.
Then \[\hat S_{s_K}\simeq\Spf\bigl(K[[x_1,\dots,x_d]]\bigr).\]
A section of $\cE_K$ over $\hat S_{s_K}$, viewed as a morphism
\[
\Spf\bigl(K[[x_1,\dots,x_d]]\bigr)\to E_K
\]
over $S_K$, is called a \emph{formal section}.
There is a unique formal section $\tilde v_K$ with
\[
\nabla\tilde v_K=0\quad\text{and}\quad\tilde v_K(s_K)=v_K,
\]
by Taylor expansion \cite[Proposition 2.2.1]{LL25} (citing \cite[\S3.4.1]{Bos01}); we call it the \emph{formal horizontal section} of $(\cE,\nabla)_K$ through $v_K$.
\end{definition}

\begin{para}
\label{para:g-integral}
Let $\mathcal{I}=(x_1,\dots,x_d)\subset R[[x_1,\dots,x_d]]$.
Define
\[
R[[x_1,\dots,x_d]]\hat{\otimes}\Q
\]
to be the completion of $R[[x_1,\dots,x_d]]\otimes_{\Z}\Q$ with respect to the $\mathcal{I}$-adic topology.
For a prime $p$ and $m\geq0$, let
\begin{equation}
\label{eq:Im-p}
\mathcal{I}^m\hat{\otimes}\Z[1/p]:=(x_1,\dots,x_d)^mR[1/p][[x_1,\dots,x_d]]\subset R[[x_1,\dots,x_d]]\hat{\otimes}\Q.
\end{equation}
It is an ideal of $R[1/p][[x_1,\dots,x_d]]$ and an $R[[x_1,\dots,x_d]]$-submodule of $R[[x_1,\dots,x_d]]\hat{\otimes}\Q$, but in general not an ideal of $R[[x_1,\dots,x_d]]\hat{\otimes}\Q$.

Let $g\colon \{\mathrm{primes}\}\to\Z_{\geq 0}$.
Consider the $R[[x_1,\dots,x_d]]$-submodule
\[
H^g:=\sum_{p\text{ prime}}\mathcal{I}^{g(p)+1}\hat{\otimes}\Z[1/p]{}=\sum_{p\text{ prime}}(x_1,\dots,x_d)^{g(p)+1}R[1/p][[x_1,\dots,x_d]]
\]
of $R[[x_1,\dots,x_d]]\hat{\otimes}\Q$.
The ideal
\[
J^g:=\bigl(\mathcal{I}+H^g\bigr)\cdot\bigl(R[[x_1,\dots,x_d]]\hat{\otimes}\Q\bigr)
\]
of $R[[x_1,\dots,x_d]]\hat{\otimes}\Q$ defines the $J^g$-adic topology on that ring.
We have
\begin{equation}
\label{eq:Jg}
J^g=\mathcal{I}\cdot\bigl(R[[x_1,\dots,x_d]]\hat{\otimes}\Q\bigr).
\end{equation}
The equality holds because every summand of $H^g$ lies in $\mathcal{I}\cdot\bigl(R[[x_1,\dots,x_d]]\hat{\otimes}\Q\bigr)$ by \eqref{eq:Im-p}, as $g(p)+1\geq1$.
Hence the $J^g$-adic and $\mathcal{I}$-adic topologies on $R[[x_1,\dots,x_d]]\hat{\otimes}\Q$ coincide.
Let
\begin{equation}
\label{eq:Bg}
B^g:=\text{the subring of }R[[x_1,\dots,x_d]]\hat{\otimes}\Q\text{ generated by }R[[x_1,\dots,x_d]]\text{ and }H^g.
\end{equation}
Define
\[
R[[x_1,\dots,x_d]]^g
\]
to be the closure of $B^g$ in the $J^g$-adic topology.
This closure is a subset, so
\[
R[[x_1,\dots,x_d]]^g\subset R[[x_1,\dots,x_d]]\hat{\otimes}\Q.
\]
We equip $R[[x_1,\dots,x_d]]^g$ with the topology induced from the $\mathcal{I}$-adic topology of the right hand side, and $\Spf\bigl(R[[x_1,\dots,x_d]]^g\bigr)$ denotes its formal spectrum for this topology.
For $n\geq0$, let
\begin{equation}
\label{eq:Rn}
R_n=R[1/p\mid p\text{ prime such that }n\geq g(p)+1]
\end{equation}
be the subring of $K$ generated by $R$ and the $1/p$ for the primes $p$ with $n\geq g(p)+1$.
An element of $R[[x_1,\dots,x_d]]^g$ is a power series over $K$ in the variables $x_1,\dots,x_d$ whose coefficients of total degree $n$ lie in $R_n$ for every $n\geq0$.
In other words, writing $x^\alpha=x_1^{\alpha_1}\cdots x_d^{\alpha_d}$ and $|\alpha|=\alpha_1+\dots+\alpha_d$ for $\alpha\in\Z_{\geq0}^d$, we have
\begin{equation}
\label{eq:g-series}
R[[x_1,\dots,x_d]]^g=\prod_{n\geq0}\ \bigoplus_{|\alpha|=n}R_n\,x^\alpha.
\end{equation}
Indeed, since the $J^g$-adic and $\mathcal{I}$-adic topologies on $R[[x_1,\dots,x_d]]\hat{\otimes}\Q$ coincide, an element lies in the closure of $B^g$ in the $J^g$-adic topology if and only if each of its truncations modulo $\mathcal{I}^{n+1}$ is the truncation of an element of $B^g$, and these truncations are exactly the polynomials whose coefficients of total degree $m\leq n$ lie in $R_m$.
\end{para}

Lam--Litt instead complete the generated subring in its own adic topology \cite[Definition 2.2.5]{LL25}.
The following lemma compares the two conventions; they give the same notion of $\omega(p)$-integrality.

\begin{lemma}[{Comparison with \cite[\S2.2]{LL25}}]
\label{lem:omega-comparison}
Let $\mathcal{I}$, $g$, $H^g$, $J^g$ and $B^g$ be as in \Cref{para:g-integral}, and use the same notation for any other function $h\colon\{\mathrm{primes}\}\to\Z_{\geq0}$.
Let
\[
J_g:=\mathcal{I}+H^g\quad\text{and}\quad\hat{B}^g:=\varprojlim_a B^g/J_g^a.
\]
This $J_g$ is the ideal of \cite[\S2.2]{LL25}.
The canonical maps to finite $\mathcal{I}$-adic quotients induce a map
\begin{equation}
\label{eq:rho-g}
\rho^g\colon\hat{B}^g\longrightarrow R[[x_1,\dots,x_d]]\hat{\otimes}\Q.
\end{equation}
Then the following hold.
\begin{enumerate}[label=\textup{(\arabic*)},ref=\arabic*]
\item\label{item:literal-to-degree} $\operatorname{Im}(\rho^g)\subset R[[x_1,\dots,x_d]]^g$.
\item\label{item:degree-to-literal} If $h(p)=\lfloor g(p)/2\rfloor$, then $R[[x_1,\dots,x_d]]^g\subset\operatorname{Im}(\rho^h)$.
\end{enumerate}
\end{lemma}

\begin{proof}
\begin{claim}\label{claim:Hg-square}
We have $(H^g)^2\subset\mathcal{I}H^g$ as $R[[x_1,\dots,x_d]]$-submodules of $R[[x_1,\dots,x_d]]\hat{\otimes}\Q$.
\end{claim}
\begin{proof}[Proof of \Cref{claim:Hg-square}]
Let $p$ and $q$ be primes, and let $f\in\mathcal{I}^{g(p)+1}\hat{\otimes}\Z[1/p]$ and $f'\in\mathcal{I}^{g(q)+1}\hat{\otimes}\Z[1/q]$.
If $p=q$, then $ff'\in\mathcal{I}\bigl(\mathcal{I}^{g(p)+1}\hat{\otimes}\Z[1/p]\bigr)\subset\mathcal{I}H^g$, since $g(q)+1\geq1$.
If $p\neq q$, then $ff'\in(x_1,\dots,x_d)^{g(p)+g(q)+2}R[1/pq][[x_1,\dots,x_d]]$, and each coefficient of $ff'$ has the form $c/(p^iq^j)$ with $c\in R$ and $i,j\geq0$.
Choose $\alpha,\beta\in\Z$ with $\alpha p^i+\beta q^j=1$.
Then
\[
\frac{c}{p^iq^j}=\frac{c\alpha}{q^j}+\frac{c\beta}{p^i}.
\]
Applying this to each coefficient gives $ff'=u+v$, where
\begin{gather*}
u\in(x_1,\dots,x_d)^{g(p)+g(q)+2}R[1/p][[x_1,\dots,x_d]],\\
v\in(x_1,\dots,x_d)^{g(p)+g(q)+2}R[1/q][[x_1,\dots,x_d]].
\end{gather*}
Since $g(p)+1\geq1$ and $g(q)+1\geq1$, we get $u\in\mathcal{I}\bigl(\mathcal{I}^{g(p)+1}\hat{\otimes}\Z[1/p]\bigr)$ and $v\in\mathcal{I}\bigl(\mathcal{I}^{g(q)+1}\hat{\otimes}\Z[1/q]\bigr)$, so $ff'\in\mathcal{I}H^g$.
\end{proof}
By \Cref{claim:Hg-square}, we have $(H^g)^2\subset\mathcal{I}H^g\subset H^g$, so $R[[x_1,\dots,x_d]]+H^g$ is closed under multiplication.
Hence
\[
B^g=R[[x_1,\dots,x_d]]+H^g.
\]
Since $H^g$ is an $R[[x_1,\dots,x_d]]$-submodule of $R[[x_1,\dots,x_d]]\hat{\otimes}\Q$ and $(H^g)^2\subset H^g$, we have $B^gJ_g\subset J_g$, i.e., $J_g$ is an ideal of $B^g$.
Moreover, \eqref{eq:Jg} gives
\begin{equation}
\label{eq:Jg-cap}
J_g=J^g\cap B^g.
\end{equation}
\textup{(\ref{item:literal-to-degree})}.
The rings $R_n$ of \eqref{eq:Rn} increase with $n$, and $R_iR_j\subset R_{i+j}$.
Thus the right hand side of \eqref{eq:g-series} is a subring of $R[[x_1,\dots,x_d]]\hat{\otimes}\Q$, and every generator of $B^g$ belongs to it.
By \eqref{eq:Jg-cap} and \eqref{eq:Jg}, $J_g\subset\mathcal{I}\cdot\bigl(R[[x_1,\dots,x_d]]\hat{\otimes}\Q\bigr)$, so each element of $J_g$ has zero constant term.
Hence truncation modulo $\mathcal{I}^{n+1}$ kills $J_g^{n+1}$ and extends to $\hat{B}^g$.
The degree-$n$ coefficients of the image are degree-$n$ coefficients of an element of $B^g$, so they belong to $R_n$.
This proves \textup{(\ref{item:literal-to-degree})} by \eqref{eq:g-series}.

\textup{(\ref{item:degree-to-literal})}.
Write
\[
f=\sum_{n\geq0}f_n\in R[[x_1,\dots,x_d]]^g,
\]
with $f_n\in(R\otimes\Q)[x_1,\dots,x_d]$ a homogeneous polynomial of total degree $n$.
Its constant term is in $R$.
For $n\geq1$, let
\[
P^g_n:=\bigl\{p\text{ prime}\mid f_n\notin(R\otimes\Z_{(p)})[x_1,\dots,x_d]\bigr\}.
\]
This set is finite by \eqref{eq:g-series}.
By \eqref{eq:g-series}, every $p\in P^g_n$ satisfies $g(p)+1\leq n$.
Since $f_n$ is a polynomial, there is an integer $N\geq1$ with
\[
Nf_n\in R[x_1,\dots,x_d].
\]
If $N$ and $N'$ are two such integers, then so is $\gcd(N,N')$.
Let $N_n$ be the smallest such integer.
A prime $p$ divides $N_n$ if and only if $p\in P^g_n$, since $p\notin P^g_n$ means that some such integer is prime to $p$.
Hence
\[
N_n=\prod_{p\in P^g_n}p^{e_p}\quad\text{with }e_p\geq1.
\]
If $P^g_n=\emptyset$, set $f_{n,0}:=f_n$.
Otherwise, the integers $N_n/p^{e_p}$ for $p\in P^g_n$ have no common prime factor, so there are $a_p\in\Z$ with
\[
\sum_{p\in P^g_n}a_p\frac{N_n}{p^{e_p}}=1,
\]
and we set $f_{n,0}:=0$ and $f_{n,p}:=a_pN_nf_n/p^{e_p}$.
In both cases, we get a finite decomposition
\[
f_n=f_{n,0}+\sum_{p\in P^g_n}f_{n,p},
\]
where
\begin{itemize}
\item $f_{n,0}\in R[x_1,\ldots,x_d]$ is homogeneous of degree $n$, and
\item every $f_{n,p}\in R[1/p][x_1,\ldots,x_d]$ is homogeneous of degree $n$.
\end{itemize}

Hence, for $p\in P^g_n$, we have $h(p)=\lfloor g(p)/2\rfloor\leq(n-1)/2$.
Every monomial $x^\alpha$ with $|\alpha|=n$ factors as $x^\alpha=x^\beta x^\gamma$ with $|\beta|=n-h(p)-1$ and $|\gamma|=h(p)+1$.
Choosing such a factorization for each monomial of $f_{n,p}$ and grouping by $\beta$ gives
\[
f_{n,p}=\sum_{|\beta|=n-h(p)-1}x^\beta g_\beta,
\]
where each $g_\beta\in R[1/p][x_1,\ldots,x_d]$ is homogeneous of degree $h(p)+1$.
Since $x^\beta\in\mathcal{I}^{n-h(p)-1}$ and $g_\beta\in\mathcal{I}^{h(p)+1}\hat{\otimes}\Z[1/p]$, we get
\[
f_{n,p}\in
\mathcal{I}^{n-h(p)-1}\bigl(\mathcal{I}^{h(p)+1}\hat{\otimes}\Z[1/p]\bigr)
\subset\mathcal{I}^{n-h(p)-1}H^h
\subset J_h^{n-h(p)}
\subset J_h^{\lceil n/2\rceil}.
\]
Also $f_{n,0}\in J_h^n$.
In particular $f_n\in B^h$ and $f_n\in J_h^{\lceil n/2\rceil}$.
Therefore the partial sums $\sum_{n=0}^Nf_n$ are $J_h$-adically Cauchy.
Their limit in $\hat{B}^h$ maps to $f$, proving \textup{(\ref{item:degree-to-literal})}.

\end{proof}

\begin{definition}[Integral and algebraic formal sections]
\label{def:integrality-notions}
Let $\tilde v_K\colon\Spf\bigl(K[[x_1,\dots,x_d]]\bigr)\to E_K$ be the formal horizontal section through $v_K$ of \Cref{def:leaf-notions}.
Consider its composite
\[
\Spf\bigl(K[[x_1,\dots,x_d]]\bigr)\xrightarrow{\tilde v_K}E_K\to E
\]
with the projection $E_K\to E$.
\begin{enumerate}
\item
The section $\tilde v_K$ is \emph{integral} if this composite factors through $\Spf\bigl(R[1/N][[x_1,\dots,x_d]]\bigr)$ for some integer $N\geq 1$.
\item
The section $\tilde v_K$ is \emph{$g$-integral} if this composite factors through $\Spf\bigl(R[[x_1,\dots,x_d]]^g\bigr)$.
\item\label{item:omega-integral}
The section $\tilde v_K$ is \emph{$\omega(p)$-integral} if it is $g$-integral for some $g$ with
\[
\lim_{p\to\infty}\frac{g(p)}{p}=\infty,
\]
where $p$ runs through the primes.
\item
The section $\tilde v_K$ is \emph{algebraic} when, for an embedding $K\subset\C$, the image of its base change $\tilde v_{\C}$ is the completion at $v_{\C}$ of a closed subscheme of $E_{\C}$ of dimension $d$.
\end{enumerate}
\end{definition}

\begin{remark}
\label{rmk:omega-comparison}
The section $\tilde v_K$ is $\omega(p)$-integral in the sense of \Cref{def:integrality-notions}\ref{item:omega-integral} if and only if it is $\omega(p)$-integral in the completed-ring convention of \cite[Definition 2.2.5]{LL25}.
\end{remark}

\begin{proof}
If $\tilde v_K$ is $\omega(p)$-integral in the completed-ring convention of \cite[Definition 2.2.5]{LL25}, i.e., the composite factors through $\Spf(\hat B^g)$ for some $g$ with $g(p)/p\to\infty$, then composing with $\rho^g$ shows that it factors through $\Spf\bigl(R[[x_1,\dots,x_d]]^g\bigr)$ by \Cref{lem:omega-comparison}\textup{(\ref{item:literal-to-degree})}, so $\tilde v_K$ is $\omega(p)$-integral in the sense of \Cref{def:integrality-notions}\ref{item:omega-integral}.

Now we prove the reverse direction.
Suppose that $\tilde v_K$ is $g$-integral for some $g$ with $g(p)/p\to\infty$, and define a function $h$ by $h(p):=\lfloor g(p)/2\rfloor$.
For \Cref{lem:omega-comparison}\textup{(\ref{item:degree-to-literal})}, we lift $\tilde v_K^*$ to $\hat{B}^h$ (here by $B^h$ we mean \eqref{eq:Bg}) as follows.
The basis $e_1,\dots,e_r$ gives $E|_{U_s}\simeq U_s\times\mathbb{A}^r$ with fibre coordinates $t_1,\dots,t_r$, so $\cO(E|_{U_s})=\cO(U_s)[t_1,\dots,t_r]$.
Choose $b_i\in\hat{B}^h$ with $\rho^h(b_i)=\tilde v_K^*t_i$ (see \eqref{eq:rho-g} for $\rho^h$), and define
\[
\begin{aligned}
\lambda\colon\cO(U_s)[t_1,\dots,t_r]&\longrightarrow\hat{B}^h,\\
a&\longmapsto a\quad\text{via }\cO(U_s)\to R[[x_1,\dots,x_d]]\subset B^h\to\hat{B}^h\quad(a\in\cO(U_s)),\\
t_i&\longmapsto b_i.
\end{aligned}
\]
Since $t_1,\dots,t_r$ are free variables, $\lambda$ is a ring homomorphism, and $\rho^h\circ\lambda$ is $f\mapsto\tilde v_K^*(f\otimes1)$ on $\cO(E|_{U_s})$.
Since $g(p)/p\to\infty$, also $h(p)/p\to\infty$.
\end{proof}

\begin{remark}
\label{rmk:leaf-implications}
For $K$ a number field, the formal section $\tilde v_K$ satisfies
\[
\text{algebraic}\ \xRightarrow{\ \text{\Cref{lem:alg-implies-int}}\ }\ \text{integral}\ \xRightarrow{\ \text{\Cref{lem:int-implies-omega}}\ }\ \omega(p)\text{-integral}.
\]
The converse implications are \cite[Conjecture 2.3.1]{LL25}.
\Cref{lem:alg-implies-int,lem:int-implies-omega} are already mentioned in \cite[\S1.1]{LL25} without proper proofs; we include them here for completeness.
\end{remark}

\begin{lemma}[Algebraic implies integral]
\label{lem:alg-implies-int}
Assume that $K$ is a number field, and let $\tilde v_K\colon\Spf\bigl(K[[x_1,\dots,x_d]]\bigr)\to E_K$ be a formal section as in \Cref{def:leaf-notions}.
If $\tilde v_K$ is algebraic, then it is integral.
\end{lemma}

This is the analogue of Eisenstein's theorem \cite[\S1.1]{LL25} (citing \cite{Eis52}).

\begin{proof}
Since $\tilde v_K$ is algebraic, there is a closed subscheme $W\subset E_{\C}$ of dimension $d$ with $\hat W_{v_{\C}}=\operatorname{Im}\tilde v_{\C}$.
This equality forces $W$ to be smooth at $v_{\C}$: since $\tilde v_{\C}$ is a section of $\tau$, it is a closed immersion, so $\hat\cO_{W,v_{\C}}\simeq\C[[x_1,\dots,x_d]]$.
It remains to find such a subscheme over $K$.
Let $Z\subset E_K$ be the closed subscheme cut out, near $v_K$, by $I:=\ker\bigl(\cO_{E_K}\to K[[x_1,\dots,x_d]]\bigr)$.
Since $K\to\C$ is flat and $K[[x_1,\dots,x_d]]\otimes_K\C\to\C[[x_1,\dots,x_d]]$ is injective, $I\otimes_K\C=\ker\bigl(\cO_{E_{\C}}\to\C[[x_1,\dots,x_d]]\bigr)$, which contains the ideal of $W$.
Hence $\operatorname{Im}\tilde v_{\C}\subset\hat Z_{\C,v_{\C}}\subset\hat W_{v_{\C}}=\operatorname{Im}\tilde v_{\C}$, so $Z_{\C}$ is smooth of dimension $d$ at $v_{\C}$, and therefore $Z$ is smooth of dimension $d$ at $v_K$.

Let $\overline{Z}$ be the closure of $Z$ in $E$.
The locus where $\overline{Z}\to\Spec R$ is smooth of relative dimension $d$ is open and contains $v_K$.
Its preimage under $v\colon\Spec R\to E$ is open and contains the generic point, so it contains all but finitely many closed points of $\Spec R$, because $\dim R=1$.
Invert an integer $N$ divisible by the residue characteristics of those finitely many points, and enlarge $N$ so that the rank $d$ module $v^*\Omega^1_{\overline{Z}/R}$ becomes free over $R[1/N]$.
Then the completion of $\overline{Z}_{R[1/N]}$ along $v$ is $\Spf\bigl(R[1/N][[x_1,\dots,x_d]]\bigr)$, and its generic fiber is the image of $\tilde v_K$, so $\tilde v_K$ is integral.
\end{proof}

\begin{lemma}[Integral implies $\omega(p)$-integral]
\label{lem:int-implies-omega}
Let $\tilde v_K$ be as in \Cref{def:leaf-notions}.
If $\tilde v_K$ is integral, then it is $\omega(p)$-integral \cite[\S1.1]{LL25}.
\end{lemma}

\begin{proof}
Suppose $\tilde v_K$ factors through $\Spf\bigl(R[1/N][[x_1,\dots,x_d]]\bigr)$.
Define
\[
g(p):=\begin{cases}0 & \text{if } p\mid N,\\ p^2 & \text{if } p\nmid N.\end{cases}
\]
Then $g(p)/p\to\infty$.
We show that $\tilde v_K$ is $g$-integral.
Let $V=\Spec A\subset E$ be an affine open subscheme containing the image of $v$.
The morphism $v\colon\Spec R\to V$ corresponds to a ring map $v^*\colon A\to R$.
The formal section $\tilde v_K$ restricts to a map $\Spf K[[x_1,\dots,x_d]]\to V_K$, which corresponds to a ring map $\tilde v_K^*\colon A\otimes_RK\to K[[x_1,\dots,x_d]]$.
Let $f\in A$.
By integrality of $\tilde v_K$,
\[
\tilde v_K^*(f\otimes1)\in R[1/N][[x_1,\dots,x_d]].
\]
By \eqref{eq:g-series}, we need to check the coefficients of $\tilde v_K^*(f\otimes1)$ for every prime $p$:
\begin{enumerate}[label=\alph*)]
\item if $p\mid N$, then $g(p)=0$, so we need the constant term to lie in $R$;
\item if $p\nmid N$, then we need that no coefficient has a denominator divisible by $p$.
\end{enumerate}
Case b) holds because the coefficients lie in $R[1/N]$.
For a), since $\tilde v_K$ passes through $v_K$, composing $\tilde v_K^*$ with $K[[x_1,\dots,x_d]]\to K$, $x_i\mapsto0$, gives $v_K^*=v^*\otimes_RK\colon A\otimes_RK\to K$:
\begin{equation}
\label{eq:const-term}
\begin{tikzcd}[ampersand replacement=\&]
A \arrow[r] \arrow[d, "v^*"'] \& A\otimes_RK \arrow[r, equal] \arrow[d, "v_K^*"] \& A\otimes_RK \arrow[d, "\tilde v_K^*"] \\
R \arrow[r] \& K \& K[[x_1,\dots,x_d]] \arrow[l, "x_i\mapsto0"']
\end{tikzcd}
\end{equation}
Hence setting $x_1=\dots=x_d=0$ in $\tilde v_K^*(f\otimes1)$ gives $v^*f\in R$.
This proves a).
\end{proof}

\begin{lemma}[Base change of $g$-integrality]
\label{lem:g-integral-base-change}
Let $K$ be a number field, let $R\subset K$ be a finitely generated $\Z$-algebra, and let $L/K$ be a finite extension.
Let $R'$ be a localization of the integral closure of $R$ in $L$.
Let $S$, $\cE$, $\nabla$, $s$, $v$ and $\tilde v_K$ be as in \Cref{def:leaf-notions}, and write $\tilde v_L$ for the formal horizontal section of $(\cE,\nabla)_L$ through $v_L$.
If $\tilde v_K$ is $g$-integral, then $\tilde v_L$ is $g$-integral for the same $g$.
\end{lemma}

\begin{proof}
Since $L/K$ is finite, $K[[x_1,\dots,x_d]]\otimes_KL=L[[x_1,\dots,x_d]]$, so $\tilde v_K\otimes_KL$ is a formal section with $\nabla(\tilde v_K\otimes_KL)=0$, and it passes through $v_L$ by $\eqref{eq:const-term}\otimes_KL$.
By the uniqueness in \Cref{def:leaf-notions} it equals $\tilde v_L$.
Let $V=\Spec A\subset E$ be an affine open subscheme containing the image of $v$, and let $f\in A$.
Then $\tilde v_L^*(f\otimes1)=\tilde v_K^*(f\otimes1)$.
By \eqref{eq:g-series}, its coefficients of total degree $n$ lie in $R_n\subset R'_n$, so by \eqref{eq:g-series} for $R'$ it lies in $R'[[x_1,\dots,x_d]]^g$.
Since the elements $f\otimes1$ generate $A\otimes_RR'$ as an $R'$-module, $\tilde v_L$ is $g$-integral.
\end{proof}

\begin{para}
\label{para:reduction}
Let $\nabla$ be an integrable connection on the vector bundle $\cE$, let $s\in S(R)$, and let $v\colon\Spec R\to E$ be a morphism with $\tau\circ v=s$.
Write $\tilde{v}_K$ for the formal horizontal section of $(\cE,\nabla)_K$ through $v_K$ (\Cref{def:leaf-notions}), and
\begin{equation}\label{eq:Is}
\cI_s:=\mathrm{Ker}\bigl(\cO_S\to s_*\cO_{\Spec R}\bigr)
\end{equation}
for the ideal sheaf of $s$.
Suppose that $\tilde{v}_K$ is $g$-integral for some $g\colon\{\mathrm{primes}\}\to\Z_{\geq 0}$.
Let $\mathfrak{p}$ be a closed point of $\mathrm{Spec}(R)$ with residue characteristic $p$.
Write
\[
S_{k(\mathfrak p)}:=S\otimes_Rk(\mathfrak p)\quad\text{and}\quad(\cE_{k(\mathfrak p)},\nabla)=(\cE\otimes_Rk(\mathfrak p),\nabla)
\]
for the fiber of $S$ over $\mathfrak p$ and the restriction of $(\cE,\nabla)$ to it, respectively.
Let $n\geq1$.
Since $s$ is a section of the smooth morphism $S\to\Spec R$, the sheaf $\cO_S/\cI_s^n$ is an extension of the finite locally free sheaves $\cI_s^k/\cI_s^{k+1}\cong\operatorname{Sym}^k(s^*\Omega^1_{S/R})$ on $\Spec R$, $0\leq k<n$.
Hence $\cO_S/\cI_s^n$ is finite locally free over $R$.
Tensoring $0\to\cI_s^n\to\cO_S\to\cO_S/\cI_s^n\to0$ with $k(\mathfrak p)$ therefore shows that $\cI_s^n\otimes_Rk(\mathfrak p)\to\cO_{S_{k(\mathfrak p)}}$ is injective.
Write
\[\cI_{s,k(\mathfrak p)}:=\cI_s\otimes_Rk(\mathfrak p)\subset\cO_{S_{k(\mathfrak p)}}.\]
Then $\cI_{s,k(\mathfrak p)}^n=\cI_s^n\otimes_Rk(\mathfrak p)$.
Set $\cI_s^n:=\cO_S$ and $\cI_{s,k(\mathfrak p)}^n:=\cO_{S_{k(\mathfrak p)}}$ for $n\leq0$.
For $n\leq g(p)$, every prime $q$ with $n\geq g(q)+1$ differs from $p$, so\[R_n\subset R_{\mathfrak p}\]($R_n$ is defined in \eqref{eq:Rn}).

Write $s_{k(\mathfrak p)}\in S_{k(\mathfrak p)}(k(\mathfrak p))$ for the base change of $s$ along $R\to k(\mathfrak p)$.
Let $U'\subset S$ be an affine open containing $s_{k(\mathfrak p)}$ and contained in the neighborhood $U_s$ of \Cref{def:leaf-notions}.
Since $s$ is a regular immersion \cite[\href{https://stacks.math.columbia.edu/tag/067R}{Tag 067R}]{SP}, $\bigl(\cO_S/\cI_s^n\bigr)(U')=R[x_1,\dots,x_d]/(x_1,\dots,x_d)^n$, and \[\Gamma\bigl(\hat S_{s_K},\cE_K\bigr)=\bigoplus_iK[[x_1,\dots,x_d]]\,e_i.\] Write \begin{equation}\label{eq:vK}\tilde v_K=\sum_if_ie_i\quad\text{with }f_i\in R[[x_1,\dots,x_d]]^g.\end{equation}
By \eqref{eq:const-term}, the constant term of $\tilde v_K$ is
\begin{equation}\label{eq:vK-const}\sum_if_i(0)\,e_i(s_K)=v_K\in s_K^*\cE_K=\bigoplus_iK\,e_i(s_K).\end{equation}
 By \eqref{eq:g-series}, the coefficients of total degree $n<g(p)$ of the $f_i$ lie in $R_n\subset R_{\mathfrak p}$.
Let\[f_i^{<g(p)}\in R_{\mathfrak p}[x_1,\dots,x_d]\]be $f_i$ with its terms of total degree $\geq g(p)$ removed, and set \begin{equation}\label{eq:vRp}\tilde v_{R_{\mathfrak p}}:=\sum_if_i^{<g(p)}e_i\in\cE(U')\otimes_RR_{\mathfrak p}.\end{equation}
If $g(p)\ge1$, then $f_i^{<g(p)}$ keeps the constant term $f_i(0)\in R_0=R$.
As $s^*\cE$ is locally free, $s^*\cE\otimes_RR_{\mathfrak p}\to s^*\cE\otimes_RK$ is injective.
So by \eqref{eq:vK-const} the constant term of $\tilde v_{R_{\mathfrak p}}$ is
\begin{equation}\label{eq:vRp-const}\sum_if_i(0)\,e_i(s)=v\in s^*\cE\otimes_RR_{\mathfrak p}.\end{equation}
 Since $U'$ contains $s(\Spec R_{\mathfrak p})$, the support of $\cO_S/\cI_s^{g(p)}\otimes_RR_{\mathfrak p}$, we have
\[\bigl(\cE/\cI_s^{g(p)}\cE\bigr)(U')\otimes_RR_{\mathfrak p}=\Gamma\bigl(S,\cE/\cI_s^{g(p)}\cE\otimes_RR_{\mathfrak p}\bigr).\]
So the image of $\tilde v_{R_{\mathfrak p}}$ lies in $\Gamma\bigl(S,\cE/\cI_s^{g(p)}\cE\otimes_RR_{\mathfrak p}\bigr)$.
Write $\tilde v_{k(\mathfrak p)}\in\Gamma\bigl(S,\cE/\cI_s^{g(p)}\cE\otimes_Rk(\mathfrak p)\bigr)$ for its reduction modulo $\mathfrak p$.
Write $U'_{k(\mathfrak p)}:=U'\otimes_Rk(\mathfrak p)$, and let
\begin{equation}\label{eq:vRp-bar}\overline{\tilde v}_{R_{\mathfrak p}}\in\Gamma(U'_{k(\mathfrak p)},\cE_{k(\mathfrak p)})\end{equation}
be the image of $\tilde v_{R_{\mathfrak p}}$ under $\cE(U')\otimes_RR_{\mathfrak p}\to\cE(U')\otimes_Rk(\mathfrak p)$.
It is a lift of $\tilde v_{k(\mathfrak p)}|_{U'_{k(\mathfrak p)}}$.
If $g(p)\ge1$, then reducing \eqref{eq:vRp-const} modulo $\mathfrak p$ gives
\begin{equation}\label{eq:vRp-bar-const}\overline{\tilde v}_{R_{\mathfrak p}}(s_{k(\mathfrak p)})=\bar v\in s^*\cE\otimes_Rk(\mathfrak p),\end{equation}
where $\bar v$ is the image of $v$.
See \Cref{fig:reduction}.
\begin{center}\begin{minipage}{\textwidth}
\centering
\newcommand{\fu}[4]{\shade[top color=#1,bottom color=white] ($(#2)+(-#3,-#4)$) rectangle ($(#2)+(#3,0)$);}
\newcommand{\fv}[4]{\shade[bottom color=#1,top color=white] ($(#2)+(-#3,0)$) rectangle ($(#2)+(#3,#4)$);}
\makebox[\linewidth]{\resizebox{\textwidth}{!}{%
\begin{tikzpicture}[>=Stealth,line cap=round,line join=round]
\def\dx{1.6}\def\dy{1.8}\def\H{3.4}
\draw[thick] (0,0) -- (8,0) node[right]{$\Spec R$};
\foreach \x in {1,5}{\fill (\x,0) circle (1.3pt);}
\fill (3,0) circle (1.8pt) node[below]{\small $\mathfrak p$};
\draw[thick,fill=white] (6.6,0) ellipse [x radius=5pt, y radius=2.6pt] node[below=2pt]{\small $\eta=\Spec K$};
\draw[orange!80!black,decorate,decoration={brace,mirror,amplitude=4pt}] (2.8,-0.55) -- (6.9,-0.55) node[midway,below=4pt]{\small $\Spec R_{\mathfrak p}=\{\mathfrak p,\eta\}$};
\begin{scope}[yshift=1.4cm]
\fill[gray!10] (0,0) -- (8,0) -- (8+\dx,\dy) -- (\dx,\dy) -- cycle;
\begin{scope}
\clip (0,0) -- (8,0) -- (8+\dx,\dy) -- (\dx,\dy) -- cycle;
\def\Up{(3.36,0.405) .. controls (4.5,-0.2) and (7.5,-0.2) .. (8.6,0.6) .. controls (9.4,1.2) and (8.8,1.75) .. (7.2,1.6) .. controls (5.5,1.45) and (4.8,1.6) .. (4.2304,1.3842) .. controls (3.3,1.05) and (2.85,0.75) .. (3.36,0.405) -- cycle}
\fill[orange!15] \Up; \draw[orange!80!black,dashed] \Up;
\foreach \w in {0.04,0.08,...,0.36}{\fill[red!80!black,opacity=0.07] (6.6-\w,0) -- (6.6+\w,0) -- (6.6+\dx+\w,\dy) -- (6.6+\dx-\w,\dy) -- cycle;}
\end{scope}
\draw[thick] (0,0) -- (8,0) -- (8+\dx,\dy) -- (\dx,\dy) -- cycle;
\foreach \x in {1,5}{\draw[gray!60] (\x,0) -- (\x+\dx,\dy);}
\draw[gray!70!black,thick] (3,0) -- (3+\dx,\dy);
\node[gray!60!black] at (3,-0.05) [below] {\scriptsize $S_{k(\mathfrak p)}$};
\node[red!70!black] at (6.6,-0.05) [below]{\scriptsize $S_K$};
\node[orange!80!black] at (8.7,0.35) {\small $U'$};
\def\sc{(0.2,0.35) .. controls (1.5,0.9) and (2.3,0.65) .. (3+0.45*\dx,0.45*\dy) .. controls (4.6,1.0) and (5.6,0.6) .. (6.6+0.45*\dx,0.45*\dy) .. controls (7.8,1.0) .. (8.6,0.95)}
\draw[blue!60!black,line width=5pt,opacity=0.18] \sc;
\draw[blue!60!black,very thick] \sc node[above right]{$s$};
\fill[blue!60!black,opacity=0.35] (3+0.45*\dx,0.45*\dy) ellipse [x radius=0.22, y radius=0.34, rotate=-41];
\fill[blue!60!black] (3+0.45*\dx,0.45*\dy) circle (1.8pt) node[below right=1pt]{\scriptsize $s_{k(\mathfrak p)}$};
\draw[red!70!black,thick,dotted] (6.6+0.45*\dx,0.45*\dy) circle (0.33);
\fill[blue!60!black] (6.6+0.45*\dx,0.45*\dy) circle (2pt) node[below right=2pt]{\scriptsize $s_K$};
\node at (0.3,1.25) {$S$};
\end{scope}
\draw[->] (4,1.3) -- (4,0.25);
\begin{scope}[yshift=5.6cm]
\def\D{1.0}
\fill[gray!12,path fading=south] (0,-\D) -- (8,-\D) -- (8+\dx,\dy-\D) -- (8+\dx,\dy) -- (8,0) -- (0,0) -- cycle;
\fill[gray!12,path fading=north] (0,0) -- (8,0) -- (8+\dx,\dy) -- (8+\dx,\dy+\H) -- (\dx,\dy+\H) -- (0,\H) -- cycle;
\fv{gray!60}{\dx,\dy}{0.2pt}{\H} \fu{gray!60}{\dx,\dy}{0.2pt}{\D}
\fill[gray,opacity=0.3,path fading=north] (3,0) -- (3+\dx,\dy) -- (3+\dx,\dy+\H) -- (3,\H) -- cycle;
\fill[gray,opacity=0.3,path fading=south] (3,-\D) -- (3+\dx,\dy-\D) -- (3+\dx,\dy) -- (3,0) -- cycle;
\fv{gray!70!black}{3,0}{0.25pt}{\H} \fv{gray!70!black}{3+\dx,\dy}{0.25pt}{\H} \fu{gray!70!black}{3,0}{0.25pt}{\D} \fu{gray!70!black}{3+\dx,\dy}{0.25pt}{\D}
\foreach \w in {0.04,0.08,...,0.36}{\fill[red!80!black,opacity=0.05,path fading=north] (6.6-\w,0) -- (6.6+\w,0) -- (6.6+\dx+\w,\dy) -- (6.6+\dx+\w,\dy+\H) -- (6.6+\dx-\w,\dy+\H) -- (6.6-\w,\H) -- cycle; \fill[red!80!black,opacity=0.05,path fading=south] (6.6-\w,-\D) -- (6.6+\w,-\D) -- (6.6+\dx+\w,\dy-\D) -- (6.6+\dx+\w,\dy) -- (6.6+\dx-\w,\dy) -- (6.6-\w,0) -- cycle;}
\draw[orange!80!black,very thick] (3.360,1.566) .. controls (3.536,1.879) and (3.633,1.884) .. (3.720,2.060) .. controls (3.851,2.290) and (4.036,2.311) .. (4.231,2.543);
\fill[orange!60,opacity=0.55] (6.6-0.3+0.12*\dx,0.12*\dy+1.25) .. controls (6.6+0.3*\dx,0.3*\dy+1.55) .. (6.6+0.3+0.85*\dx,0.85*\dy+1.35) -- (6.6+0.3+0.85*\dx,0.85*\dy+1.15) .. controls (6.6+0.3*\dx,0.3*\dy+1.35) .. (6.6-0.3+0.12*\dx,0.12*\dy+1.05) -- cycle;
\def\vc{(0.2,0.35+1.0) .. controls (1.5,0.9+1.5) and (2.3,0.65+0.9) .. (3+0.45*\dx,0.45*\dy+1.25) .. controls (4.6,1.0+1.6) and (5.6,0.6+0.9) .. (6.6+0.45*\dx,0.45*\dy+1.3) .. controls (7.8,1.0+1.5) .. (8.6,0.95+1.2)}
\draw[blue!60!black,very thick] \vc node[right]{$v$};
\draw[blue!80!black,line width=3.2pt,opacity=0.75] (3+0.37*\dx,0.37*\dy+1.17) .. controls (3+0.45*\dx,0.45*\dy+1.3) .. (3+0.53*\dx,0.53*\dy+1.3);
\fill[blue!60!black] (3+0.45*\dx,0.45*\dy+1.25) circle (1.8pt);
\fill[red!70!black,opacity=0.55] (6.6+0.45*\dx,0.45*\dy+1.3) ellipse [x radius=0.36, y radius=0.16, rotate=12];
\fill[blue!60!black] (6.6+0.45*\dx,0.45*\dy+1.3) circle (2pt) node[below right=2pt,blue!60!black]{\scriptsize $v_K$};
\fv{black}{0,0}{0.4pt}{\H} \fv{black}{8,0}{0.4pt}{\H} \fv{black}{8+\dx,\dy}{0.4pt}{\H} \fu{black}{0,0}{0.4pt}{\D} \fu{black}{8,0}{0.4pt}{\D} \fu{black}{8+\dx,\dy}{0.4pt}{\D}
\draw[thick] (0,0) -- (8,0) -- (8+\dx,\dy); \draw[thin,dashed,gray!60!black] (0,0) -- (\dx,\dy) -- (8+\dx,\dy);
\node at (0.35,1.95) {$E$};
\node[gray!60!black] at (3+\dx,\dy+2.6) [right] {\small $E\otimes_Rk(\mathfrak p)$};
\node[red!70!black] at (6.6+\dx+0.3,\dy+2.6) [right] {\small $E_K$};
\end{scope}
\draw[->] (5.6,4.5) -- node[right]{$\tau$} (5.6,3.35);
\begin{scope}[font=\small]
\node[red!70!black,align=left,anchor=west] (b) at (10.4,11.0) {$\tilde v_K=\sum_if_ie_i\in\Gamma(\hat S_{s_K},\cE_K)$};
\draw[->,red!70!black] (b.west) .. controls (7.50,11.15) and (7.45,9.50) .. (7.45,7.85);
\node[orange!80!black,align=left,anchor=west] (c) at (10.4,9.3) {$\tilde v_{R_{\mathfrak p}}=\sum_if_i^{<g(p)}e_i\in\cE(U')\otimes_RR_{\mathfrak p}$};
\draw[->,orange!80!black] (c.west) .. controls (9.0,9.4) and (7.812,9.0) .. (7.812,8.071);
\draw[->,orange!80!black] (c.west) .. controls (9.65,8.70) and (8.80,7.35) .. (4.077,7.993);
\node[blue!60!black,align=left,anchor=west] (a) at (10.4,8.15) {$\tilde v_{k(\mathfrak p)}\in\Gamma\bigl(S,\cE/\cI_s^{g(p)}\cE\otimes_Rk(\mathfrak p)\bigr)$,\\ image of $\tilde v_{R_{\mathfrak p}}$ mod $(\cI_s^{g(p)},\mathfrak p)$};
\draw[->,blue!60!black] (a.west) .. controls (9.6,6.70) and (5.3,6.60) .. (3.83,7.83);
\node[orange!80!black,align=left,anchor=west] (c2) at (10.4,7.0) {$\overline{\tilde v}_{R_{\mathfrak p}}\in\Gamma(U'_{k(\mathfrak p)},\cE_{k(\mathfrak p)})$,\\ image of $\tilde v_{R_{\mathfrak p}}$ mod $\mathfrak p$};
\draw[->,orange!80!black] (c2.west) .. controls (9.0,5.95) and (4.9,5.95) .. (4.077,7.993);
\node[blue!60!black,align=left,anchor=west] (vbar) at (10.4,5.85) {$\bar v\in\Gamma(\Spec k(\mathfrak p),s_{k(\mathfrak p)}^*\cE_{k(\mathfrak p)})$,\\ image of $v$ mod $\mathfrak p$};
\draw[->,blue!60!black,shorten >=2.2pt] (vbar.west) .. controls (9.0,5.05) and (5.0,5.15) .. (3+0.45*\dx,5.6+0.45*\dy+1.25);
\node[blue!60!black,align=left,anchor=west] (d) at (10.4,1.2) {blue band: the support\\ $s(\Spec R)$ of $\cO_S/\cI_s^{g(p)}$};
\draw[->,blue!60!black] (d.west) .. controls (9.0,0.9) and (6.0,1.2) .. (5.2,2.2);
\node[blue!60!black,align=right,anchor=east] (g) at (-0.3,0.9) {blue disk around $s_{k(\mathfrak p)}$:\\ $\Spec\bigl(\cO_{S_{k(\mathfrak p)}}/\cI_{s,k(\mathfrak p)}^{g(p)}\bigr)$};
\draw[->,blue!60!black] (g.east) .. controls (2.2,0.55) and (3.5,0.6) .. (3.62,1.98);
\node[red!70!black,anchor=west] (hatS) at (10.4,2.45) {$\hat S_{s_K}=\Spf K[[x_1,\dots,x_d]]$};
\draw[->,red!70!black,shorten >=0.35cm] (hatS.west) -- (6.6+0.45*\dx,1.4+0.45*\dy);
\end{scope}
\path (487.24pt,0);
\end{tikzpicture}%
}}
\captionsetup{labelformat=simple,labelsep=period}

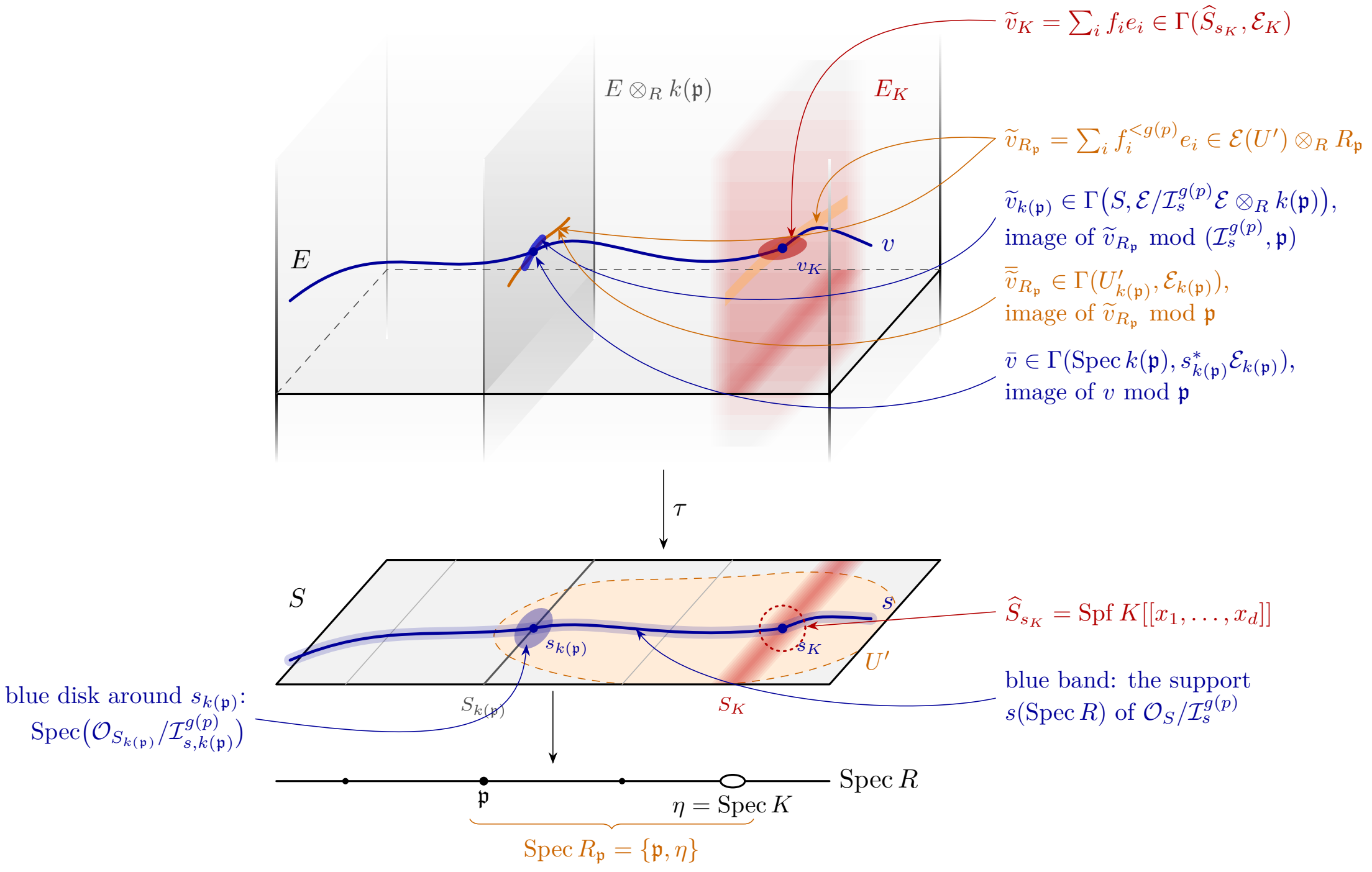
\captionof{figure}{The schemes and sections in \ref{para:reduction}.}
\label{fig:reduction}
\end{minipage}\end{center}

\end{para}

\begin{lemma}
\label{lem:omega-reduction}
In the situation of \ref{para:reduction}, for every open $U\subset S_{k(\mathfrak p)}$ containing $s_{k(\mathfrak p)}$, every lift $w\in\cE_{k(\mathfrak p)}(U)$ of $\tilde v_{k(\mathfrak p)}|_U\in\bigl(\cE/\cI_{s,k(\mathfrak p)}^{g(p)}\cE\bigr)(U)$, and every $\partial\in T_{S_{k(\mathfrak p)}}(U)$, we have
\[
\psi_p(\partial)(w)\in\cI_{s,k(\mathfrak p)}^{g(p)-p}\cE_{k(\mathfrak p)}(U).
\]
\end{lemma}

\begin{proof}
The claim is local, so with $U'$ as in \ref{para:reduction} and $\tilde v_{R_{\mathfrak p}}$ as in \eqref{eq:vRp} we may replace $U$ by $U\cap U'$.
By the Leibniz rule, for every $n$,
\begin{equation}\label{eq:leibniz-I}
\nabla\bigl(\cI_s^n\cE\bigr)\subset\cI_s^{n-1}\bigl(\Omega^1_{S/R}\otimes_{\cO_S}\cE\bigr),\qquad\nabla_\partial\bigl(\cI_s^n\cE\bigr)\subset\cI_s^{n-1}\cE,
\end{equation}
and likewise on $S_{k(\mathfrak p)}$ with $\cI_{s,k(\mathfrak p)}$ in place of $\cI_s$ and $\cE_{k(\mathfrak p)}$ in place of $\cE$.
By horizontality of $\tilde v_K$, $\nabla\tilde v_K=0$, so we have $\nabla\tilde v_{R_{\mathfrak p}}=-\nabla(\tilde v_K-\tilde v_{R_{\mathfrak p}})$.
As \[\tilde v_K-\tilde v_{R_{\mathfrak p}}=\sum_i(f_i-f_i^{<g(p)})e_i\quad\text{with }f_i-f_i^{<g(p)}\in(x_1,\dots,x_d)^{g(p)}K[[x_1,\dots,x_d]],\] we have $\tilde v_K-\tilde v_{R_{\mathfrak p}}\in\cI_s^{g(p)}\cE$, so \eqref{eq:leibniz-I} gives \[\nabla\tilde v_{R_{\mathfrak p}}\in\cI_s^{g(p)-1}\bigl(\Omega^1_{S/R}\otimes_{\cO_S}\cE\bigr).\] Hence $\nabla\overline{\tilde v}_{R_{\mathfrak p}}\in\cI_{s,k(\mathfrak p)}^{g(p)-1}\bigl(\Omega^1_{S_{k(\mathfrak p)}/k(\mathfrak p)}\otimes\cE_{k(\mathfrak p)}\bigr)$ for $\overline{\tilde v}_{R_{\mathfrak p}}$ from \eqref{eq:vRp-bar}.
As $w-\overline{\tilde v}_{R_{\mathfrak p}}\in\cI_{s,k(\mathfrak p)}^{g(p)}\cE_{k(\mathfrak p)}(U)$ (both are lifts of $\tilde v_{k(\mathfrak p)}|_U$), \eqref{eq:leibniz-I} on $S_{k(\mathfrak p)}$ gives \[\nabla\bigl(w-\overline{\tilde v}_{R_{\mathfrak p}}\bigr)\in\cI_{s,k(\mathfrak p)}^{g(p)-1}\bigl(\Omega^1_{S_{k(\mathfrak p)}/k(\mathfrak p)}\otimes\cE_{k(\mathfrak p)}\bigr)(U).\] So \[\nabla w=\nabla\overline{\tilde v}_{R_{\mathfrak p}}+\nabla\bigl(w-\overline{\tilde v}_{R_{\mathfrak p}}\bigr)\in\cI_{s,k(\mathfrak p)}^{g(p)-1}\bigl(\Omega^1_{S_{k(\mathfrak p)}/k(\mathfrak p)}\otimes\cE_{k(\mathfrak p)}\bigr)(U).\] The map $\nabla_\partial$ is the composition of $\nabla$ with $\partial\otimes\mathrm{id}_{\cE}\colon\Omega^1_{S_{k(\mathfrak p)}/k(\mathfrak p)}\otimes\cE_{k(\mathfrak p)}\to\cE_{k(\mathfrak p)}$, and the latter is $\cO$-linear (do not be misled by the symbol: $\partial\otimes\mathrm{id}_{\cE}$ evaluates forms on $\partial$ and does not differentiate), so it sends $\cI_{s,k(\mathfrak p)}^{g(p)-1}\bigl(\Omega^1_{S_{k(\mathfrak p)}/k(\mathfrak p)}\otimes\cE_{k(\mathfrak p)}\bigr)(U)$ into $\cI_{s,k(\mathfrak p)}^{g(p)-1}\cE_{k(\mathfrak p)}(U)$; the same holds for $\partial^p\in T_{S_{k(\mathfrak p)}}(U)$.
Hence \[\nabla_\partial(w),\nabla_{\partial^p}(w)\in\cI_{s,k(\mathfrak p)}^{g(p)-1}\cE_{k(\mathfrak p)}(U),\] so $\nabla_\partial^p(w)=\nabla_\partial^{p-1}\bigl(\nabla_\partial(w)\bigr)\in\cI_{s,k(\mathfrak p)}^{g(p)-p}\cE_{k(\mathfrak p)}(U)$, and hence
\[
\psi_p(\partial)(w)\in\cI_{s,k(\mathfrak p)}^{g(p)-p}\cE_{k(\mathfrak p)}(U).
\]
\end{proof}

\begin{lemma}[Algebraicity criterion]
\label{lem:cover-algebraic}
In the situation of \Cref{def:leaf-notions}, suppose there are the following.
\begin{enumerate}[label=\textup{(\alph*)}]
\item An \'{e}tale morphism $h\colon S'\to S_K$ with $S'$ connected and of finite type over $K$.
\item A closed point $s'\in S'$ over $s_K$.
\item \emph{Algebraicity up to \'{e}tale extension.} An element $u\in H^0(S',\cE|_{S'})=\Hom_{S_K}(S',E_K)$ with $\nabla u=0$ such that the diagram
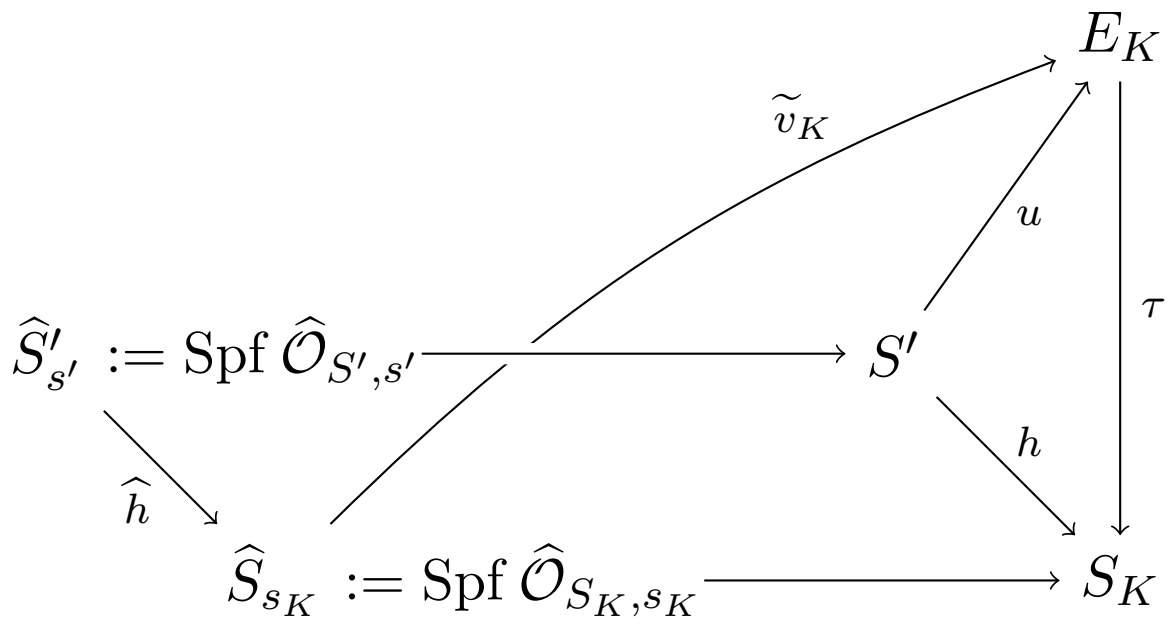
\begin{equation}\label{eq:cover-square}
\begin{tikzpicture}[>=to,baseline=(current bounding box.center)]
\node (A) at (0,1.5) {$\hat S'_{s'}$};
\node[anchor=west,inner sep=0pt] (A2) at (A.east) {$:=\Spf\hat\cO_{S',s'}$};
\node (P) at (5.6,1.5) {$S'$};
\node (B) at (1.5,0) {$\hat S_{s_K}$};
\node[anchor=west,inner sep=0pt] (B2) at (B.east) {$:=\Spf\hat\cO_{S_K,s_K}$};
\node (Q) at (7.1,0) {$S_K$};
\node (E) at (7.1,3.6) {$E_K$};
\draw[->] (B) to[bend left=12] node[pos=0.75,above left=-2pt]{$\scriptstyle\tilde v_K$} (E);
\draw[->,preaction={draw,white,line width=4pt,shorten <=2pt,shorten >=6pt}] (A2) -- (P);
\draw[->] (A) -- node[below left=-2pt]{$\scriptstyle\hat h$} (B);
\draw[->] (B2) -- (Q);
\draw[->] (P) -- node[above right=-2pt]{$\scriptstyle h$} (Q);
\draw[->] (P) -- node[below right=-2pt]{$\scriptstyle u$} (E);
\draw[->] (E) -- node[right]{$\scriptstyle\tau$} (Q);
\end{tikzpicture}
\end{equation}
commutes, where $\hat h$ is the morphism induced by $h$.
\end{enumerate}
Let $Z\subset E_K$ be the Zariski closure of $u(S')$, with its reduced structure.
Then the following hold.
\begin{enumerate}
\item\label{item:coveralg-Z}
$Z$ is the Zariski closure of the image of the morphism $\Spec\hat\cO_{S_K,s_K}\to E_K$ induced by $\tilde v_K$.
In particular, $Z$ does not depend on the choice of $S'$, $s'$, and $u$.
\item\label{item:coveralg-1}
$\dim Z=d$.
\item\label{item:coveralg-2}
The closed formal subschemes $\Spf\hat\cO_{Z,v_K}$ and $\tilde v_K\bigl(\Spf\hat\cO_{S_K,s_K}\bigr)$ of $\Spf\hat\cO_{E_K,v_K}$ coincide.
\item\label{item:coveralg-3}
The formal horizontal section $\tilde v_K$ is algebraic.
\end{enumerate}
\end{lemma}

\begin{proof}
Write $L:=k(s')$, a finite extension of $K$.
Let $Z_0\subset E_K$ be the Zariski closure of the image of $\Spec\hat\cO_{S_K,s_K}\to E_K$, with its reduced structure.
Since $h$ is \'etale at $s'$,
\[\hat{\cO}_{S',s'}\simeq\hat{\cO}_{S_K,s_K}\otimes_KL\simeq L[[x_1,\dots,x_d]],\]
which is finite free over $K[[x_1,\dots,x_d]]$.
So $\Spec L[[x_1,\dots,x_d]]\to\Spec K[[x_1,\dots,x_d]]$ is surjective.
By \eqref{eq:cover-square}, the composite
\[\Spec L[[x_1,\dots,x_d]]\to\Spec K[[x_1,\dots,x_d]]\to E_K\]
factors through $u\colon S'\to E_K$.
Hence $Z_0\subset Z$.
Since $S'$ is smooth over $K$ and connected, it is irreducible, so $Z$ is irreducible.
As $\tau\circ u=h$ is \'etale, $u$ is quasi-finite, so $\dim Z\le\dim S'=d$.
The composite $\cO_{S_K,s_K}\to\cO_{E_K,v_K}\to K[[x_1,\dots,x_d]]$ is injective, 
so $Z_0$ dominates the connected component $S_K^0$ of $S_K$ through $s_K$, and $\dim Z_0\geq\dim S_K^0=d$. 
As $Z_0\subset Z$ is closed, $Z$ is irreducible, and $\dim Z_0\ge d\ge\dim Z$, we get $Z_0=Z$ and $\dim Z=d$.
This proves \ref{item:coveralg-Z} and \ref{item:coveralg-1}.

Let $\cI_Z\subset\cO_{E_K}$ be the ideal sheaf of $Z$.
By \ref{item:coveralg-Z}, for an affine open $V=\Spec A\subset E_K$ containing $ v_K$,
\begin{equation}\label{eq:cov-ker}
\cI_Z(V)=\mathrm{Ker}\bigl(\tilde v_K^*\colon A\to K[[x_1,\dots,x_d]]\bigr).
\end{equation}
The reduced version of the two ideals are clearly equal.
Since the $\cI_Z(V)$ is reduced by definition, it suffices to show the right hand side is reduced as well. 
This is clear: because $K[[x_1,\dots,x_d]]$ is a domain.
\begin{center}\begin{minipage}{\textwidth}
\centering
\newcommand{\fu}[4]{\shade[top color=#1,bottom color=white] ($(#2)+(-#3,-#4)$) rectangle ($(#2)+(#3,0)$);}
\newcommand{\fv}[4]{\shade[bottom color=#1,top color=white] ($(#2)+(-#3,0)$) rectangle ($(#2)+(#3,#4)$);}
\resizebox{\textwidth}{!}{%
\begin{tikzpicture}[>=Stealth,line cap=round,line join=round]
\def\dx{1.6}\def\dy{1.8}\def\H{4.2}
\begin{scope}[xshift=-4cm,yshift=1.4cm]
\fill[green!6] (1,0) -- (6.5,0) -- (6.5+\dx,\dy) -- (1+\dx,\dy) -- cycle;
\draw[thick] (1,0) -- (6.5,0) -- (6.5+\dx,\dy) -- (1+\dx,\dy) -- cycle;
\coordinate (sp) at (2.4+0.45*\dx,0.45*\dy);
\draw[red!70!black,thick,dotted] (sp) circle (0.33);
\fill[green!40!black] (sp) circle (2pt) node[below right=2pt]{\scriptsize $s'$};
\node[red!70!black] at ($(sp)+(-0.58,0.24)$) {\scriptsize $\hat S'_{s'}$};
\node at (1.35,1.25) {$S'$};
\end{scope}
\begin{scope}[xshift=5cm]
\fill[red!4] (0,0) -- (8,0) -- (8+\dx,\dy) -- (\dx,\dy) -- cycle;
\draw[thick] (0,0) -- (8,0) -- (8+\dx,\dy) -- (\dx,\dy) -- cycle;
\draw[gray!60!black,line width=4.4pt,dash pattern=on 2.5pt off 1.5pt,line cap=butt] plot[smooth] coordinates {(5.200,0.000) (5.838,0.324) (6.526,0.648) (7.100,0.900) (7.656,1.188) (8.194,1.512) (8.600,1.800)};
\draw[gray!60!black,line width=4.4pt,dash pattern=on 2.5pt off 1.5pt,line cap=butt] plot[smooth] coordinates {(7.500,0.000) (7.420,0.360) (7.276,0.648) (7.100,0.900) (7.006,1.188) (7.094,1.512) (7.200,1.800)};
\draw[white,line width=3.6pt,line cap=butt] plot[smooth] coordinates {(5.200,0.000) (5.838,0.324) (6.526,0.648) (7.100,0.900) (7.656,1.188) (8.194,1.512) (8.600,1.800)};
\draw[white,line width=3.6pt,line cap=butt] plot[smooth] coordinates {(7.500,0.000) (7.420,0.360) (7.276,0.648) (7.100,0.900) (7.006,1.188) (7.094,1.512) (7.200,1.800)};
\node[gray!60!black,font=\scriptsize,below] at (6.4,-0.05) {$S_K^0\setminus h(S')$};
\draw[thick] (0,0) -- (8,0) -- (8+\dx,\dy) -- (\dx,\dy) -- cycle;
\coordinate (sK) at (2.4+0.45*\dx,0.45*\dy);
\draw[gray!60!black,dashed] ($(sK)+(0.25,0)$) ellipse [x radius=1.7, y radius=0.66];
\node[gray!60!black,font=\scriptsize] at ($(sK)+(1.3,0.3)$) {$U_{s,K}$};
\draw[red!70!black,thick,dotted] (sK) circle (0.33);
\fill[blue!60!black] (sK) circle (2pt) node[below right=2pt]{\scriptsize $s_K$};
\node[red!70!black] at ($(sK)+(-0.8,0.3)$) {\scriptsize $\hat S_{s_K}$};
\node at (2.2,-0.05) [below] {\small $S_K^0$};
\draw[thick,fill=gray!10] (10.6,0) -- (11.4,0) -- (11.4+0.5*\dx,0.5*\dy) -- (10.6+0.5*\dx,0.5*\dy) -- cycle;
\node[gray!60!black] at (11.4,-0.05) [below] {\scriptsize other components};
\draw[decorate,decoration={brace,mirror,amplitude=4pt}] (0,-0.6) -- (12.2,-0.6) node[midway,below=4pt]{$S_K$};
\draw[->] (9.5,4.4) -- node[right]{$\tau$} (9.5,2.0);
\begin{scope}[yshift=3.4cm]
\def\D{1.0}
\fill[gray!12,path fading=south] (0,-\D) -- (8,-\D) -- (8+\dx,\dy-\D) -- (8+\dx,\dy) -- (8,0) -- (0,0) -- cycle;
\fill[gray!12,path fading=north] (0,0) -- (8,0) -- (8+\dx,\dy) -- (8+\dx,\dy+\H) -- (\dx,\dy+\H) -- (0,\H) -- cycle;
\fv{gray!60}{\dx,\dy}{0.2pt}{\H} \fu{gray!60}{\dx,\dy}{0.2pt}{\D}
\def\shA{1.3+0.25*sin(60*\x)}
\fill[blue!40,opacity=0.35] plot[domain=0:8,smooth,samples=60] (\x,{\shA}) -- plot[domain=8:0,smooth,samples=60] ({\x+\dx},{\dy+\shA}) -- cycle;
\draw[blue!60!black] plot[domain=0:8,smooth,samples=60] (\x,{\shA});
\draw[blue!60!black,thin] plot[domain=0:8,smooth,samples=60] ({\x+\dx},{\dy+\shA});
\draw[black,line width=4.4pt,dash pattern=on 2.5pt off 1.5pt,line cap=butt] plot[smooth] coordinates {(5.200,1.114) (5.838,1.511) (6.526,1.935) (7.100,2.277) (7.656,2.635) (8.194,3.006) (8.600,3.317)};
\draw[black,line width=4.4pt,dash pattern=on 2.5pt off 1.5pt,line cap=butt] plot[smooth] coordinates {(7.500,1.550) (7.420,1.888) (7.276,2.115) (7.100,2.277) (7.006,2.475) (7.094,2.747) (7.200,2.998)};
\draw[white,line width=3.6pt,line cap=butt] plot[smooth] coordinates {(5.200,1.114) (5.838,1.511) (6.526,1.935) (7.100,2.277) (7.656,2.635) (8.194,3.006) (8.600,3.317)};
\draw[white,line width=3.6pt,line cap=butt] plot[smooth] coordinates {(7.500,1.550) (7.420,1.888) (7.276,2.115) (7.100,2.277) (7.006,2.475) (7.094,2.747) (7.200,2.998)};
\draw[blue!60!black] plot[domain=0:8,smooth,samples=60] (\x,{\shA});
\draw[blue!60!black,thin] plot[domain=0:8,smooth,samples=60] ({\x+\dx},{\dy+\shA});
\node[blue!60!black] at (0.9+0.7*\dx,{0.7*\dy+1.55}) {$Z$};
\coordinate (vK) at (2.4+0.45*\dx,{0.45*\dy+1.3+0.25*sin(144)});
\fill[red!70!black,opacity=0.6] (vK) ellipse [x radius=0.36, y radius=0.16, rotate=20];
\fill[blue!60!black] (vK) circle (2pt) node[below right=1pt]{\scriptsize $v_K$};
\fv{black}{0,0}{0.4pt}{\H} \fv{black}{8,0}{0.4pt}{\H} \fv{black}{8+\dx,\dy}{0.4pt}{\H} \fu{black}{0,0}{0.4pt}{\D} \fu{black}{8,0}{0.4pt}{\D} \fu{black}{8+\dx,\dy}{0.4pt}{\D}
\draw[thick] (0,0) -- (8,0) -- (8+\dx,\dy); \draw[thin,dashed,gray!60!black] (0,0) -- (\dx,\dy) -- (8+\dx,\dy);
\node at (0.35,\H-0.4) {$E_K$};
\fill[gray!12,path fading=south] (10.6,-\D) -- (11.4,-\D) -- (11.4+0.5*\dx,0.5*\dy-\D) -- (11.4+0.5*\dx,0.5*\dy) -- (11.4,0) -- (10.6,0) -- cycle;
\fill[gray!12,path fading=north] (10.6,0) -- (11.4,0) -- (11.4+0.5*\dx,0.5*\dy) -- (11.4+0.5*\dx,0.5*\dy+\H) -- (10.6+0.5*\dx,0.5*\dy+\H) -- (10.6,\H) -- cycle;
\fv{gray!60}{10.6+0.5*\dx,0.5*\dy}{0.2pt}{\H} \fu{gray!60}{10.6+0.5*\dx,0.5*\dy}{0.2pt}{\D}
\fv{black}{10.6,0}{0.4pt}{\H} \fv{black}{11.4,0}{0.4pt}{\H} \fv{black}{11.4+0.5*\dx,0.5*\dy}{0.4pt}{\H} \fu{black}{10.6,0}{0.4pt}{\D} \fu{black}{11.4,0}{0.4pt}{\D} \fu{black}{11.4+0.5*\dx,0.5*\dy}{0.4pt}{\D}
\draw[thick] (10.6,0) -- (11.4,0) -- (11.4+0.5*\dx,0.5*\dy); \draw[thin,dashed,gray!60!black] (10.6,0) -- (10.6+0.5*\dx,0.5*\dy) -- (11.4+0.5*\dx,0.5*\dy);
\end{scope}
\end{scope}
\draw[->] (3.0,2.0) -- node[above]{$h$} (6.6,1.44);
\draw[->,red!70!black,dashed,shorten <=10pt,shorten >=10pt] (sp) -- node[below]{\scriptsize $\hat h$} (sK);
\draw[->,red!70!black,dashed,shorten <=10pt,shorten >=6pt] (sK) to[bend left=25] node[left]{\scriptsize $\tilde v_K$} (vK);
\draw[->,green!40!black] (3.3,3.1) .. controls (3.8,4.0) and (4.4,5.0) .. node[above left=-1pt]{$u$, $\nabla u=0$} (5.5,5.4);
\begin{scope}[font=\small,xshift=5cm]
\node[align=left,anchor=west] (g) at (13.0,4.0) {$\overline{u(S')}\setminus u(S')$,\\ of dimension $<d$};
\draw[->] (g.west) .. controls (10.0,4.0) and (8.4,4.4) .. (5.96,4.91);
\end{scope}
\end{tikzpicture}%
}
\captionsetup{labelformat=simple,labelsep=period}

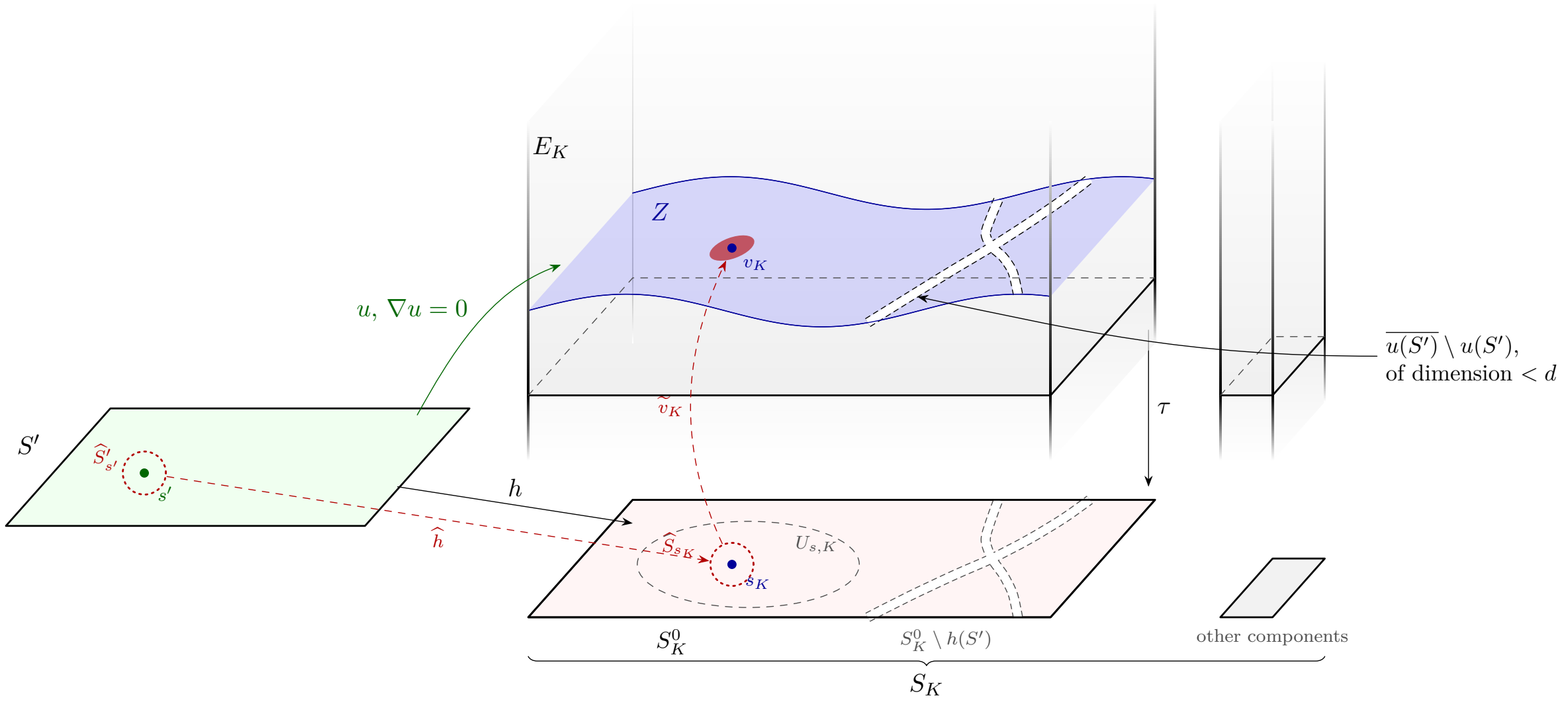
\captionof{figure}{The schemes in the proof of \Cref{lem:cover-algebraic}.}
\label{fig:cover-algebraic}
\end{minipage}\end{center}

{
Let $e'_j$ be the formal extension of $e_j(s_K)$ with $\nabla e'_j=0$, which exists uniquely by \cite[Proposition 2.2.1]{LL25}.
\begin{itemize}
\item $e_1,\dots,e_r$ is a basis of $\cE_K$ over $U_{s,K}:=U_s\times_RK$, and its dual basis $t_1,\dots,t_r$ of $\cE_K^\vee$ gives fiber coordinates of $E_K$ over $U_{s,K}$;
\item $e'_1,\dots,e'_r$ is a horizontal basis of $\cE_K$ over $\Spf\bigl(K[[x_1,\dots,x_d]]\bigr)$, and its dual basis $t'_1,\dots,t'_r$ gives fiber coordinates of $E_K$ over $\Spf\bigl(K[[x_1,\dots,x_d]]\bigr)$.
\end{itemize}
Write
\[
\tilde v_K=\sum_jc_je'_j\in\Gamma\bigl(\hat S_{s_K},\cE_K\bigr)=\bigoplus_jK[[x_1,\dots,x_d]]\,e'_j.
\]
A priori, $c_j\in K[[x_1,\dots,x_d]]$.
Since $\nabla\tilde v_K=0$ and $\nabla e'_j=0$ for all $j$, we have $dc_j=0$.
As $K$ has characteristic $0$, this is only possible when $c_j\in K$.
We shift the coordinates: Set
\[
z_j:=t'_j-c_j\in\Gamma\bigl(E_K|_{\hat S_{s_K}},\cO\bigr)=K[[x_1,\dots,x_d]][t'_1,\dots,t'_r].
\]
Since all $c_j\in K$, $z_j$ is again a fiber coordinate, shifted so that $v_K$ sits at $z=0$.
Write $x=(x_1,\dots,x_d)$ and $z=(z_1,\dots,z_r)$.
Then
\[
\hat{\cO}_{E_K,v_K}=K[[x,z]].
\]
The map $\tilde v_K\colon\Spf\hat\cO_{S_K,s_K}\rightarrow E_K$ factors uniquely through 
$\Spf\hat\cO_{S_K,s_K}\rightarrow \Spf\hat\cO_{E_K,v_K}$, which we again denote by $\tilde v_K$.
By construction of the coordinates $z_j$, the map
\[\tilde v_K^*\colon \Gamma(\Spf\hat\cO_{E_K,v_K},\cO)=K[[x,z]]\to \Gamma(\Spf\hat\cO_{S_K,s_K},\cO)=K[[x]]\]
sends $x_i\mapsto x_i$ and $z_j\mapsto0$.
It therefore induces a closed immersion of formal schemes 
\[\tilde v_K: \Spf\hat\cO_{S_K,s_K}\simeq \{z=0\}\hookrightarrow\Spf\hat\cO_{E_K,v_K}.\]
To put $\Spf\hat{\cO}_{Z,v_K}$ into the picture, first note that $\cO_{Z,v_K}=\cO_{E_K,v_K}/\cI_{Z,v_K}$, so
 \[\hat{\cO}_{Z,v_K}=K[[x,z]]/\hat{\cI}_Z\quad\text{ with }\hat{\cI}_Z:=\cI_ZK[[x,z]]. \]
To identify $\tilde v_K(\Spf\hat\cO_{S_K,s_K})$ with $\Spf\hat{\cO}_{Z,v_K}$, it is equivalent to show
\eq{lem:cover-algebraic-cl}{
\hat{\cI}_Z=(z_1,\dots,z_r).
}

\begin{center}\begin{minipage}{\textwidth}
\centering
\begin{tikzpicture}[>=Stealth,line cap=round,line join=round]
\def\W{6}\def\dx{1.6}\def\dy{1.8}
\fill[red!4] (0,0) -- (\W,0) -- (\W+\dx,\dy) -- (\dx,\dy) -- cycle;
\draw[thick] (0,0) -- (\W,0) -- (\W+\dx,\dy) -- (\dx,\dy) -- cycle;
\coordinate (sK) at (3+0.5*\dx,0.5*\dy);
\fill[blue!60!black] (sK) circle (2pt) node[below right=1pt]{\scriptsize $s_K$};
\node[anchor=north,font=\scriptsize] at (3.4,-0.1) {$\Spf\hat\cO_{S_K,s_K}=\Spf K[[x_1,x_2]]$};
\begin{scope}[yshift=4.2cm]
\def\Dn{1.3}\def\Up{2.4}
\draw[gray!60] (0,-\Dn) -- (0,\Up); \draw[gray!60] (\W,-\Dn) -- (\W,\Up);
\draw[gray!60] (\W+\dx,\dy-\Dn) -- (\W+\dx,\dy+\Up); \draw[gray!60,dashed] (\dx,\dy-\Dn) -- (\dx,\dy+\Up);
\draw[gray!60] (0,-\Dn) -- (\W,-\Dn) -- (\W+\dx,\dy-\Dn);
\draw[gray!60] (0,\Up) -- (\W,\Up) -- (\W+\dx,\dy+\Up) -- (\dx,\dy+\Up) -- cycle;
\draw[gray!60,dashed] (0,-\Dn) -- (\dx,\dy-\Dn) -- (\W+\dx,\dy-\Dn);
\fill[red!70!black,opacity=0.25] (0,0) -- (\W,0) -- (\W+\dx,\dy) -- (\dx,\dy) -- cycle;
\draw[red!70!black,thick] (0,0) -- (\W,0) -- (\W+\dx,\dy) -- (\dx,\dy) -- cycle;
\coordinate (vK) at (3+0.5*\dx,0.5*\dy);
\draw[->,gray!30!black] (vK) -- ++(0,2.0) node[above]{\scriptsize $z$};
\draw[->,gray!30!black] (vK) -- ++(1.6,0) node[right]{\scriptsize $x_1$};
\draw[->,gray!30!black] (vK) -- ++(0.55,0.62) node[above right=-2pt]{\scriptsize $x_2$};
\draw[->,green!40!black,thick] (1.0,0.35) -- (2.0,0.35);
\draw[->,green!40!black,thick] (1.0,0.35) -- (1.45,0.86);
\node[green!40!black,font=\scriptsize,anchor=west] at (2.05,0.3) {$\partial/\partial x_i=(\partial/\partial x_i)^h$};
\fill[blue!60!black] (vK) circle (2pt) node[below right=1pt]{\scriptsize $v_K$};
\node[anchor=north west,font=\scriptsize] at (\dx+0.15,\dy+\Up-0.25) {$\Spf\hat\cO_{E_K,v_K}=\Spf K[[x_1,x_2,z]]$};
\end{scope}
\node[red!70!black,font=\scriptsize,align=left,anchor=west] (lab) at (\W+\dx+0.6,5.1) {$\{z=0\}=\tilde v_K\bigl(\Spf\hat\cO_{S_K,s_K}\bigr)$\\ $\phantom{\{z=0\}}\stackrel{!}{=}
\Spf\hat\cO_{Z,v_K}$};
\draw[->,red!70!black] (lab.west) -- (\W+0.6*\dx,4.2+0.6*\dy);
\draw[->] (1.2,2.6) -- node[left]{\scriptsize $\tau$} (1.2,1.0);
\draw[->,red!70!black,dashed] (\W-0.6,0.4) -- node[pos=0.2,right]{\scriptsize $\tilde v_K$} (\W-0.6,4.5);
\end{tikzpicture}
\captionsetup{labelformat=simple,labelsep=period}

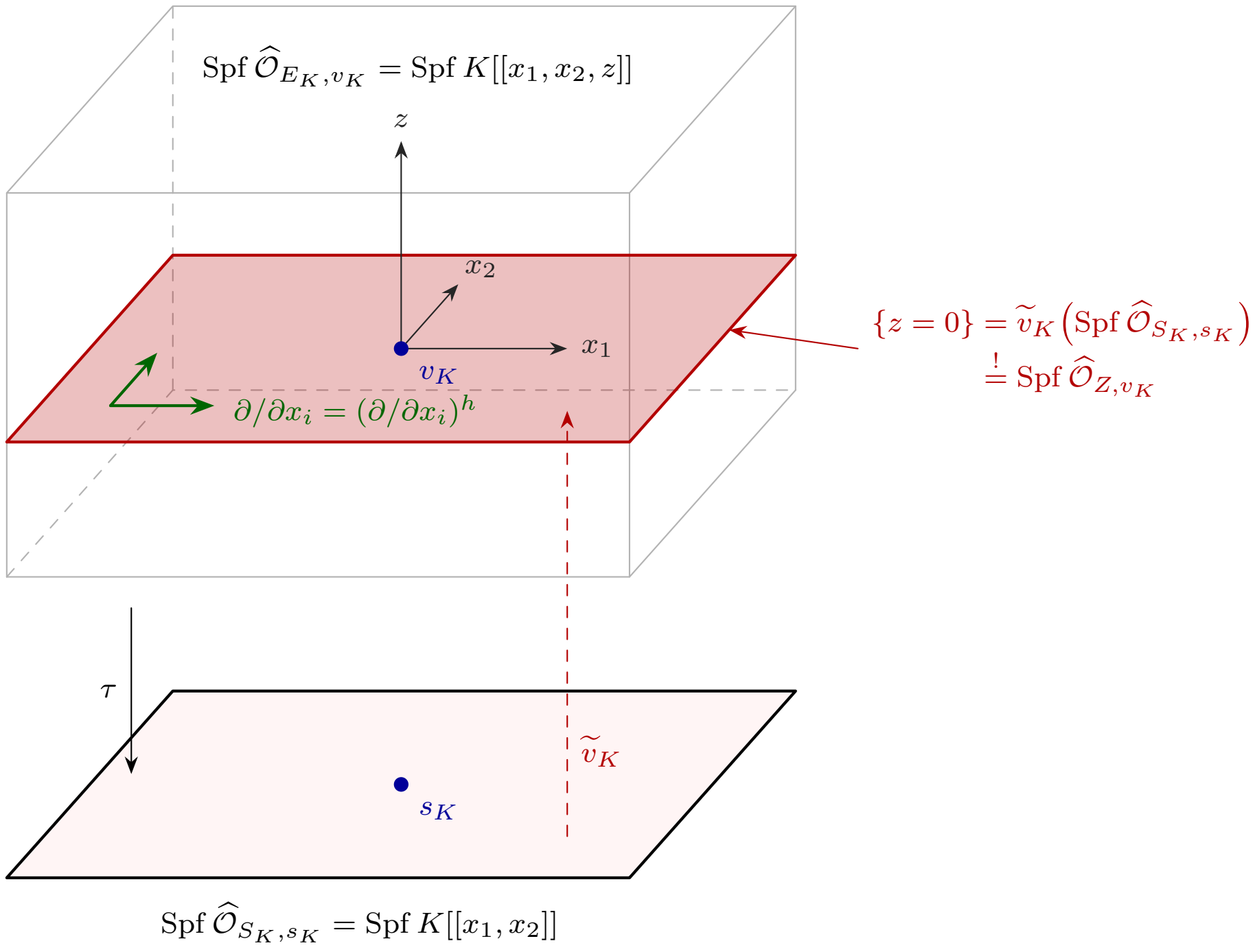
\captionof{figure}{The formal schemes in the proof of \Cref{lem:cover-algebraic}.}
\label{fig:cover-formal}
\end{minipage}\end{center}

\begin{claim}\label{claim:cover-gen}
$$\hat{\cI}_Z=\bigl(\hat{\cI}_Z\cap K[[z]]\bigr)K[[x,z]]$$
as ideals in $K[[x,z]]$.
\end{claim}
\begin{proof}[Proof of \Cref{claim:cover-gen}]
For $\alpha\in\Z_{\geq0}^d$, write $x^\alpha:=x_1^{\alpha_1}\cdots x_d^{\alpha_d}$, $\alpha!:=\alpha_1!\cdots\alpha_d!$ and $\partial_x^\alpha:=(\partial/\partial x_1)^{\alpha_1}\circ\dots\circ(\partial/\partial x_d)^{\alpha_d}$, acting on $K[[x,z]]$.

We first show that $\hat{\cI}_Z$ is stable under every $\partial_x^\alpha$. 
Writing
\[
\nabla e'_j=\sum_i\omega'_{ij}\otimes e'_i\quad\text{with }\omega'_{ij}\in\bigoplus_{k=1}^dK[[x]]\,dx_k
\]
as in \Cref{lem:connection-splitting}(1), all $\omega'_{ij}$ vanish as $\nabla e'_j=0$, so \eqref{eq:horizontal-lift} gives 
\[(\partial/\partial x_i)^h=\partial/\partial x_i.\] Moreover, for a vector field $\partial\in T_{S_K}$ near $s_K$, let $\partial^h$ be as in \Cref{lem:connection-splitting}. 
For a local section $f$ of $\cI_Z$, 
\Cref{lem:connection-splitting}(2) applied to $\tilde v_K$ gives 
\[\tilde v_K^*(\partial^h f)=\partial(\tilde v_K^*f)=0,\]
so $\partial^h f\in \Ker(\tilde v_K^*)=\cI_Z$ by \eqref{eq:cov-ker}.
In particular, $(\partial/\partial x_i)^h=\partial/\partial x_i$ maps $\hat{\cI}_Z$ into itself.
As derivations extend continuously to completions, $\hat{\cI}_Z$ is stable under every $\partial/\partial x_i$, and hence under every $\partial_x^\alpha$.

Next, let $y=(y_1,\dots,y_d)$ be further variables.
As $K$ has characteristic $0$, for every $f\in K[[x,z]]$ we have the Taylor expansion of $f$ around $x$,
\[
f(y,z)=\sum_{\alpha\in\Z_{\geq0}^d}\frac{(y-x)^\alpha}{\alpha!}\,(\partial_x^\alpha f)(x,z)\quad\text{in }K[[x,y,z]].
\]
Setting $y=0$ gives
\[
f(0,z)=\sum_{\alpha\in\Z_{\geq0}^d}\frac{(-x)^\alpha}{\alpha!}\,(\partial_x^\alpha f)(x,z),
\]
and exchanging the roles of $x$ and $y$ and then setting $y=0$, i.e., expanding around $0$, gives
\[
f(x,z)=\sum_{\alpha\in\Z_{\geq0}^d}\frac{x^\alpha}{\alpha!}\,(\partial_x^\alpha f)(0,z).
\]
Both sums converge $x$-adically.

Now let $f\in\hat{\cI}_Z$.
Every term of the expansion of $f(0,z)$ lies in $\hat{\cI}_Z$ by the stability above, so $f(0,z)\in\hat{\cI}_Z\cap K[[z]]$.
Applying this to $\partial_x^\alpha f\in\hat{\cI}_Z$ shows that every coefficient $(\partial_x^\alpha f)(0,z)$ of the expansion of $f$ around $0$ lies in $\hat{\cI}_Z\cap K[[z]]$. 
By Krull intersection theorem applied to $A/\mathfrak a$, every ideal $\mathfrak a$ of a noetherian ring $A$ is closed in the $\mathfrak b$-adic topology for every ideal $\mathfrak b$ contained in the Jacobson radical of $A$.
Hence the ideal $\bigl(\hat{\cI}_Z\cap K[[z]]\bigr)K[[x,z]]$ is closed  in the $(x)$-adic topology, and therefore it contains $f$ as it is a limit point.
The reverse inclusion is clear.

\end{proof}

\begin{claim}\label{claim:cover-z}
\[
\hat{\cI}_Z\cap K[[z]]=(z_1,\dots,z_r)
\]
as ideals in $K[[z]]$.
\end{claim}
\begin{proof}[Proof of \Cref{claim:cover-z}]
It follows from \eqref{eq:cov-ker} that  $\hat{\cI}_Z\cap K[[z]]\subset(z_1,\dots,z_r)$. 
It remains to show the converse.
Set
\[B:=K[[z]]/(\hat{\cI}_Z\cap K[[z]]).\]
By \Cref{claim:cover-gen},
\[
\hat{\cO}_{Z,v_K}=K[[x,z]]/\bigl(\hat{\cI}_Z\cap K[[z]]\bigr)K[[x,z]].
\]
As $K[[z]]$ is noetherian, $\hat{\cI}_Z\cap K[[z]]$ is finitely generated, so 
\[\bigl(\hat{\cI}_Z\cap K[[z]]\bigr)K[[x,z]]
= (\hat{\cI}_Z\cap K[[z]])[[x]]
\]  
by \cite[\href{https://stacks.math.columbia.edu/tag/00MA}{Tag 00MA}(3)]{SP} (applied to the finite $K[[z]][x]$-module $B[x]$ and the ideal $(x)$).
Hence $\hat{\cO}_{Z,v_K}\simeq B[[x]]$.
So it remains to show that $B=K$.

On the one hand, as $v_K$ is a closed point of the scheme $Z$ of finite type over $K$, which is irreducible as the closure of the image of the irreducible scheme $\Spec K[[x_1,\dots,x_d]]$, $\dim\hat{\cO}_{Z,v_K}=\dim Z=d$.
As $\dim B[[x]]=\dim B+d$, we get $\dim B=0$.

On the other hand, the local ring $\cO_{Z,v_K}$ is reduced and, being essentially of finite type over $K$, Nagata, so $\hat{\cO}_{Z,v_K}$ is reduced \cite[\href{https://stacks.math.columbia.edu/tag/07NZ}{Tag 07NZ}]{SP}.
So $B$, as a subring of $B[[x]]\simeq \hat \cO_{Z,v_K}$, is also reduced.
A reduced local $K$-algebra of dimension $0$ with residue field $K$ is $K$, so $B=K$, i.e., $\hat{\cI}_Z\cap K[[z]]=(z_1,\dots,z_r)$.
\end{proof}
\Cref{claim:cover-gen,claim:cover-z} together prove \eqref{lem:cover-algebraic-cl}.
We therefore establish \ref{item:coveralg-2}.

Finally, let $K\subset\C$ be an embedding.
As $v_K$ is a $K$-point, $\hat\cO_{Z_{\C},v_{\C}}=\hat\cO_{Z,v_K}\hat\otimes_K\C$.
So by \ref{item:coveralg-2}, the image of $\tilde v_{\C}$ is the completion at $v_{\C}$ of the closed subscheme $Z_{\C}$ of $E_{\C}$, and $\dim Z_{\C}=\dim Z=d$ by \ref{item:coveralg-1}.
This proves \ref{item:coveralg-3}.}
\end{proof}

\subsection{Hodge and conjugate filtrations}
\label{subsec:katz}

For a morphism of schemes $\pi\colon X\to S$, the Hodge filtration on $R^m\pi_*\Omega^\bullet_{X/S}$ is
\[F^i_{\mathrm{Hdg}}R^m\pi_*\Omega^\bullet_{X/S}:=\operatorname{Im}\bigl(R^m\pi_*\Omega^{\geq i}_{X/S}\to R^m\pi_*\Omega^\bullet_{X/S}\bigr).\]
Let $\pi\colon X\to S$ be a smooth morphism, with $S$ of characteristic $p>0$.
The conjugate filtration on $R^{m}\pi_*\Omega^\bullet_{X/S}$ is
\[F^i_{\conj}R^m\pi_*\Omega^\bullet_{X/S}:=\operatorname{Im}\bigl(R^m\pi_*\tau_{\leq m-i}\Omega^\bullet_{X/S}\to R^m\pi_*\Omega^\bullet_{X/S}\bigr),\]
where $\tau_{\leq\bullet}$ is the canonical filtration of $\Omega^\bullet_{X/S}$; see \cite[\S2.3]{Kat72} and \cite[\S4.2]{LL25}.

\begin{lemma}
\label{lem:GM-filtrations}
Let $\pi\colon X\to S$ be a smooth morphism, with $S$ smooth over a ring $R$, and let $\nabla_{\GM}$ be the Gauss--Manin connection on $R^m\pi_*\Omega^\bullet_{X/S}$.
\begin{enumerate}[label=\textup{(\arabic*)}]
\item \emph{Griffiths transversality:}
\[\nabla_{\GM}\bigl(F^i_{\mathrm{Hdg}}R^m\pi_*\Omega^\bullet_{X/S}\bigr)\subset F^{i-1}_{\mathrm{Hdg}}R^m\pi_*\Omega^\bullet_{X/S}\otimes_{\cO_S}\Omega^1_{S/R}.\]
\item If $R$ has characteristic $p>0$, then
\[\nabla_{\GM}\bigl(F^i_{\conj}R^m\pi_*\Omega^\bullet_{X/S}\bigr)\subset F^i_{\conj}R^m\pi_*\Omega^\bullet_{X/S}\otimes_{\cO_S}\Omega^1_{S/R}.\]
The $p$-curvature $\psi_p$ of $\nabla_{\GM}$ satisfies
\[\psi_p\bigl(F^i_{\conj}R^m\pi_*\Omega^\bullet_{X/S}\bigr)\subset F_{\mathrm{abs}}^*\Omega^1_{S/R}\otimes_{\cO_S}F^{i+1}_{\conj}R^m\pi_*\Omega^\bullet_{X/S},\]
or equivalently, for every local vector field $\partial\in T_S$,
\begin{equation}
\label{eq:psi-conj}
\psi_p(\partial)\bigl(F^i_{\conj}R^m\pi_*\Omega^\bullet_{X/S}\bigr)\subset F^{i+1}_{\conj}R^m\pi_*\Omega^\bullet_{X/S}.
\end{equation}
\end{enumerate}
\end{lemma}
\begin{proof}
(1) is \cite[(1.4.1.6)]{Kat72}.
For (2), $F^\bullet_{\conj}$ is horizontal with respect to $\nabla_{\GM}$ by \cite[(2.3.0)]{Kat72}.
The $p$-curvature statement is \cite[\S3.1, (3.1.2)]{Kat72}.
\end{proof}

\begin{lemma}
\label{lem:dR-flat-bc}
Let $\pi\colon X\to S$ be a qcqs morphism, locally of finite type.
\begin{enumerate}[label=\textup{(\arabic*)},ref=\arabic*]
\item\label{item:dR-bc} \emph{Flat base change of the de Rham cohomology.} Let $g\colon S'\to S$ be a flat morphism, and let $\pi'\colon X':=X\times_SS'\to S'$ be the base change of $\pi$.
Then for every $m$,
\[
g^*R^m\pi_*\Omega^\bullet_{X/S}\simeq R^m\pi'_*\Omega^\bullet_{X'/S'}.
\]
\item\label{item:dR-const} \emph{Smooth base change trivializes $\nabla_{\GM}$.}
If moreover $\pi$ and $g$ are smooth, then for every $m$ there is an isomorphism
\[
(R^m\pi'_*\Omega^\bullet_{X'/S'},\nabla_{\GM})\simeq(g^*R^m\pi_*\Omega^\bullet_{X/S},1\otimes d)
\]
of quasi-coherent $\cO_{S'}$-modules with connection relative to $S$, where
\[
\nabla_{\GM}\colon R^m\pi'_*\Omega^\bullet_{X'/S'}\to R^m\pi'_*\Omega^\bullet_{X'/S'}\otimes_{\cO_{S'}}\Omega^1_{S'/S}\quad\text{and}\quad d=d_{S'/S}\colon\cO_{S'}\to\Omega^1_{S'/S}
\]
are the Gauss--Manin connection of $\pi'$ relative to $S$ and the differential of $S'/S$, respectively.
\end{enumerate}
\end{lemma}

\begin{proof}
(\ref{item:dR-bc}) Write $g'\colon X'\to X$ for the projection.
As $\Omega^j_{X'/S'}=g'^*\Omega^j_{X/S}$, flat base change \cite[\href{https://stacks.math.columbia.edu/tag/02KH}{Tag 02KH}]{SP} gives $g^*R^q\pi_*\Omega^j_{X/S}\simeq R^q\pi'_*\Omega^j_{X'/S'}$ for all $q,j$.
As $g$ is flat, $g^*$ is exact, so applying $g^*$ to the long exact sequence induced by $0\to\Omega^{\geq j+1}_{X/S}\to\Omega^{\geq j}_{X/S}\to\Omega^j_{X/S}[-j]\to0$ gives a long exact sequence.
The assertion is local on $S$, so we may assume that $S$ is affine.
Then $X$ is quasi-compact, and as $\pi$ is locally of finite type, $\Omega^{\geq N}_{X/S}=0$ for some $N$, hence $\Omega^{\geq N}_{X'/S'}=0$, so descending induction on $j$ and the five lemma give
\[
g^*R^m\pi_*\Omega^{\geq i}_{X/S}\xrightarrow{\sim}R^m\pi'_*\Omega^{\geq i}_{X'/S'}
\]
for all $i$ and $m$.
The case $i=0$ is the assertion.

(\ref{item:dR-const}) The sheaf isomorphism is (\ref{item:dR-bc}).
We have, for every $n$,
\begin{equation}\label{eq:kunneth-omega}
\Omega^n_{X'/S}=\bigoplus_{a+b=n}g'^*\Omega^a_{X/S}\otimes_{\cO_{X'}}\pi'^*\Omega^b_{S'/S},
\end{equation}
and, by the Leibniz rule, for local sections $\alpha$ of $\Omega^a_{X/S}$ and $\beta$ of $\Omega^b_{S'/S}$,
\[
d(g'^*\alpha\wedge\pi'^*\beta)=g'^*(d\alpha)\wedge\pi'^*\beta+(-1)^a\,g'^*\alpha\wedge\pi'^*(d\beta).
\]
On the other hand, by \cite[\S3]{Kat70}, $\nabla_{\GM}$ is the connecting map on $R^m\pi'_*$ of the short exact sequence of complexes
\begin{equation}\label{eq:katz-oda}
0\to\pi'^*\Omega^1_{S'/S}\otimes\Omega^{\bullet-1}_{X'/S'}\to\Omega^\bullet_{X'/S}/(\pi'^*\Omega^2_{S'/S}\wedge\Omega^{\bullet-2}_{X'/S})\to\Omega^\bullet_{X'/S'}\to0.
\end{equation}
By \eqref{eq:kunneth-omega} and the Leibniz rule, \eqref{eq:katz-oda} is obtained by applying $g'^{-1}\Omega^\bullet_{X/S}\otimes-$ to
\[
0\to\Omega^1_{S'/S}[-1]\to(\cO_{S'}\xrightarrow{d}\Omega^1_{S'/S})\to\cO_{S'}\to0.
\]
Hence, for a local cocycle $\omega$ of $\Omega^\bullet_{X'/S'}$ and a lift $\widetilde\omega\in\Omega^\bullet_{X'/S}/(\pi'^*\Omega^2_{S'/S}\wedge\Omega^{\bullet-2}_{X'/S})$ of $\omega$, we have
\begin{equation}\label{eq:GM-lift}
\nabla_{\GM}([\omega])=[d\widetilde\omega]\in R^{m+1}\pi'_*(\pi'^*\Omega^1_{S'/S}\otimes\Omega^{\bullet-1}_{X'/S'}).
\end{equation}

By (\ref{item:dR-bc}), we have
\[
R^m\pi'_*\Omega^\bullet_{X'/S'}\simeq g^*R^m\pi_*\Omega^\bullet_{X/S}=\cO_{S'}\otimes_{g^{-1}\cO_S}g^{-1}R^m\pi_*\Omega^\bullet_{X/S},
\]
so locally a class in $R^m\pi'_*\Omega^\bullet_{X'/S'}$ is a sum of classes $[\alpha]\otimes f$ with $[\alpha]\in R^m\pi_*\Omega^\bullet_{X/S}$ and $f\in\cO_{S'}$, and
we have the following equalities of local sections of $R^m\pi'_*\Omega^\bullet_{X'/S'}\otimes_{\cO_{S'}}\Omega^1_{S'/S}$:
\begin{align*}
\nabla_{\GM}([\alpha]\otimes f)&=[d(f\,g'^*\alpha)]&&(\text{by \eqref{eq:GM-lift} with $\widetilde\omega=f\,g'^*\alpha$})\\
&=[df\wedge g'^*\alpha+f\,g'^*(d\alpha)]&&(\text{Leibniz rule})\\
&=[df\wedge g'^*\alpha]&&(d\alpha=0)\\
&=[\alpha]\otimes df&&(\text{projection formula}),
\end{align*}
where by the projection formula we mean the isomorphism
\[
R^{m+1}\pi'_*(\pi'^*\Omega^1_{S'/S}\otimes\Omega^{\bullet-1}_{X'/S'})\simeq R^m\pi'_*\Omega^\bullet_{X'/S'}\otimes_{\cO_{S'}}\Omega^1_{S'/S}.
\]
So $\nabla_{\GM}=1\otimes d$.
\end{proof}

\begin{lemma}
\label{lem:filt-degen}
Let $K$ be a complex in an abelian category with a finite decreasing filtration $F$, i.e., $F^aK=K$ and $F^bK=0$ for some $a\leq b$.
Let \[E_1^{l,j}=H^{l+j}(\mathrm{gr}^l_FK)\Rightarrow H^{l+j}(K)\] be its spectral sequence.
Then the following are equivalent.
\begin{enumerate}[label=\textup{(\arabic*)}]
\item The spectral sequence degenerates at $E_1$.
\item All maps $H^m(F^lK)\to H^m(K)$ are injective.
\item For all $l$ and $m$, the sequence \[0\to H^m(F^{l+1}K)\to H^m(F^lK)\to H^m(\mathrm{gr}^l_FK)\to0\] induced by $0\to F^{l+1}K\to F^lK\to\mathrm{gr}^l_FK\to0$ is exact.
\end{enumerate}
\end{lemma}

\begin{proof}
By the Freyd--Mitchell embedding theorem, we may assume that we work in a module category, and may hence argue with elements.
\begin{itemize}
\item[(1)$\Leftrightarrow$(2)] This is \cite[Proposition (1.3.2)]{Del71}.
There, (2) is stated as the strict compatibility of the differential $d$ of $K$ with $F$.
A filtered map $f\colon(A,F)\to(B,F)$ is \emph{strict} if $f(F^lA)=f(A)\cap F^lB$ for all $l$ \cite[(1.1.5)]{Del71}.
For $d\colon K^m\to K^{m+1}$, this says $d(F^lK^m)=dK^m\cap F^lK^{m+1}$, i.e., a boundary lying in $F^l$ is the boundary of something in $F^l$.
This is exactly the injectivity of $H(F^lK)\to H(K)$.
Indeed, a cocycle $x\in F^lK$ is $0$ in $H(K)$ if and only if $x=dy$ for some $y\in K$.
It is $0$ in $H(F^lK)$ if and only if $x=dy'$ for some $y'\in F^lK$.
So the injectivity for all $l$ is the inclusion $dK\cap F^lK\subset d(F^lK)$, and the reverse inclusion always holds.
\item[(1)$\Leftrightarrow$(3)] Write the long exact sequences associated with $0\to F^{l+1}K\to F^lK\to\mathrm{gr}^l_FK\to0$ as \[\cdots\to H^m(F^{l+1}K)\xrightarrow{\alpha}H^m(F^lK)\xrightarrow{\beta}H^m(\mathrm{gr}^l_FK)\xrightarrow{\delta}H^{m+1}(F^{l+1}K)\to\cdots.\] An element of $E_r^{l,j}$ is represented by some $x\in H^{l+j}(\mathrm{gr}^l_FK)$ such that $\delta(x)=\alpha^{r-1}(y)$ for some $y\in H^{l+j+1}(F^{l+r}K)$, and then $d_r(x)$ is the class of $\beta(y)$. 
By exactness, (3) holds if and only if $\delta=0$, if and only if all maps $\alpha$ are injective.

(3)$\Rightarrow$(1): If (3) holds, then $\delta=0$.
Hence $d_r=0$ for all $r\geq1$.

(1)$\Rightarrow$(3): Suppose $d_r=0$ for all $r\geq1$, so that $E_r=E_1$ for all $r$.
Let $x\in H^{l+j}(\mathrm{gr}^l_FK)$.
Suppose that $\delta(x)=\alpha^{r-1}(y_r)$ for some $y_r\in H^{l+j+1}(F^{l+r}K)$; for $r=1$, take $y_1:=\delta(x)$.
Then $\beta(y_r)$ represents $d_r(x)=0$ in $E_r=E_1$, so $\beta(y_r)=0$.
By exactness, $y_r=\alpha(y_{r+1})$ for some $y_{r+1}$, so $\delta(x)=\alpha^{r}(y_{r+1})$.
Since $F^{l+r}K=0$ for $l+r\geq b$, this gives $\delta(x)=0$.
Hence $\delta=0$, i.e., (3) holds.
\end{itemize}
\end{proof}

\begin{corollary}
\label{lem:HdR-degen}
Let $\pi\colon X\to S$ be a smooth morphism of relative dimension $d$.
Then the following are equivalent.
\begin{enumerate}[label=\textup{(\arabic*)}]
\item The Hodge--de Rham spectral sequence
\[E_1^{i,j}=R^j\pi_*\Omega^i_{X/S}\Rightarrow R^{i+j}\pi_*\Omega^\bullet_{X/S}\]
degenerates at $E_1$.
\item All maps $R^m\pi_*\Omega^{\geq i}_{X/S}\to R^m\pi_*\Omega^\bullet_{X/S}$ are injective.
\item For all $i$ and $m$, the sequence
\begin{equation}
\label{eq:Hdg-fil-seq}
0\to R^{m}\pi_*\Omega^{\geq i+1}_{X/S}\to R^{m}\pi_*\Omega^{\geq i}_{X/S}\to R^{m-i}\pi_*\Omega^i_{X/S}\to0
\end{equation}
induced by $0\to\Omega^{\geq i+1}_{X/S}\to\Omega^{\geq i}_{X/S}\to\Omega^i_{X/S}[-i]\to0$ is exact.
\end{enumerate}
\end{corollary}

\begin{proof}
Let $K$ be $\pi_*$ of a filtered injective resolution of $\Omega^\bullet_{X/S}$ with the filtration by the $\Omega^{\geq i}_{X/S}$, and let $F^iK$ be $\pi_*$ of the $i$-th filtered piece.
This filtration is finite, as $\Omega^{\geq0}_{X/S}=\Omega^\bullet_{X/S}$ and $\Omega^{\geq d+1}_{X/S}=0$.
Then \[H^m(F^iK)=R^m\pi_*\Omega^{\geq i}_{X/S},\quad H^m(K)=R^m\pi_*\Omega^\bullet_{X/S}\quad\text{and}\quad H^m(\mathrm{gr}^i_FK)=R^{m-i}\pi_*\Omega^i_{X/S},\] and the spectral sequence of $(K,F)$ is the Hodge--de Rham spectral sequence.
So the claim is \Cref{lem:filt-degen}.
\end{proof}

\begin{definition}[Higgs field]
\label{def:higgs}
Let $S$ be smooth over a base ring $R$, let $\nabla$ be a connection on a vector bundle $\cE$ on $S$, and let $F_{\mathrm{Hdg}}^\bullet$ be a filtration of $\cE$ such that $\nabla$ sends $F_{\mathrm{Hdg}}^i$ into $F_{\mathrm{Hdg}}^{i-1}\otimes_{\cO_S}\Omega^1_{S/R}$.
The \emph{Higgs field} is the associated graded of $\nabla$,
\[
\theta:=\mathrm{gr}_{F_{\mathrm{Hdg}}}\nabla\colon \bigoplus_i\mathrm{gr}^i_{F_{\mathrm{Hdg}}}\cE\longrightarrow \bigoplus_i\mathrm{gr}^{i-1}_{F_{\mathrm{Hdg}}}\cE\otimes_{\cO_S}\Omega^1_{S/R}.
\]
\end{definition}

\begin{definition}[Order of vanishing]
\label{def:order-vanishing}
Let $S$ be a smooth $R$-scheme, let $s\in S(R)$, and write $\cI_s$ for the ideal sheaf of $s$.
Let $\mathfrak{p}$ be a closed point of $\operatorname{Spec} R$ of residue characteristic $p$, and let $\cE$ be a coherent sheaf on $S\otimes_R k(\mathfrak{p})$ (or a morphism of such sheaves).
\begin{enumerate}[label=\textup{(\alph*)}]
\item
A section of $\cE$ \emph{vanishes to order $N$ at $s$ modulo $\mathfrak{p}$} if it vanishes identically on $V(\cI_s^N,\mathfrak{p})\subset S\otimes_R k(\mathfrak{p})$, or equivalently if it lies in $\cI_s^N\cdot\cE$.
\item
A section of a coherent sheaf on $S$ (or after base change to a fiber of $S$) \emph{vanishes to order $N$ at $s$} if it vanishes on $V(\cI_s^N)$, or equivalently if it lies in $\cI_s^N\cdot\cE$.
\end{enumerate}
If $g\colon\{\mathrm{primes}\}\to\Z_{\geq 0}$, then ``vanishes to order $g(p)$ at $s$ modulo $\mathfrak{p}$'' means the case $N=g(p)$ in (a).
\end{definition}

A variant of the following lemma already appears in the proof of \cite[Theorem 4.1.1]{LL25} for $\mathrm{Sym}^{l}$ of $R^1\pi_*\Omega^\bullet$; here it is stated for $R^n\pi_*\Omega^\bullet$.

\begin{lemma}[Katz]
\label{prop:katz-F2}
Let $\pi\colon X\to S$ be a smooth projective morphism, with $S$ smooth over a ring $R$ of characteristic $p>0$.
Assume:
\begin{enumerate}[label=\textup{(\arabic*)},ref=\arabic*,labelindent=0pt,labelwidth=\widthof{\textup{(LFBS)}},labelsep=0.5em,leftmargin=!,align=left]
\item[\textup{(LFBS)}]\phantomsection
\label{item:lf-hodge}
$R^a\pi_*\Omega^b_{X/S}$ are locally free for all $a,b\geq0$, and their formation commutes with base change.%
\item[\textup{(HdR)}]\phantomsection
\label{item:lf-degen}
The Hodge--de Rham spectral sequence degenerates at $E_1$.
\end{enumerate}
Let $q,j\geq0$ and $n=q+j$.
If $\mathrm{gr}^q_{F_{\conj}}\psi_p$ of $\nabla_{\GM}$ on $R^n\pi_*\Omega^\bullet_{X/S}$ vanishes to order $N$ at a point $s$ of $S$, then the restriction of $\theta$ to $\mathrm{gr}^j_{F_{\mathrm{Hdg}}}R^n\pi_*\Omega^\bullet_{X/S}$ vanishes to order $\lfloor N/p\rfloor$ at $s$.
\end{lemma}

\begin{proof}
Write $F_{\mathrm{abs}}\colon S\to S$ for absolute Frobenius.
The following diagram is \cite[Theorem 4.2.1]{LL25}, due to Katz.
By \cite[Theorem 3.2]{Kat72} for $(a,b)=(q,j)$ and \cite[(1.4.1.7.2)]{Kat72}, the diagram
\begin{equation}
\label{eq:katz-square}
\begin{tikzcd}[column sep=4.4em, row sep=2.6em, ampersand replacement=\&]
F_{\mathrm{abs}}^*\mathrm{gr}^{j}_{F_{\mathrm{Hdg}}}R^{n}\pi_*\Omega^\bullet_{X/S}
\arrow[r, "\scriptstyle -F_{\mathrm{abs}}^*\theta"]
\arrow[d, "C^{-1}"', "\wr"]
\&
F_{\mathrm{abs}}^*\Omega^1_{S/R}\otimes F_{\mathrm{abs}}^*\mathrm{gr}^{j-1}_{F_{\mathrm{Hdg}}}R^{n}\pi_*\Omega^\bullet_{X/S}
\arrow[d, "1\otimes C^{-1}", "\wr"']
\\
\mathrm{gr}^{q}_{F_{\conj}}R^{n}\pi_*\Omega^\bullet_{X/S}
\arrow[r, "\scriptstyle \psi_p"']
\&
F_{\mathrm{abs}}^*\Omega^1_{S/R}\otimes\mathrm{gr}^{q+1}_{F_{\conj}}R^{n}\pi_*\Omega^\bullet_{X/S}
\end{tikzcd}
\end{equation}
commutes,\footnote{We follow the sign of \cite[Theorem 4.2.1]{LL25}.
In \cite[Theorem 3.2]{Kat72}, the lower horizontal arrow is $(-1)^{b+1}F_{\mathrm{abs}}^*\rho$, where $\rho$ is the cup product with the Kodaira--Spencer class, and $\mathrm{gr}\,\nabla=\rho$ by \cite[(1.4.1.7.2)]{Kat72}; this differs from the sign above when $j$ is odd.
This sign error in \cite[Theorem 3.2]{Kat72} is also noted in \cite[\S4.2.4 and footnote~11]{Dri25} and \cite[Remark~1.11]{Bar26}.} where $C^{-1}$ is the inverse of the Cartier isomorphism \cite[(2.3.1.2.1)]{Kat72} composed with $F_{\mathrm{abs}}^*$ of the isomorphism $\mathrm{gr}^j_{F_{\mathrm{Hdg}}}R^n\pi_*\Omega^\bullet_{X/S}\simeq R^q\pi_*\Omega^j_{X/S}$ (and similarly on the right).
Hence $F_{\mathrm{abs}}^*$ of the restriction of $\theta$ to $\mathrm{gr}^j_{F_{\mathrm{Hdg}}}R^n\pi_*\Omega^\bullet_{X/S}$ vanishes to order $N$ at $s$.
The lemma is just a restatement of the module structures of these two horizontal maps.
\end{proof}

\begin{lemma}
\label{lem:conj-bc}
Let $\pi\colon X\to S$ be as in \Cref{prop:katz-F2}, and let $m\geq0$.
Then the following hold.
\begin{enumerate}[label=\textup{(\arabic*)},ref=\arabic*]
\item
The conjugate spectral sequence
\[
{}_{\conj}E_2^{a,b}=F_{\mathrm{abs}}^*R^a\pi_*\Omega^b_{X/S}\Longrightarrow R^{a+b}\pi_*\Omega^\bullet_{X/S}
\]
degenerates at $E_2$.
\item
Each $F^a_{\conj}R^m\pi_*\Omega^\bullet_{X/S}$ is locally a direct summand, and the Cartier isomorphism gives
\[
\mathrm{gr}^a_{F_{\conj}}R^m\pi_*\Omega^\bullet_{X/S}\cong F_{\mathrm{abs}}^*R^a\pi_*\Omega^{m-a}_{X/S}.
\]
\item
The formation of $F^\bullet_{\conj}$ commutes with base change along any morphism $T\to S$.
\end{enumerate}
\end{lemma}
\begin{proof}
By \hyperref[item:lf-hodge]{(LFBS)} and \hyperref[item:lf-degen]{(HdR)}, (1), the isomorphism in (2), and (3) follow from \cite[(2.3.2)]{Kat72}.
The $F^a_{\conj}$ are locally direct summands by \cite[(1.11), Definition~1.13]{Wed08}.
\end{proof}

The following lemma is a variant of \cite[(5.0.1)]{Kat72}.
We prove it without mixed Hodge theory.

\begin{lemma}
\label{rmk:local-freeness}
Let $R$ be a finitely generated $\Z$-algebra that is a domain of characteristic $0$, let $S$ be a smooth $R$-scheme of finite type, and let $\pi\colon X\to S$ be a smooth projective morphism.
There exists a nonzero element $a\in R$ such that, after replacing $R$ by $R[1/a]$ and base changing $S$ and $X$ along this localization, the following hold.
\begin{enumerate}[label=\textup{(\arabic*)},ref=\arabic*]
\item
\label{item:lf-smooth}
$R$ is smooth over $\Z$.
\item
\label{item:lf-221-lfbs}
Condition \hyperref[item:lf-hodge]{(LFBS)} of \Cref{prop:katz-F2} holds for $\pi$.
\item
\label{item:lf-221-hdr}
Condition \hyperref[item:lf-degen]{(HdR)} of \Cref{prop:katz-F2} holds for $\pi$ and for the base changes $\pi_{k(\mathfrak p)}$ below, i.e., the Hodge--de Rham spectral sequence of $\pi$
\begin{equation}
\label{eq:HdR-f}
E_1^{i,j}=R^j\pi_*\Omega^i_{X/S}\Rightarrow R^{i+j}\pi_*\Omega^\bullet_{X/S},
\end{equation}and that of its base change $\pi_{k(\mathfrak p)}\colon X_{k(\mathfrak p)}\to S_{k(\mathfrak p)}$ along $\Spec k(\mathfrak{p})\to\Spec R$, for $\mathfrak{p}$ the generic point or a maximal ideal of $R$
\begin{equation}
\label{eq:HdR-fp}
E_1^{i,j}=R^j\pi_{k(\mathfrak p)*}\Omega^i_{X_{k(\mathfrak p)}/S_{k(\mathfrak p)}}\Rightarrow R^{i+j}\pi_{k(\mathfrak p)*}\Omega^\bullet_{X_{k(\mathfrak p)}/S_{k(\mathfrak p)}},
\end{equation}degenerate at $E_1$.
\item
\label{item:lf-221-summand}
Each $R^m\pi_*\Omega^\bullet_{X/S}$ is locally free, and each $F^i_{\mathrm{Hdg}}R^m\pi_*\Omega^\bullet_{X/S}=R^m\pi_*\Omega^{\geq i}_{X/S}$ is locally a direct summand of it.
For every maximal ideal $\mathfrak{p}$ of $R$, the restriction of $R^m\pi_*\Omega^\bullet_{X/S}$ to $S_{k(\mathfrak p)}$ is $R^m\pi_{k(\mathfrak p)*}\Omega^\bullet_{X_{k(\mathfrak p)}/S_{k(\mathfrak p)}}$.
\end{enumerate}
\end{lemma}

\begin{proof}
As $R$ is a domain of characteristic $0$, $\operatorname{Spec}R\to\operatorname{Spec}\Z$ is smooth at the generic point; as the smooth locus is open, $R[1/a]$ is smooth over $\Z$ for some $0\neq a\in R$.

For $t\in S$ put
\[
h^{i,j}(t):=\dim_{k(t)}H^j(X_t,\Omega^i_{X_t});
\]
each $h^{i,j}$ is upper semicontinuous on $S$ \cite[\S5, Corollary (a)]{Mum74}.
Put $S_{\Q}:=S\otimes_{\Z}\Q$, smooth over $\Q$.
By \cite[Th\'eor\`eme (5.5)(i)]{Del68}, each $R^j\pi_*\Omega^i_{X/S}|_{S_{\Q}}$ is locally free and its formation commutes with base change.
Hence each $h^{i,j}$ is locally constant on $S_{\Q}$.
Let $C\subset S$ be a connected component, with generic point $\eta$.
As $S$ is smooth over $\Z$, $C$ is integral and $\eta\in S_{\Q}$; so $C\cap S_{\Q}$ is irreducible and
\begin{equation}
\label{eq:hij-const}
h^{i,j}(t)=h^{i,j}(\eta)\quad\text{for all }t\in C\cap S_{\Q}.
\end{equation}
By upper semicontinuity, $h^{i,j}(t)\geq h^{i,j}(\eta)$ for all $t\in C$, and
\[
Z:=\bigcup_{i,j}\bigl\{t\in C\bigm| h^{i,j}(t)>h^{i,j}(\eta)\bigr\}
\]
is closed.
By \eqref{eq:hij-const}, $Z\cap S_{\Q}=\emptyset$.
The generic points of the finitely many irreducible components of $Z$ have residue characteristics $p_1,\dots,p_r>0$, so $Z\subset V(p_1\cdots p_r)$.
As $S$ is noetherian, there are only finitely many connected components $C$.
Let $N$ be the product of the $p_1\cdots p_r$ over all these finitely many $C$, replace $a$ by $Na$, and replace $\pi$ by its base change along $\Spec R[1/a]\to\Spec R$.
Then $S$ is replaced by $S[1/N]=\bigsqcup_C C'$, where
\[
C':=C\cap D(N).
\]
Each $C'$ is disjoint from $Z$, hence
\[
h^{i,j}(t)=h^{i,j}(\eta)\quad\text{for all }t\in C'.
\]
Each $C'$ is a nonempty open subset of the integral scheme $C$, hence irreducible, so the $C'$ are the connected components of the new $S$.
Thus every $h^{i,j}$ is constant on every connected component of $S$.
As $S$ is reduced, \cite[\S5, Corollary 2]{Mum74} applied on each connected component shows that each $R^j\pi_*\Omega^i_{X/S}$ is locally free and commutes with base change to the points of $S$; by \cite[Theorem (8.0)]{Kat70} its formation then commutes with arbitrary base change, which proves (\ref{item:lf-221-lfbs}), while (\ref{item:lf-smooth}) remains valid after this replacement.

As $S_{k(\mathfrak p)}$ for $\mathfrak{p}=(0)$ is a scheme of characteristic $0$, \eqref{eq:HdR-fp} for $\mathfrak{p}=(0)$ degenerates at $E_1$ by \cite[Th\'eor\`eme (5.5)(ii)]{Del68}.

We show that \eqref{eq:HdR-f} degenerates at $E_1$.
We show by induction on $r\geq1$ that
\[
E_r^{i,j}=R^j\pi_*\Omega^i_{X/S}\quad\text{for all }i,j.
\]
For $r=1$ this holds by \eqref{eq:HdR-f}.
Suppose it holds for some $r\geq1$.
Then by \hyperref[item:lf-hodge]{(LFBS)} the source and target of $d_r$ are locally free.
Since $R\to K$ is flat, $d_r\otimes_RK$ is the differential $d_r$ of \eqref{eq:HdR-fp} for $\mathfrak{p}=(0)$.
This vanishes, as \eqref{eq:HdR-fp} for $\mathfrak{p}=(0)$ degenerates at $E_1$.
As $S$ is flat over the domain $R$, every locally free $\mathcal{O}_S$-module injects into its tensor product with $K$, so $d_r=0$.
Since $E_{r+1}$ is the cohomology of $(E_r,d_r)$, it follows that $E_{r+1}^{i,j}=E_r^{i,j}=R^j\pi_*\Omega^i_{X/S}$.
This completes the induction, so \eqref{eq:HdR-f} degenerates at $E_1$.

We show that \eqref{eq:HdR-fp} degenerates at $E_1$ for $\mathfrak{p}$ a maximal ideal of $R$.
By \Cref{lem:HdR-degen}:
\begin{itemize}
\item the degeneration of \eqref{eq:HdR-f} at $E_1$ is equivalent to the exactness of the sequences \eqref{eq:Hdg-fil-seq} for $\pi$ and all $i$ and $m$;
\item the degeneration of \eqref{eq:HdR-fp} at $E_1$ is equivalent to the exactness of the sequences \eqref{eq:Hdg-fil-seq} for $\pi_{k(\mathfrak p)}$ and all $i$ and $m$.
\end{itemize}
By \hyperref[item:lf-hodge]{(LFBS)} and descending induction on $i$, all terms of \eqref{eq:Hdg-fil-seq} are locally free, so \eqref{eq:Hdg-fil-seq} is locally split.
Let $\mathfrak{p}$ be a maximal ideal of $R$ and $g\colon S_{k(\mathfrak p)}\to S$ the base change morphism.
Then $g^*$\eqref{eq:Hdg-fil-seq} is exact.
Applying \cite[Theorem (8.0)]{Kat70} to $\Omega^{\geq i}_{X/S}$ and applying \hyperref[item:lf-hodge]{(LFBS)} to $\Omega^i_{X/S}$, we have
\[
g^*R^{m}\pi_*\Omega^{\geq i}_{X/S}\simeq R^{m}\pi_{k(\mathfrak p)*}\Omega^{\geq i}_{X_{k(\mathfrak p)}/S_{k(\mathfrak p)}},\qquad g^*R^{m-i}\pi_*\Omega^i_{X/S}\simeq R^{m-i}\pi_{k(\mathfrak p)*}\Omega^i_{X_{k(\mathfrak p)}/S_{k(\mathfrak p)}},
\]
compatibly with the maps in \eqref{eq:Hdg-fil-seq}.
Hence $g^*$\eqref{eq:Hdg-fil-seq} is the sequence \eqref{eq:Hdg-fil-seq} for $\pi_{k(\mathfrak p)}$, which is therefore exact for all $i$ and $m$, so \eqref{eq:HdR-fp} degenerates at $E_1$.

For (\ref{item:lf-221-summand}), the sequences \eqref{eq:Hdg-fil-seq} for $\pi$ are exact with locally free terms, as shown above, so they are locally split.
By \Cref{lem:HdR-degen}, $F^i_{\mathrm{Hdg}}R^m\pi_*\Omega^\bullet_{X/S}=R^m\pi_*\Omega^{\geq i}_{X/S}$.
The case $i=0$ gives local freeness of $R^m\pi_*\Omega^\bullet_{X/S}$, and the local splittings give the direct summands.
The base change statement is the case $i=0$ of the first isomorphism in the proof of (\ref{item:lf-221-hdr}) above.
\end{proof}

\section{\texorpdfstring{Ordinary $K3$ surfaces and opposite flags}{Ordinary K3 surfaces and opposite flags}}
\label{sec:ordinary}

Let $k$ be a perfect field of characteristic $p>0$, and let $Y$ be a smooth proper $k$-scheme, with structure morphism $\pi\colon Y\to\Spec k$.

\begin{definition}[$F$-zip]
\label{def:fzip}
\begin{enumerate}
\item
Let $\sigma$ be the Frobenius endomorphism of $k$.
An \emph{$F$-zip} over $k$ \cite[Definition 1.7]{Wed08} is a tuple
\[
(M,C^\bullet,D_\bullet,\varphi_\bullet)
\]
as follows.
\begin{enumerate}
\item
$M$ is a finite-dimensional $k$-vector space.
\item
$C^\bullet$ is a descending filtration of $\sigma^*M$, by which we mean $k$-subspaces $C^i\subset\sigma^*M$ with $C^i\supset C^{i+1}$, equal to $\sigma^*M$ for $i\ll 0$ and equal to $0$ for $i\gg 0$.
Write $\mathrm{gr}^i_C=C^i/C^{i+1}$.
\item
$D_\bullet$ is an ascending filtration of $M$, by which we mean $k$-subspaces $D_i\subset M$ with $D_i\subset D_{i+1}$, equal to $0$ for $i\ll 0$ and equal to $M$ for $i\gg 0$.
Write $\mathrm{gr}_i D=D_i/D_{i-1}$.
\item
For each $i\in\Z$, $\varphi_i\colon\mathrm{gr}^i_C\xrightarrow{\sim}\mathrm{gr}_i D$ is a $k$-linear isomorphism.
\end{enumerate}
\item
\label{def:family-fzip}
Let $S$ be a scheme.
For a prime $p$, write $S_{\F_p}=S\times_{\Spec\Z}\Spec\F_p$, and $\sigma_p=F_{\mathrm{abs}}$ for the absolute Frobenius of $S_{\F_p}$.
A \emph{family of $F$-zips} over $S$, or simply an \emph{$F$-zip over $S$}, is a tuple
\[
\bigl(M,(C_p^\bullet,D_{p,\bullet},\varphi_{p,\bullet})_p\bigr)
\]
as follows.
\begin{enumerate}
\item
$M$ is a locally free $\cO_S$-module of finite rank.
\item
For each prime $p$, the tuple $(M|_{S_{\F_p}},C^\bullet_p,D_{p,\bullet},\varphi_{p,\bullet})$ is an $F$-zip over $S_{\F_p}$ \cite[(1.11), Definition~1.13]{Wed08}.
That is, $C_p^\bullet$ and $D_{p,\bullet}$ are filtrations of $\sigma_p^*(M|_{S_{\F_p}})$ and $M|_{S_{\F_p}}$ by locally direct summands, and the $\varphi_{p,i}\colon\mathrm{gr}^i_{C_p}\xrightarrow{\sim}\mathrm{gr}_iD_p$ are $\cO_{S_{\F_p}}$-linear isomorphisms.
\end{enumerate}
\end{enumerate}
\end{definition}

If $S$ is an $\F_p$-scheme, then $S_{\F_p}=S$ and $S_{\F_\ell}=\emptyset$ for $\ell\neq p$.
So a family of $F$-zips over $S$ is the same as an $F$-zip over $S$ in the sense of \cite[(1.11) and Definition~1.13]{Wed08}.

\begin{definition}[Ordinary $F$-zip]
\label{def:ordinary-fzip}
An $F$-zip $(M,C^\bullet,D_\bullet,\varphi_\bullet)$ over $k$ \cite[Definition 1.12]{Wed08} is \emph{ordinary} if
\begin{equation}
\label{eq:ordinary-fzip}
\sigma^*(D_{i-1})\oplus C^i=\sigma^*(M)
\end{equation}
for all $i\in\Z$.
\end{definition}

\begin{example}
\label{ex:ordinary-scheme}
\begin{enumerate}
\item
Suppose that $Y$ satisfies Hodge--de Rham degeneration at $E_1$.
For every $m$, the $F$-zip structure on $M=R^{m}\pi_*\Omega^\bullet_{Y/k}$ is given by
\[
C^i=\sigma^*\bigl(F_{\mathrm{Hdg}}^iR^{m}\pi_*\Omega^\bullet_{Y/k}\bigr)\quad\text{and}\quad D_i=F^{m-i}_{\conj}R^{m}\pi_*\Omega^\bullet_{Y/k},
\]
and the Cartier isomorphisms
\[
\varphi_i\colon C^i/C^{i+1}=\sigma^*\mathrm{gr}^i_{F_{\mathrm{Hdg}}}M\xrightarrow{\sim}\mathrm{gr}^{m-i}_{F_{\conj}}M=D_i/D_{i-1}.
\]
See \cite[(1.7)(b)(c)(d)]{Wed08}.
We say that $Y$ is \emph{ordinary} if all these de Rham $F$-zips are ordinary \cite[Definition 1.12]{Wed08}.
\item
More generally, let $S$ be a scheme, and let $f\colon X\to S$ be smooth and proper, such that its Hodge--de Rham spectral sequence degenerates at $E_1$ and all $R^bf_*\Omega^a_{X/S}$ are locally free.
\begin{enumerate}
\item
Then $M=R^mf_*\Omega^\bullet_{X/S}$ is locally free.
\item
By cohomology and base change, the formation of the $R^bf_*\Omega^a_{X/S}$ and of $M$ commutes with base change.
So for each prime $p$, the base change $f_p\colon X_{\F_p}\to S_{\F_p}$ satisfies the same conditions, and $M|_{S_{\F_p}}=R^mf_{p*}\Omega^\bullet_{X_{\F_p}/S_{\F_p}}$.
\item
Take $C_p^\bullet$, $D_{p,\bullet}$ and $\varphi_{p,\bullet}$ to be the $F$-zip structure of $f_p$ over $S_{\F_p}$ \cite[(1.11)]{Wed08}, given by the same formulas as in~(1).
\end{enumerate}
This makes $M$ a family of $F$-zips over $S$ in the sense of \Cref{def:fzip}\ref{def:family-fzip}.
\end{enumerate}
\end{example}

\begin{example}
\label{ex:K3-fzip}
If $Y$ is a $K3$ surface, then its Hodge--de Rham spectral sequence degenerates at $E_1$.
Over $\bar k$, this is \cite[Proposition 1.1(a)]{Del81}.
It descends to $k$ by \Cref{lem:HdR-degen}(3), since $k\to\bar k$ is faithfully flat and the sequences \eqref{eq:Hdg-fil-seq} are compatible with flat base change.
See also \cite[(2.6)]{Wed08}.
\begin{enumerate}
\item Suppose that $Y$ is a $K3$ surface, and let $m=2$. The Hodge numbers
\[
h^{2,0}=h^{0,2}=1\quad\text{and}\quad h^{1,1}=20
\]
give \cite[(2.6)]{Wed08}
\begin{align*}
\dim C^1&=21,&\dim C^2&=1,\\
\dim D_0&=1,&\dim D_1&=21.
\end{align*}
Thus $C^2=\sigma^*\bigl(F_{\mathrm{Hdg}}^2(R^2\pi_*\Omega^\bullet_{Y/k})\bigr)$ and $D_1=F^1_{\conj}\bigl(R^2\pi_*\Omega^\bullet_{Y/k}\bigr)$.
\item \label{ex:K3-family}Similarly, let $\pi\colon X\to S$ be a smooth projective morphism whose fibers are $K3$ surfaces, satisfying \hyperref[item:lf-hodge]{(LFBS)} and \hyperref[item:lf-degen]{(HdR)}, and let $m=2$.
By \Cref{ex:ordinary-scheme}, $M=R^2\pi_*\Omega^\bullet_{X/S}$ is an $F$-zip over $S$.
The $R^j\pi_*\Omega^i_{X/S}$ with $i+j=2$ are locally free, of ranks
\[
\rk_{\cO_S}\pi_*\Omega^2_{X/S}=\rk_{\cO_S} R^2\pi_*\cO_X=1\quad\text{and}\quad\rk_{\cO_S} R^1\pi_*\Omega^1_{X/S}=20.
\]
These give, for every prime $p$,
\begin{align*}
\rk_{\cO_{S_{\F_p}}} C_p^1&=21,&\rk_{\cO_{S_{\F_p}}} C_p^2&=1,\\
\rk_{\cO_{S_{\F_p}}} D_{p,0}&=1,&\rk_{\cO_{S_{\F_p}}} D_{p,1}&=21.
\end{align*}
The $C_p^i$ and $D_{p,i}$ are locally direct summands \cite[(1.11), Definition~1.13]{Wed08}, so they and their graded pieces are locally free.
Thus $C_p^2=\sigma_p^*\bigl(F_{\mathrm{Hdg}}^2M|_{S_{\F_p}}\bigr)$ with $F^2_{\mathrm{Hdg}}M\simeq\pi_*\Omega^2_{X/S}$ a line bundle, and $\varphi_{p,2}$ is the Cartier isomorphism
\[
C_p^2\xrightarrow{\sim}\mathrm{gr}^0_{F_{\conj}}(M|_{S_{\F_p}}).
\]
So $\mathrm{gr}^0_{F_{\conj}}(M|_{S_{\F_p}})$ is a line bundle on $S_{\F_p}$, isomorphic to $\sigma_p^*\bigl(\pi_*\Omega^2_{X/S}|_{S_{\F_p}}\bigr)$.
\end{enumerate}
\end{example}

\begin{remark}
\label{rmk:ordinary-K3}
Assume that $Y$ is a $K3$ surface, and let $h$ be the height of its formal Brauer group \cite[pp.~90--91]{AM77}.
By \cite[II, 7.2]{Ill79} and \cite[(2.6)]{Wed08}, $Y$ is ordinary in the sense of \Cref{ex:ordinary-scheme} if and only if $h=1$.
For finite $k$, this is the definition of ordinary in \cite[\S1.2]{BZ09}.
\end{remark}

\begin{lemma}
\label{lem:opposite}
Suppose that $Y$ satisfies Hodge--de Rham degeneration at $E_1$, and let $m\geq0$ be such that $R^{m}\pi_*\Omega^\bullet_{Y/k}$, with the $F$-zip structure of \Cref{ex:ordinary-scheme}, is ordinary.
Then for every $i\in\Z$,
\[
F_{\mathrm{Hdg}}^{i}\bigl(R^{m}\pi_{*}\Omega^{\bullet}_{Y/k}\bigr)\cap F^{m-i+1}_{\conj}\bigl(R^{m}\pi_{*}\Omega^{\bullet}_{Y/k}\bigr)=0.
\]
\end{lemma}

\begin{proof}
Write $M=R^{m}\pi_*\Omega^\bullet_{Y/k}$ and $(M,C^\bullet,D_\bullet,\varphi_\bullet)$ for its $F$-zip.
By \cite[(1.7)(b)(c)]{Wed08}, as in \Cref{ex:ordinary-scheme},
\[
C^i=\sigma^*\bigl(F_{\mathrm{Hdg}}^i M\bigr)\quad\text{and}\quad D_{i-1}=F^{m-i+1}_{\conj}M.
\]
Since the $F$-zip is ordinary, \eqref{eq:ordinary-fzip} gives
\[
\sigma^*(D_{i-1})\oplus C^i=\sigma^*(M),
\]
so $\sigma^*\bigl(F^{m-i+1}_{\conj}M\bigr)\cap\sigma^*\bigl(F_{\mathrm{Hdg}}^iM\bigr)=0$ in $\sigma^*(M)$.
For $k$ perfect, $\sigma$ is an automorphism of $k$, so $\sigma^*$ is an autoequivalence of finite-dimensional $k$-vector spaces and commutes with intersections of subspaces.
Hence $F_{\mathrm{Hdg}}^iM\cap F^{m-i+1}_{\conj}M=0$.
\end{proof}

In particular, if $Y$ satisfies Hodge--de Rham degeneration at $E_1$ and $R^2\pi_*\Omega^\bullet_{Y/k}$ is ordinary, then
\begin{align}
\label{eq:opposite-n2-i2}
F_{\mathrm{Hdg}}^2\bigl(R^2\pi_*\Omega^\bullet_{Y/k}\bigr)\cap F^1_{\conj}\bigl(R^2\pi_*\Omega^\bullet_{Y/k}\bigr)&=0,\\
\label{eq:opposite-n2-i1}
F_{\mathrm{Hdg}}^1\bigl(R^2\pi_*\Omega^\bullet_{Y/k}\bigr)\cap F^2_{\conj}\bigl(R^2\pi_*\Omega^\bullet_{Y/k}\bigr)&=0.
\end{align}
If $Y$ is ordinary in the sense of \Cref{ex:ordinary-scheme}, then \Cref{lem:opposite} applies to every $m\geq0$.

\begin{definition}[Density one]
\label{def:density-one}
Let $L$ be a number field, and let $\Sigma$ be a set of non-archimedean places of $L$.
Write $N(w):=\#k(w)$ for the cardinality of the residue field $k(w)$ of a place $w$.
We say that $\Sigma$ has \emph{density one} if
\[
\lim_{s\to 1^+}\frac{\sum_{w\in\Sigma} N(w)^{-s}}{\sum_{w} N(w)^{-s}}=1,
\]
where both sums run over the non-archimedean places of $L$.
\end{definition}

\begin{theorem}[{\cite[Theorem 0.1]{BZ09}}]
\label{thm:bz}
Let $X$ be a $K3$ surface over a number field $K$.
Then there exists a finite extension $L/K$ and a set of non-archimedean places of $L$ of density one at which $X_L$ has ordinary good reduction.
\end{theorem}

Ordinary here is height one of the formal Brauer group of the reduction, equivalently Newton slopes $0,1,2$, equivalently an ordinary $F$-zip on $R^2\pi_*\Omega^\bullet$ of the reduction, as above.

\section{The theorem}
\label{sec:theorem}

\begin{theorem}
\label{thm:main}
Let $K$ be a number field and $R\subset K$ a finitely generated $\Z$-algebra.
Let $S$ be a smooth integral $R$-scheme of finite type, and let $\pi\colon X\to S$ be a smooth projective morphism whose fibers are $K3$ surfaces.
Let
\[
(\cE,\nabla)=\bigl(R^2\pi_*\Omega^\bullet_{X/S},\nabla_{\GM}\bigr).
\]
Assume:
\begin{enumerate}[label=\textup{(\alph*)}]
\item $\theta\colon F_{\mathrm{Hdg}}^2\cE\to\bigl(F_{\mathrm{Hdg}}^1\cE/F_{\mathrm{Hdg}}^2\cE\bigr)\otimes_{\cO_S}\Omega^1_{S/R}$ is not identically zero.
\item $s\in S(R)$ and $v\colon\Spec R\to E$ is a morphism over $s$ that factors through $F_{\mathrm{Hdg}}^2E$ and is not the zero section.
\end{enumerate}
Then the formal horizontal section of $(\cE,\nabla)_K$ through $v_K$ is not $\omega(p)$-integral.
\end{theorem}

\begin{proof}
Suppose that the formal horizontal section $\tilde{v}_K$ through $v_K$ is $\omega(p)$-integral, and fix
\[
g\colon\{\mathrm{primes}\}\to\Z_{\geq 0}\quad\text{with}\quad g(p)/p\to\infty
\]
such that $\tilde{v}_K$ is $g$-integral with respect to the neighborhood $U_s$, the coordinates $x_1,\dots,x_d$ and the basis $e_1,\dots,e_r$ of \Cref{def:leaf-notions}.
The localizations in Step~0 do not affect this assumption:
\begin{itemize}
\item If $R'\subset K$ is a localization of $R$, then $R_n\subset R'_n$ (notation defined in \eqref{eq:Rn}) for every $n$, so $\tilde{v}_K$ remains $g$-integral over $R'$ by \eqref{eq:g-series}.
\item Replacing $S$ by an open subscheme containing $s(\Spec R)$ does not change $\tilde{v}_K$.
\end{itemize}

\textbf{Step 0: localizations.}
\begin{enumerate}[label=\textup{(\arabic*)}]
\item Replace $R$ by $R[1/a]$ as in \Cref{rmk:local-freeness}.
\end{enumerate}
By \Cref{rmk:local-freeness}(\ref{item:lf-221-lfbs}) and (\ref{item:lf-221-hdr}), $\pi$ satisfies \hyperref[item:lf-hodge]{(LFBS)}, and $\pi$ and its reductions $\pi_{k(\mathfrak p)}$ satisfy \hyperref[item:lf-degen]{(HdR)}.
In particular, \Cref{prop:katz-F2} applies to the reduction of $\pi$ modulo every maximal ideal of $R$.
By \Cref{rmk:local-freeness}(\ref{item:lf-221-summand}), $\cE$ is locally free and the Hodge filtration consists of local direct summands.
By \Cref{ex:K3-fzip}\ref{ex:K3-family}, $F^2_{\mathrm{Hdg}}\cE\simeq\pi_*\Omega^2_{X/S}$ is a line bundle, so $s^*F_{\mathrm{Hdg}}^2\cE$ is an invertible sheaf on $\Spec R$ with global section $v$.
As $v$ is not the zero section and $R$ is a domain, $v$ generates $s^*F_{\mathrm{Hdg}}^2\cE$ over $\Spec R[1/b]$ for some $0\neq b\in R$.
\begin{enumerate}[label=\textup{(\arabic*)},resume]
\item Replace $R$ by $R[1/b]$.
\end{enumerate}
In particular, the reduction of $v$ remains nonzero modulo every maximal ideal of $R$.

Since $s$ is a section, $S\to\Spec R$ is dominant.
The generic point of the integral scheme $S$ therefore lies on $S_K$, so $S_K$ is irreducible.
The zero locus of $\theta$ is closed, and $\theta$ is not identically zero, so $\theta$ remains not identically zero on any open neighborhood of $s(\Spec R)$.
By \Cref{rmk:local-freeness}(\ref{item:lf-221-summand}), the sheaves $\cE$, $F_{\mathrm{Hdg}}^i\cE$ and $\mathrm{gr}^i_{F_{\mathrm{Hdg}}}\cE$ for all $i$, and $\Omega^1_{S/R}$ are locally free.
Let $U\subset U_s$ (defined in \Cref{def:leaf-notions}) be an affine open neighborhood of the image of the generic point of $\Spec R$ on which all of them are free.
Then $s^{-1}(U)$ contains $\Spec R[1/c]$ for some $0\neq c\in R$.
Write $\cI_s$ for the ideal sheaf of $s$, as in \eqref{eq:Is}.
\begin{enumerate}[label=\textup{(\arabic*)},resume]
\item Replace $R$ by $R[1/c]$ and $S$ by $U\otimes_R R[1/c]$.
\end{enumerate}
Then $s\in S(R)$, $S$ is affine, and $S$ is contained in the original $U_s$.
So $x_1,\dots,x_d$ and $e_1,\dots,e_r$ are defined on all of $S$, and we may take $U_s=S$ in \Cref{def:leaf-notions}.
On $S$, $x_1,\dots,x_d$ generate $\cI_s$.
By \Cref{def:leaf-notions}, $\hat{\cO}_{S,s}\simeq R[[x_1,\dots,x_d]]$.
After these replacements, the following hold.
\begin{enumerate}[label=\textup{(0.\alph*)},ref=\textup{(0.\alph*)},widest=d,align=left,labelindent=0pt]
\item\label{item:step0-a} $\pi$ satisfies \hyperref[item:lf-hodge]{(LFBS)}, and $\pi$ and every $\pi_{k(\mathfrak p)}$ satisfy \hyperref[item:lf-degen]{(HdR)}. In particular, \Cref{prop:katz-F2} applies to the reduction of $\pi$ modulo every maximal ideal of $R$.
\item\label{item:step0-b} $v$ generates $s^*F_{\mathrm{Hdg}}^2\cE$. In particular, the reduction of $v$ modulo every maximal ideal of $R$ is nonzero.
\item\label{item:step0-c} $S$ is affine, and $U_s=S$.
\item\label{item:step0-d} $\cE$, $F_{\mathrm{Hdg}}^i\cE$ and $\mathrm{gr}^i_{F_{\mathrm{Hdg}}}\cE$ for all $i$, and $\Omega^1_{S/R}$ are free on $S$.
\end{enumerate}
Fix the following on $S$:
\begin{itemize}
\item a generator $e$ of $F_{\mathrm{Hdg}}^2\cE$ with $s^*e=v$;
\item a basis $h_1,\dots,h_{20}$ of $F_{\mathrm{Hdg}}^1\cE/F_{\mathrm{Hdg}}^2\cE$;
\item a basis $\omega_1,\dots,\omega_d$ of $\Omega^1_{S/R}$.
\end{itemize}
Then we can write the map $\theta\colon F_{\mathrm{Hdg}}^2\cE\to\bigl(F_{\mathrm{Hdg}}^1\cE/F_{\mathrm{Hdg}}^2\cE\bigr)\otimes_{\cO_S}\Omega^1_{S/R}$ as
\begin{equation}\label{eq:theta-kl}
{e\mapsto}\sum_{k=1}^{20}\sum_{l=1}^{d}\theta_{kl}\,{h_k}\otimes\omega_l,\qquad{\text{with each }\theta_{kl}\in\Gamma(S,\cO_S)}.
\end{equation}
The $\theta_{kl}$ are the coefficient functions of $\theta$.
Via $\Gamma(S,\cO_S)\to\hat{\cO}_{S,s}\simeq R[[x_1,\dots,x_d]]$, we view each $\theta_{kl}$ as a power series in $x_1,\dots,x_d$ with coefficients in $R$.

\textbf{Step 1: the good primes.}
Let $X_s$ denote the fiber of $\pi$ over $s$, with generic fiber $X_{s,K}$.
By \Cref{thm:bz} there is a finite extension $L/K$ and a set
\[
\Sigma:=\{\text{non-archimedean places }w\text{ of }L\mid X_{s,L}\text{ has ordinary good reduction at }w\}
\]
of density one, where $\cO_w$ and $k(w)$ denote the valuation ring and the residue field of $w$.
Following \cite[\S2.1 and \S3.1]{BZ09}, by ordinary good reduction at $w$ we mean that there is a smooth projective model of $X_{s,L}$ over $\cO_w$ whose special fiber is an ordinary $K3$ surface over $k(w)$ in the sense of \cite[\S1.2]{BZ09}.
In particular, $\Sigma$ is infinite.
\begin{claim}\label{claim:good-primes}
For all but finitely many places $w\in\Sigma$, with maximal ideal $\mathfrak{m}_w$ of $\cO_w$, the ideal $\mathfrak{p}:=\mathfrak{m}_w\cap R$ is a maximal ideal of $R$, and the reduction $X_{s,k(\mathfrak{p})}$ is an ordinary $K3$ surface over the finite field $k(\mathfrak{p})$.
\end{claim}
In fact, for all but finitely many $w\in\Sigma$, we have $R\subset\cO_w$.
For such $w$, the ideal $\mathfrak{p}=\mathfrak{m}_w\cap R$ is a maximal ideal of $R$.
The scheme $X_s\otimes_R\cO_w$ is a smooth projective model of $X_{s,L}$ over $\cO_w$.
By \cite[Lemma~1.1 and \S3.1]{BZ09} ordinarity of the special fiber depends only on the characteristic polynomial of geometric Frobenius, which is the same for every smooth projective model.
Thus the reduction $X_{s,k(\mathfrak{p})}$ is an ordinary $K3$ surface over the finite field $k(\mathfrak{p})$.
This finishes the proof of \Cref{claim:good-primes}.

Let $\Sigma'\subset\Sigma$ be the subset of places $w$ for which the conclusion of \Cref{claim:good-primes} holds, and let
\[
\mathcal{P}:=\{\mathfrak{m}_w\cap R\mid w\in\Sigma'\}.
\]
Only finitely many places of $L$ contract to a given maximal ideal of $R$: they are the maximal ideals of the integral closure of $R_{\mathfrak{p}}$ in $L$.
Therefore $\mathcal{P}$ is infinite.
Since $\dim R=1$, only finitely many maximal ideals of $R$ have a given residue characteristic, so the residue characteristics of this set are unbounded.

\textbf{Step 2: vanishing of $\theta$ to order $\lfloor g(p)/p\rfloor-1$ modulo $\mathfrak p$.}
Let $\mathfrak{p}\in\mathcal{P}$ be of residue characteristic $p$ with $g(p)>p$.
We use the notation of \ref{para:reduction} for $\mathfrak p$.
By \ref{item:step0-a}, the reduction of $\pi$ modulo $\mathfrak{p}$ satisfies the hypotheses of \Cref{prop:katz-F2}.
We record the following facts on $\pi_{k(\mathfrak p)}$.
\begin{itemize}
\item By \Cref{rmk:local-freeness}(\ref{item:lf-221-summand}), $\cE_{k(\mathfrak p)}$ is $R^2\pi_{k(\mathfrak p)*}\Omega^\bullet_{X_{k(\mathfrak p)}/S_{k(\mathfrak p)}}$.
\item By \Cref{lem:conj-bc}(2) for $\pi_{k(\mathfrak p)}$ with $m=2$ and $a=0$,
\[
\mathrm{gr}^0_{F_{\conj}}\cE_{k(\mathfrak p)}\cong F_{\mathrm{abs}}^*\pi_{k(\mathfrak p)*}\Omega^2_{X_{k(\mathfrak p)}/S_{k(\mathfrak p)}}.
\]
This is free of rank $1$ on $S_{k(\mathfrak p)}$, since $\pi_{k(\mathfrak p)*}\Omega^2_{X_{k(\mathfrak p)}/S_{k(\mathfrak p)}}$ is the restriction of $\pi_*\Omega^2_{X/S}\simeq F^2_{\mathrm{Hdg}}\cE$ by \hyperref[item:lf-hodge]{(LFBS)}, which is free of rank $1$ by \ref{item:step0-d}.
\item By \Cref{lem:conj-bc}(3) for $\pi_{k(\mathfrak p)}$, the fiber of $\mathrm{gr}^0_{F_{\conj}}\cE_{k(\mathfrak p)}$ at $s_{k(\mathfrak p)}$ is $\mathrm{gr}^0_{F_{\conj}}\bigl(R^2\pi_*\Omega^\bullet_{X_{s,k(\mathfrak{p})}/k(\mathfrak{p})}\bigr)$.
\item By \eqref{eq:psi-conj}, $\psi_p$ carries $F^i_{\conj}$ into $F^{i+1}_{\conj}$, and $\mathrm{gr}^0_{F_{\conj}}\psi_p$ is well-defined as an $\cO$-linear map
\[
\mathrm{gr}^0_{F_{\conj}}\psi_p\colon\mathrm{gr}^0_{F_{\conj}}\cE_{k(\mathfrak p)}\to F_{\mathrm{abs}}^*\Omega^1_{S_{k(\mathfrak p)}/k(\mathfrak p)}\otimes_{\cO_{S_{k(\mathfrak p)}}}\mathrm{gr}^1_{F_{\conj}}\cE_{k(\mathfrak p)}.
\]
Equivalently, for every local vector field $\partial$ on $S_{k(\mathfrak p)}$, it gives
\begin{equation}\label{eq:gr0psi-partial}
(\mathrm{gr}^0_{F_{\conj}}\psi_p)(\partial)\colon\mathrm{gr}^0_{F_{\conj}}\cE_{k(\mathfrak p)}\to\mathrm{gr}^1_{F_{\conj}}\cE_{k(\mathfrak p)}.
\end{equation}
\end{itemize}

Since $U_s=S$ by \ref{item:step0-c}, we may take $U'=S$ in \ref{para:reduction}.
Then $\overline{\tilde v}_{R_{\mathfrak p}}\in\Gamma(S_{k(\mathfrak p)},\cE_{k(\mathfrak p)})$ from \eqref{eq:vRp-bar} is a lift of $\tilde v_{k(\mathfrak{p})}$ to a section of $\cE_{k(\mathfrak p)}$ on all of $S_{k(\mathfrak p)}$; see \Cref{fig:step2}.
\begin{center}\begin{minipage}{\textwidth}
\centering
\newcommand{\fu}[4]{\shade[top color=#1,bottom color=white] ($(#2)+(-#3,-#4)$) rectangle ($(#2)+(#3,0)$);}
\newcommand{\fv}[4]{\shade[bottom color=#1,top color=white] ($(#2)+(-#3,0)$) rectangle ($(#2)+(#3,#4)$);}
\makebox[\linewidth]{\resizebox{\textwidth}{!}{%
\begin{tikzpicture}[>=Stealth,line cap=round,line join=round]
\def\dx{1.6}\def\dy{1.8}\def\H{3.4}
\draw[thick] (0,0) -- (8,0) node[right]{$\Spec R$};
\foreach \x in {1,5}{\fill (\x,0) circle (1.3pt);}
\fill (3,0) circle (1.8pt) node[below]{\small $\mathfrak p$};
\draw[thick,fill=white] (6.6,0) ellipse [x radius=5pt, y radius=2.6pt] node[below=2pt]{\small $\eta=\Spec K$};
\draw[orange!80!black,decorate,decoration={brace,mirror,amplitude=4pt}] (2.8,-0.55) -- (6.9,-0.55) node[midway,below=4pt]{\small $\Spec R_{\mathfrak p}=\{\mathfrak p,\eta\}$};
\begin{scope}[yshift=1.4cm]
\fill[gray!10] (0,0) -- (8,0) -- (8+\dx,\dy) -- (\dx,\dy) -- cycle;
\begin{scope}
\clip (0,0) -- (8,0) -- (8+\dx,\dy) -- (\dx,\dy) -- cycle;
\def\Up{(3.36,0.405) .. controls (4.5,-0.2) and (7.5,-0.2) .. (8.6,0.6) .. controls (9.4,1.2) and (8.8,1.75) .. (7.2,1.6) .. controls (5.5,1.45) and (4.8,1.6) .. (4.2304,1.3842) .. controls (3.3,1.05) and (2.85,0.75) .. (3.36,0.405) -- cycle}
\fill[orange!15] \Up; \draw[orange!80!black,dashed] \Up;
\foreach \w in {0.04,0.08,...,0.36}{\fill[red!80!black,opacity=0.07] (6.6-\w,0) -- (6.6+\w,0) -- (6.6+\dx+\w,\dy) -- (6.6+\dx-\w,\dy) -- cycle;}
\end{scope}
\draw[thick] (0,0) -- (8,0) -- (8+\dx,\dy) -- (\dx,\dy) -- cycle;
\foreach \x in {1,5}{\draw[gray!60] (\x,0) -- (\x+\dx,\dy);}
\draw[gray!70!black,thick] (3,0) -- (3+\dx,\dy);
\node[gray!60!black] at (3,-0.05) [below] {\scriptsize $S_{k(\mathfrak p)}$};
\node[red!70!black] at (6.6,-0.05) [below]{\scriptsize $S_K$};
\node[orange!80!black] at (8.7,0.35) {\small $U'$};
\draw[orange!80!black,line width=3pt,opacity=0.6] (3+0.2*\dx,0.2*\dy) -- (3+0.75*\dx,0.75*\dy);
\def\sc{(0.2,0.35) .. controls (1.5,0.9) and (2.3,0.65) .. (3+0.45*\dx,0.45*\dy) .. controls (4.6,1.0) and (5.6,0.6) .. (6.6+0.45*\dx,0.45*\dy) .. controls (7.8,1.0) .. (8.6,0.95)}
\draw[blue!60!black,line width=5pt,opacity=0.18] \sc;
\draw[blue!60!black,very thick] \sc node[above right]{$s$};
\fill[blue!60!black,opacity=0.35] (3+0.45*\dx,0.45*\dy) ellipse [x radius=0.22, y radius=0.34, rotate=-41];
\fill[blue!60!black] (3+0.45*\dx,0.45*\dy) circle (1.8pt) node[below right=1pt]{\scriptsize $s_{k(\mathfrak p)}$};
\draw[red!70!black,thick,dotted] (6.6+0.45*\dx,0.45*\dy) circle (0.33);
\fill[blue!60!black] (6.6+0.45*\dx,0.45*\dy) circle (2pt) node[below right=2pt]{\scriptsize $s_K$};
\node at (0.3,1.25) {$S$};
\end{scope}
\draw[->] (4,1.3) -- (4,0.25);
\begin{scope}[yshift=5.6cm]
\def\D{1.0}
\fill[gray!12,path fading=south] (0,-\D) -- (8,-\D) -- (8+\dx,\dy-\D) -- (8+\dx,\dy) -- (8,0) -- (0,0) -- cycle;
\fill[gray!12,path fading=north] (0,0) -- (8,0) -- (8+\dx,\dy) -- (8+\dx,\dy+\H) -- (\dx,\dy+\H) -- (0,\H) -- cycle;
\fv{gray!60}{\dx,\dy}{0.2pt}{\H} \fu{gray!60}{\dx,\dy}{0.2pt}{\D}
\fill[gray,opacity=0.3,path fading=north] (3,0) -- (3+\dx,\dy) -- (3+\dx,\dy+\H) -- (3,\H) -- cycle;
\fill[gray,opacity=0.3,path fading=south] (3,-\D) -- (3+\dx,\dy-\D) -- (3+\dx,\dy) -- (3,0) -- cycle;
\fv{gray!70!black}{3,0}{0.25pt}{\H} \fv{gray!70!black}{3+\dx,\dy}{0.25pt}{\H} \fu{gray!70!black}{3,0}{0.25pt}{\D} \fu{gray!70!black}{3+\dx,\dy}{0.25pt}{\D}
\draw[gray!60!black,dashed] (3+0.45*\dx,0.45*\dy-\D) -- (3+0.45*\dx,0.45*\dy+\H-0.6);
\fill[black!35] (3,0) -- (3+\dx,\dy) -- (3+\dx,\dy+0.77) .. controls (3+0.65*\dx,0.65*\dy+0.42) and (3+0.3*\dx,0.3*\dy+0.95) .. (3,0.56) -- cycle;
\fill[black!35,path fading=south] (3,0) -- (3+\dx,\dy) -- (3+\dx,\dy-0.77) .. controls (3+0.65*\dx,0.65*\dy-0.42) and (3+0.3*\dx,0.3*\dy-0.95) .. (3,-0.56) -- cycle;
\draw[black!50] (3,0.56) .. controls (3+0.3*\dx,0.3*\dy+0.95) and (3+0.65*\dx,0.65*\dy+0.42) .. (3+\dx,\dy+0.77);
\draw[black!50,path fading=south] (3,-0.56) .. controls (3+0.3*\dx,0.3*\dy-0.95) and (3+0.65*\dx,0.65*\dy-0.42) .. (3+\dx,\dy-0.77);
\fill[black!60] (3,0) -- (3+\dx,\dy) -- (3+\dx,\dy+0.25) .. controls (3+0.65*\dx,0.65*\dy+0.38) and (3+0.3*\dx,0.3*\dy+0.1) .. (3,0.21) -- cycle;
\fill[black!60,path fading=south] (3,0) -- (3+\dx,\dy) -- (3+\dx,\dy-0.25) .. controls (3+0.65*\dx,0.65*\dy-0.38) and (3+0.3*\dx,0.3*\dy-0.1) .. (3,-0.21) -- cycle;
\draw[black!75] (3,0.21) .. controls (3+0.3*\dx,0.3*\dy+0.1) and (3+0.65*\dx,0.65*\dy+0.38) .. (3+\dx,\dy+0.25);
\draw[black!75,path fading=south] (3,-0.21) .. controls (3+0.3*\dx,0.3*\dy-0.1) and (3+0.65*\dx,0.65*\dy-0.38) .. (3+\dx,\dy-0.25);
\foreach \w in {0.04,0.08,...,0.36}{\fill[red!80!black,opacity=0.05,path fading=north] (6.6-\w,0) -- (6.6+\w,0) -- (6.6+\dx+\w,\dy) -- (6.6+\dx+\w,\dy+\H) -- (6.6+\dx-\w,\dy+\H) -- (6.6-\w,\H) -- cycle; \fill[red!80!black,opacity=0.05,path fading=south] (6.6-\w,-\D) -- (6.6+\w,-\D) -- (6.6+\dx+\w,\dy-\D) -- (6.6+\dx+\w,\dy) -- (6.6+\dx-\w,\dy) -- (6.6-\w,0) -- cycle;}
\draw[orange!80!black,very thick] (3.360,1.566) .. controls (3.536,1.879) and (3.633,1.884) .. (3.720,2.060) .. controls (3.851,2.290) and (4.036,2.311) .. (4.231,2.543);
\fill[orange!60,opacity=0.55] (6.6-0.3+0.12*\dx,0.12*\dy+1.25) .. controls (6.6+0.3*\dx,0.3*\dy+1.55) .. (6.6+0.3+0.85*\dx,0.85*\dy+1.35) -- (6.6+0.3+0.85*\dx,0.85*\dy+1.15) .. controls (6.6+0.3*\dx,0.3*\dy+1.35) .. (6.6-0.3+0.12*\dx,0.12*\dy+1.05) -- cycle;
\draw[orange!80!black,very thick,dashed] (3,0.9) .. controls (3+0.06*\dx,0.06*\dy+0.95) .. (3+0.12*\dx,0.12*\dy+1.0) .. controls (3.255,1.364) and (3.311,1.477) .. (3.360,1.566);
\draw[orange!80!black,very thick,dashed] (4.231,2.543) .. controls (4.274,2.593) and (4.317,2.654) .. (4.360,2.730) .. controls (3+0.92*\dx,0.92*\dy+1.25) .. (3+\dx,\dy+1.15);
\def\vc{(0.2,0.35+1.0) .. controls (1.5,0.9+1.5) and (2.3,0.65+0.9) .. (3+0.45*\dx,0.45*\dy+1.25) .. controls (4.6,1.0+1.6) and (5.6,0.6+0.9) .. (6.6+0.45*\dx,0.45*\dy+1.3) .. controls (7.8,1.0+1.5) .. (8.6,0.95+1.2)}
\draw[blue!60!black,very thick] \vc node[right]{$v$};
\draw[blue!80!black,line width=3.2pt,opacity=0.75] (3+0.37*\dx,0.37*\dy+1.17) .. controls (3+0.45*\dx,0.45*\dy+1.3) .. (3+0.53*\dx,0.53*\dy+1.3);
\fill[blue!60!black] (3+0.45*\dx,0.45*\dy+1.25) circle (1.8pt);
\fill[red!70!black,opacity=0.55] (6.6+0.45*\dx,0.45*\dy+1.3) ellipse [x radius=0.36, y radius=0.16, rotate=12];
\fill[blue!60!black] (6.6+0.45*\dx,0.45*\dy+1.3) circle (2pt) node[below right=2pt,blue!60!black]{\scriptsize $v_K$};
\fv{black}{0,0}{0.4pt}{\H} \fv{black}{8,0}{0.4pt}{\H} \fv{black}{8+\dx,\dy}{0.4pt}{\H} \fu{black}{0,0}{0.4pt}{\D} \fu{black}{8,0}{0.4pt}{\D} \fu{black}{8+\dx,\dy}{0.4pt}{\D}
\draw[thick] (0,0) -- (8,0) -- (8+\dx,\dy); \draw[thin,dashed,gray!60!black] (0,0) -- (\dx,\dy) -- (8+\dx,\dy);
\node at (0.35,1.95) {$E$};
\node[gray!60!black] at (3+\dx,\dy+2.6) [right] {\small $E\otimes_Rk(\mathfrak p)$};
\node[red!70!black] at (6.6+\dx+0.3,\dy+2.6) [right] {\small $E_K$};
\end{scope}
\draw[->] (5.6,4.5) -- node[right]{$\tau$} (5.6,3.35);
\begin{scope}[font=\small]
\node[red!70!black,align=left,anchor=west] (b) at (10.4,11.0) {$\tilde v_K=\sum_if_ie_i\in\Gamma(\hat S_{s_K},\cE_K)$};
\draw[->,red!70!black] (b.west) .. controls (7.50,11.15) and (7.45,9.50) .. (7.45,7.85);
\node[orange!80!black,align=left,anchor=west] (c) at (10.4,9.3) {$\tilde v_{R_{\mathfrak p}}=\sum_if_i^{<g(p)}e_i\in\cE(U')\otimes_RR_{\mathfrak p}$};
\draw[->,orange!80!black] (c.west) .. controls (9.0,9.4) and (7.812,9.0) .. (7.812,8.071);
\draw[->,orange!80!black] (c.west) .. controls (9.65,8.70) and (8.80,7.35) .. (4.077,7.993);
\node[blue!60!black,align=left,anchor=west] (a) at (10.4,8.15) {$\tilde v_{k(\mathfrak p)}\in\Gamma\bigl(S,\cE/\cI_s^{g(p)}\cE\otimes_Rk(\mathfrak p)\bigr)$,\\ image of $\tilde v_{R_{\mathfrak p}}$ mod $(\cI_s^{g(p)},\mathfrak p)$};
\draw[->,blue!60!black] (a.west) .. controls (9.6,6.70) and (5.3,6.60) .. (3.83,7.83);
\node[orange!80!black,align=left,anchor=west] (c2) at (10.4,7.0) {$\overline{\tilde v}_{R_{\mathfrak p}}\in\Gamma(U'_{k(\mathfrak p)},\cE_{k(\mathfrak p)})$,\\ image of $\tilde v_{R_{\mathfrak p}}$ mod $\mathfrak p$};
\draw[->,orange!80!black] (c2.west) .. controls (9.0,5.95) and (4.9,5.95) .. (4.077,7.993);
\node[blue!60!black,align=left,anchor=west] (vbar) at (10.4,5.85) {$\bar v\in\Gamma(\Spec k(\mathfrak p),s_{k(\mathfrak p)}^*\cE_{k(\mathfrak p)})$,\\ image of $v$ mod $\mathfrak p$};
\draw[->,blue!60!black,shorten >=2.2pt] (vbar.west) .. controls (9.0,5.05) and (5.0,5.15) .. (3+0.45*\dx,5.6+0.45*\dy+1.25);
\node[blue!60!black,align=left,anchor=west] (d) at (10.4,1.2) {blue band: the support\\ $s(\Spec R)$ of $\cO_S/\cI_s^{g(p)}$};
\draw[->,blue!60!black] (d.west) .. controls (9.0,0.9) and (6.0,1.2) .. (5.2,2.2);
\node[blue!60!black,align=right,anchor=east] (g) at (-0.3,0.9) {blue disk around $s_{k(\mathfrak p)}$:\\ $\Spec\bigl(\cO_{S_{k(\mathfrak p)}}/\cI_{s,k(\mathfrak p)}^{g(p)}\bigr)$};
\draw[->,blue!60!black] (g.east) .. controls (2.2,0.55) and (3.5,0.6) .. (3.62,1.98);
\node[red!70!black,anchor=west] (hatS) at (10.4,2.45) {$\hat S_{s_K}=\Spf K[[x_1,\dots,x_d]]$};
\draw[->,red!70!black,shorten >=0.35cm] (hatS.west) -- (6.6+0.45*\dx,1.4+0.45*\dy);
\node[orange!80!black,align=right,anchor=east] (c3) at (-0.3,7.15) {$\overline{\tilde v}_{R_{\mathfrak p}}\in\Gamma(S_{k(\mathfrak p)},\cE_{k(\mathfrak p)})$,\\ image of $\tilde v_{R_{\mathfrak p}}$ mod $\mathfrak p$};
\draw[->,orange!80!black] (c3.east) .. controls (1.2,6.85) and (2.4,6.7) .. (3.096,6.658);
\node[black!35,anchor=east] (e) at (-0.3,6.4) {$F^1_{\conj}$};
\node[black!60,anchor=east] (f) at (-0.3,5.8) {$F^2_{\conj}$};
\draw[->,black!60] (f.east) .. controls (1.5,5.8) and (2.6,5.95) .. (3.55,6.32);
\draw[->,black!35] (e.east) .. controls (1.5,6.3) and (2.6,6.25) .. (3.45,6.4);
\begin{scope}
\clip[yshift=5.6cm] (3,0) -- (3+\dx,\dy) -- (3+\dx,\dy+0.77) .. controls (3+0.65*\dx,0.65*\dy+0.42) and (3+0.3*\dx,0.3*\dy+0.95) .. (3,0.56) -- cycle;
\draw[->,white] (e.east) .. controls (1.5,6.3) and (2.6,6.25) .. (3.45,6.4);
\end{scope}
\begin{scope}
\clip[yshift=5.6cm] (3,0) -- (3+\dx,\dy) -- (3+\dx,\dy+0.25) .. controls (3+0.65*\dx,0.65*\dy+0.38) and (3+0.3*\dx,0.3*\dy+0.1) .. (3,0.21) -- cycle;
\draw[->,white] (f.east) .. controls (1.5,5.8) and (2.6,5.95) .. (3.55,6.32);
\end{scope}
\node[orange!80!black,align=left,anchor=west] (h) at (10.4,3.4) {open neighborhood of $s_{k(\mathfrak p)}$ in $S_{k(\mathfrak p)}$ where\\ the image of $\overline{\tilde v}_{R_{\mathfrak p}}$ generates $\mathrm{gr}^0_{F_{\conj}}\cE_{k(\mathfrak p)}$};
\draw[->,orange!80!black] (h.west) .. controls (8.0,3.2) and (6.0,2.6) .. (4.024,2.552);
\end{scope}
\end{tikzpicture}%
}}
\captionsetup{labelformat=simple,labelsep=period}
\captionof{figure}{Proof of \Cref{thm:main}.}
\label{fig:step2}
\end{minipage}\end{center}
By \Cref{lem:omega-reduction}, for every vector field $\partial$ on $S_{k(\mathfrak p)}$,
\begin{equation}\label{eq:psi-vanish}
\psi_p(\partial)(\overline{\tilde v}_{R_{\mathfrak p}})\in\cI_{s,k(\mathfrak{p})}^{g(p)-p}\cE_{k(\mathfrak p)}.
\end{equation}
As $g(p)\ge 1$, $\overline{\tilde v}_{R_{\mathfrak p}}(s_{k(\mathfrak p)})=\bar v$ by \eqref{eq:vRp-bar-const}, with $s^*\cE\otimes_Rk(\mathfrak p)=R^2\pi_*\Omega^\bullet_{X_{s,k(\mathfrak{p})}/k(\mathfrak{p})}$.
So $\overline{\tilde v}_{R_{\mathfrak p}}$ and $\bar v$ have the same image in the fiber of $\mathrm{gr}^0_{F_{\conj}}\cE_{k(\mathfrak p)}$ at $s_{k(\mathfrak p)}$.
By \ref{item:step0-b}, $\bar v$ generates $F_{\mathrm{Hdg}}^2\bigl(R^2\pi_*\Omega^\bullet_{X_{s,k(\mathfrak{p})}/k(\mathfrak{p})}\bigr)$, which has rank $1$.
Hence $\bar v\neq0$.
As $\mathfrak p\in\mathcal P$, $X_{s,k(\mathfrak p)}$ is an ordinary $K3$ surface by \Cref{claim:good-primes}.
It satisfies (HdR) by \Cref{ex:K3-fzip}, so it is ordinary in the sense of \Cref{ex:ordinary-scheme} (by \Cref{rmk:ordinary-K3}).
As $k(\mathfrak{p})$ is perfect, \eqref{eq:opposite-n2-i2} shows that $\bar v\notin F^1_{\conj}$, so this image is nonzero.
A section of a line bundle that is nonzero in the fiber at a point generates the stalk there by Nakayama's lemma, and the locus where a section generates an invertible sheaf is open.
Therefore,
\begin{itemize}
\item the image of $\overline{\tilde v}_{R_{\mathfrak p}}$ in $\mathrm{gr}^0_{F_{\conj}}\cE_{k(\mathfrak p)}$ generates this line bundle on an open neighborhood of $s_{k(\mathfrak p)}$, and
\item the same section $\psi_p(\partial)(\overline{\tilde v}_{R_{\mathfrak p}})$ lies in $F^1_{\conj}\cE_{k(\mathfrak p)}$ by \eqref{eq:psi-conj}.
\end{itemize}
As $F^1_{\conj}\cE_{k(\mathfrak p)}$ is locally free by \Cref{lem:conj-bc}(1) for $\pi_{k(\mathfrak p)}$,
\begin{equation}\label{eq:F1-intersect}
F^1_{\conj}\cE_{k(\mathfrak p)}\cap\cI_{s,k(\mathfrak{p})}^{g(p)-p}\cE_{k(\mathfrak p)}=\cI_{s,k(\mathfrak{p})}^{g(p)-p}F^1_{\conj}\cE_{k(\mathfrak p)}
\end{equation}
on that neighborhood.
Hence $(\mathrm{gr}^0_{F_{\conj}}\psi_p)(\partial)$ sends the image of $\overline{\tilde v}_{R_{\mathfrak p}}$ into $\cI_{s,k(\mathfrak{p})}^{g(p)-p}\mathrm{gr}^1_{F_{\conj}}\cE_{k(\mathfrak p)}$ by \eqref{eq:psi-vanish} and \eqref{eq:F1-intersect}.
The map \eqref{eq:gr0psi-partial} is $\cO$-linear and the image of $\overline{\tilde v}_{R_{\mathfrak p}}$ generates the source there, so on that neighborhood, for every vector field $\partial$,
\[
(\mathrm{gr}^0_{F_{\conj}}\psi_p)(\partial)\bigl(\mathrm{gr}^0_{F_{\conj}}\cE_{k(\mathfrak p)}\bigr)\subset\cI_{s,k(\mathfrak{p})}^{g(p)-p}\mathrm{gr}^1_{F_{\conj}}\cE_{k(\mathfrak p)}.
\]
As $T_{S_{k(\mathfrak p)}/k(\mathfrak p)}$ is free by \ref{item:step0-d}, applying this to a basis $\partial_1,\dots,\partial_d$ shows that $\mathrm{gr}^0_{F_{\conj}}\psi_p$ vanishes to order $g(p)-p$ at $s$ modulo $\mathfrak{p}$ (\Cref{def:order-vanishing}).
\Cref{prop:katz-F2} for $(q,j)=(0,2)$ gives that $\theta$ vanishes to order $\lfloor g(p)/p\rfloor-1$ at $s$ modulo $\mathfrak{p}$.

\textbf{Step 3: $\theta=0$.}
Under the isomorphism $\hat{\cO}_{S,s}\simeq R[[x_1,\dots,x_d]]$, the ideal $\cI_s$ corresponds to $(x_1,\dots,x_d)$.
By the vanishing estimate of $\theta$ in Step 2 and \Cref{def:order-vanishing}(a), since the $h_k\otimes\omega_l$ in \eqref{eq:theta-kl} form a basis, we have
\[
\overline\theta_{kl}\in\cI_{s,k(\mathfrak p)}^{\lfloor g(p)/p\rfloor-1},
\]
where $\overline\theta_{kl}$ is the image of $\theta_{kl}$ in $\Gamma(S_{k(\mathfrak p)},\cO_{S_{k(\mathfrak p)}})$.
In other words, if we write $\theta_{kl}=\sum_{\alpha\in\Z_{\geq0}^d}a_{kl,\alpha}x^\alpha$ with $a_{kl,\alpha}\in R$, then
\[
a_{kl,\alpha}\in\mathfrak p\quad\text{for all multi-indices }\alpha\text{ with }|\alpha|<\lfloor g(p)/p\rfloor-1.
\]

Fix $k,l$ and a multi-index $\alpha$, and put $n:=|\alpha|$.
The set
\stepcounter{equation}
\begin{equation*}\tag*{(\theequation)$_n$}\label{eq:bad-primes}
\{p\text{ prime number}\mid g(p)\le p\text{ or }\lfloor g(p)/p\rfloor-1\le n\}
\end{equation*}
is finite.
Since only finitely many maximal ideals of $R$ lie over each prime and $\mathcal{P}$ is infinite, there are infinitely many $\mathfrak p\in\mathcal{P}$ that have residue characteristic outside the set \ref{eq:bad-primes}.
Let
\[
\mathcal{Q}_n:=\{\mathfrak p\in\mathcal{P}\mid\text{the residue characteristic of }\mathfrak p\text{ lies outside the set \ref{eq:bad-primes}}\}.
\]
Then by the analysis above,
\[
a_{kl,\alpha}\in\mathfrak p\quad\text{for all }\mathfrak p\in\mathcal{Q}_n.
\]
But a nonzero element of $R$ lies in only finitely many maximal ideals, because $\dim R=1$.
Thus
\[
a_{kl,\alpha}=0\quad\text{in }R.
\]
Every coefficient of each $\theta_{kl}$ in $R[[x_1,\dots,x_d]]$ is therefore zero.
Since $\Gamma(S,\cO_S)\to\cO_{S,s}$ is injective by integrality of $S$ at $s$, and $\cO_{S,s}\to\hat{\cO}_{S,s}$ is injective by the Krull intersection theorem \cite[\href{https://stacks.math.columbia.edu/tag/00IQ}{Tag 00IQ}]{SP}, we have
\[
\theta_{kl}=0\quad\text{in }\Gamma(S,\cO_S).
\]
So $\theta=0$ on $S$, contrary to hypothesis.
\end{proof}

\section{Some deformation theory}
\label{sec:deformation-theory}

Our first aim is to prove \Cref{thm:ks-isotrivial}.
Over an algebraically closed field, a version of the second half of its proof is recorded in \cite[Proposition 2.6.10]{Ser06}; see also \cite[Proposition 6.3]{Poc25}.

For a local ring $B$, $\mathfrak{m}_B$ denotes its maximal ideal; for a point $s$ of a scheme $S$, we write $\mathfrak{m}_s:=\mathfrak{m}_{\cO_{S,s}}$.
For a scheme $Y$ over a field $K$, $T_Y:=\mathcal{H}om_{\cO_Y}(\Omega^1_{Y/K},\cO_Y)$.

\begin{lemma}[{\cite[Exp.~221, \S4.c]{Gro61}}]
\label{lem:isom-representable}
Let $S$ be a noetherian scheme, and let $X\to S$ and $Y\to S$ be flat projective morphisms.
\begin{enumerate}[label=\textup{(\arabic*)},ref=\arabic*]
\item
\label{item:isom-rep-mor}
The functor
\begin{align*}
\mathfrak{Mor}_S(Y,X)\colon(\mathrm{LNSch}/S)^{\mathrm{op}}&\longrightarrow\mathrm{Sets}
\\
T&\longmapsto\{T\text{-morphisms }Y\times_S T\to X\times_S T\}
\end{align*}
where $\mathrm{LNSch}/S$ is the category of locally noetherian $S$-schemes, is represented by an open subscheme $\mathrm{Mor}_S(Y,X)$ of $\mathrm{Hilb}_{Y\times_S X/S}$.
\item
\label{item:isom-rep-isom}
The functor
\begin{align*}
\mathfrak{Isom}_S(Y,X)\colon(\mathrm{LNSch}/S)^{\mathrm{op}}&\longrightarrow\mathrm{Sets}
\\
T&\longmapsto\{T\text{-isomorphisms }Y\times_S T\to X\times_S T\},
\end{align*}
is represented by an open subscheme of $\mathrm{Mor}_S(Y,X)$.
In particular, it is represented by a scheme locally of finite type over $S$.
\end{enumerate}
\end{lemma}

\begin{proof}
\begin{enumerate}[label=\textup{(\arabic*)},ref=\arabic*]
\item
Put
\[
Z:=Y\times_S X.
\]
The projection $Z\to Y$ is projective, as a base change of the projective morphism $X\to S$.
Since $S$ is noetherian, it is qcqs, so the composite $Z\to Y\to S$ is projective \cite[\href{https://stacks.math.columbia.edu/tag/0C4P}{Tag 0C4P}]{SP}.
Choose a closed immersion $i\colon Z\to\mathbf{P}(\mathcal{E})$ over $S$, with $\mathcal{E}$ a quasi-coherent $\cO_S$-module of finite type \cite[\href{https://stacks.math.columbia.edu/tag/01W8}{Tag 01W8}]{SP}.
Then $L:=i^*\cO_{\mathbf{P}(\mathcal{E})}(1)$ is a relatively very ample line bundle on $Z$ over $S$ \cite[\href{https://stacks.math.columbia.edu/tag/01VM}{Tag 01VM}]{SP}.
Since
\[
\mathrm{Hilb}_{Z/S}=\mathrm{Quot}_{\cO_Z/Z/S}=\coprod_{\Phi\in\Q[t]}\mathrm{Quot}^{\Phi,L}_{\cO_Z/Z/S}
\]
\cite[\S\S5.1.3--5.1.4, pp.~109--110]{Nit05} and each $\mathrm{Quot}^{\Phi,L}_{\cO_Z/Z/S}$ is a projective $S$-scheme \cite[Theorem 5.14]{Nit05}, $\mathrm{Hilb}_{Z/S}$ is locally of finite type over $S$.
By \cite[Theorem 5.23]{Nit05},\footnote{\cite[Theorem 5.23]{Nit05} is stated on all $S$-schemes, but its proof uses \cite[Theorem 5.22(b)]{Nit05}, which only covers locally noetherian ones.
Also \cite[Exp.~221, \S4.c]{Gro61} is only stated for locally noetherian $S$-schemes \cite[Exp.~221, p.~221-11]{Gro61}.} the functor $\mathfrak{Mor}_S(Y,X)$ is represented by an open subscheme
\[
M:=\mathrm{Mor}_S(Y,X)
\]
of $\mathrm{Hilb}_{Z/S}$.
In particular, $M$ is locally of finite type over $S$.
\item
Let $u\colon Y\times_S M\to X\times_S M$ be the universal $M$-morphism, i.e., the $M$-morphism that corresponds to $\mathrm{id}_M$ when $M$ represents $\mathfrak{Mor}_S(Y,X)$, so that for every locally noetherian $S$-scheme $T$ the bijection is
\begin{align*}
\Hom_S(T,M)&\xlongrightarrow{\simeq}\mathfrak{Mor}_S(Y,X)(T)
\\
h&\mapsto u_T:=u\times_{M,h}T.
\end{align*}
Then $\mathfrak{Isom}_S(Y,X)$ is a subfunctor of $\mathfrak{Mor}_S(Y,X)$, and a morphism $h\colon T\to M$ lies in $\mathfrak{Isom}_S(Y,X)(T)$ if and only if $u_T$ is an isomorphism.

Let $U\subset M$ be an affine open.
The scheme $U$ is noetherian, since $M$ is locally of finite type over the noetherian $S$ by (\ref{item:isom-rep-mor}) \cite[\href{https://stacks.math.columbia.edu/tag/01T6}{Tag 01T6}]{SP}.
As base changes of the projective and flat morphisms $X\to S$ and $Y\to S$, the morphisms $X_U\to U$ and $Y_U\to U$ are projective and flat, hence proper.
In particular, $X_U\to U$ is separated, so $u_U\colon Y_U\to X_U$ is projective \cite[\href{https://stacks.math.columbia.edu/tag/0C4Q}{Tag 0C4Q}]{SP}.
By \cite[Theorem 5.22(b)]{Nit05}, there is an open $U_2\subset U$ such that, for every locally noetherian $U$-scheme $T$, $u_T$ is an isomorphism if and only if $T\to U$ factors through $U_2$.
For affine opens $U,U'\subset M$ and an affine open $W\subset U\cap U'$, testing the universal properties of $U_2$, $U'_2$ and $W_2$ on locally noetherian $W$-schemes gives $U_2\cap W=W_2=U'_2\cap W$.
So the $U_2$ glue to an open subscheme $M_2\subset M$, which represents $\mathfrak{Isom}_S(Y,X)$ because being an isomorphism is local on $T$.
As an open subscheme of $M$, which is locally of finite type over $S$ by (\ref{item:isom-rep-mor}), the scheme $M_2$ is locally of finite type over $S$.
\end{enumerate}
\end{proof}

\begin{lemma}
\label{lem:tangent-aut}
Let $K$ be a field, let $Y$ be a projective $K$-scheme, let $\mathrm{Aut}(Y)$ be the $K$-scheme representing $\mathfrak{Isom}_K(Y,Y)$ (\Cref{lem:isom-representable}(\ref{item:isom-rep-isom}) with $S=\Spec K$), and let $e:=\mathrm{id}_Y\in\mathrm{Aut}(Y)(K)$.
Put $S_{0,\varepsilon}:=\Spec(K[\varepsilon]/\varepsilon^2)$, let $\iota\colon\Spec K\to S_{0,\varepsilon}$ be the closed immersion, and let
\begin{equation}
\label{eq:rho}
\begin{aligned}
\rho\colon\mathrm{Aut}_{S_{0,\varepsilon}}\bigl(Y\times_{\Spec K}S_{0,\varepsilon}\bigr)&\longrightarrow\mathrm{Aut}_K(Y)
\\
g&\longmapsto g\times_{S_{0,\varepsilon}}\Spec K.
\end{aligned}
\end{equation}
Then
\[
T_e\mathrm{Aut}(Y)=\mathrm{Ker}(\rho)=H^0(Y,T_Y).
\]
\end{lemma}

\begin{proof}
By \Cref{lem:isom-representable}(\ref{item:isom-rep-isom}), naturally in locally noetherian $K$-schemes $T$,
\begin{equation}
\label{eq:aut-rep}
\mathrm{Hom}_K\bigl(T,\mathrm{Aut}(Y)\bigr)=\mathrm{Aut}_T\bigl(Y\times_{\Spec K}T\bigr),
\end{equation}
where $e$ corresponds to $\mathrm{id}_Y$ for $T=\Spec K$, and naturality along $\iota$ turns $t\mapsto t\circ\iota$ into $\rho$.
Let $\pi\colon Y\to\Spec K$ be the structure morphism, and let $\underline{K}=\pi^{-1}\cO_{\Spec K}$ be the constant sheaf $K$ on $Y$.
Let $\mathrm{NL}_{Y/K}$ be the naive cotangent complex of $Y$ over $K$, i.e.,
\[
\mathrm{NL}_{Y/K}=\Bigl(\underset{\text{degree }-1}{\mathcal{I}/\mathcal{I}^2}\longrightarrow\underset{\text{degree }0}{\Omega^1_{\underline{K}[\cO_Y]/\underline{K}}\otimes_{\underline{K}[\cO_Y]}\cO_Y}\Bigr)\quad\text{(\cite[\href{https://stacks.math.columbia.edu/tag/08TJ}{Tag 08TJ}]{SP})},
\]
where $\underline{K}[\cO_Y]$ is the sheafification of $U\mapsto\underline{K}(U)[\cO_Y(U)]$, and $\mathcal{I}:=\mathrm{Ker}\bigl(\underline{K}[\cO_Y]\to\cO_Y,\ [b]\mapsto b\bigr)$.
The ideal of the first order thickening $\Spec K\subset S_{0,\varepsilon}$ is $\varepsilon K[\varepsilon]/\varepsilon^2$, so
\[
\pi^*\bigl(\varepsilon K[\varepsilon]/\varepsilon^2\bigr)\simeq\pi^*\cO_{\Spec K}=\cO_Y.
\]
As $k(e)=K$,
\begin{align*}
T_e\mathrm{Aut}(Y)&=\mathrlap{\left\{t\in\mathrm{Hom}_K\bigl(S_{0,\varepsilon},\mathrm{Aut}(Y)\bigr)\,\middle|\,
\begin{tikzcd}[column sep=0.5em, row sep=0.8em, cells={nodes={font=\scriptsize}}, ampersand replacement=\&, baseline=(current bounding box.center)]
\& \Spec K \arrow[dl, "\iota"'] \arrow[dr, "e"] \&
\\
S_{0,\varepsilon} \arrow[rr, "t"'] \& \& \mathrm{Aut}(Y)
\end{tikzcd}
\text{ commutes}\right\}}
\\
&&&\text{(\cite[\href{https://stacks.math.columbia.edu/tag/0B2C}{Tag 0B2C}]{SP})}
\\
&=\bigl\{g\in\mathrm{Aut}_{S_{0,\varepsilon}}\bigl(Y\times_{\Spec K}S_{0,\varepsilon}\bigr)\bigm| g\times_{S_{0,\varepsilon}}\Spec K=\mathrm{id}_Y\bigr\}&&\text{(by \eqref{eq:aut-rep})}
\\
&=\mathrm{Ker}(\rho)
\\
&=\mathrm{Ext}^0_{\cO_Y}\bigl(\mathrm{NL}_{Y/K},\cO_Y\bigr)&&\text{(\cite[\href{https://stacks.math.columbia.edu/tag/0D14}{Tag 0D14}]{SP})}
\\
&=\mathrm{Hom}_{\cO_Y}(\Omega^1_{Y/K},\cO_Y)&&\text{(\cite[\href{https://stacks.math.columbia.edu/tag/0DY9}{Tag 0DY9}]{SP})}
\\
&=H^0(Y,T_Y).
\end{align*}
\end{proof}

\begin{lemma}
\label{lem:tangent-aut-etale}
Let $K$, $Y$, $\mathrm{Aut}(Y)$, and $e$ be as in \Cref{lem:tangent-aut}.
Let $S$ be a noetherian $K$-scheme, let $s\in S(K)$, let $\pi\colon X\to S$ be a flat projective morphism with $X_s=Y$, write $Y_S:=Y\times_K S$, let $I$ be the $S$-scheme representing $\mathfrak{Isom}_S(Y_S,X)$ (\Cref{lem:isom-representable}(\ref{item:isom-rep-isom})), and let $e\in I(K)$ also denote the point $\mathrm{id}_Y$ over $s$.
If
\begin{enumerate}[label=\textup{(\alph*)}]
\item $H^0(Y,T_Y)=0$, and
\item $I\to S$ admits a formal section through $e$ (i.e., a local $K$-algebra map $\sigma\colon\hat{\cO}_{I,e}\to\hat{\cO}_{S,s}$ with $\sigma\circ(\hat{\cO}_{S,s}\to\hat{\cO}_{I,e})=\mathrm{id}$),
\end{enumerate}
then $I\to S$ is \'{e}tale at $e$.
\end{lemma}

\begin{proof}
As $Y_S\to S$ is a base change of $\pi_s\colon Y\to\Spec K$, it is flat and projective, so by \Cref{lem:isom-representable}(\ref{item:isom-rep-isom}) $I$ is locally of finite type over $S$.
As $X_s=Y$, the fiber of $I$ over $s$ is $\mathrm{Aut}(Y)$.
We claim that
\begin{equation}
\label{eq:unram-e}
\Omega^1_{I/S}\otimes k(e)=0.
\end{equation}
In fact,
\begin{align*}
\mathrm{Hom}_{\cO_{I,e}}\bigl((\Omega^1_{I/S})_e,k(e)\bigr)&=T_{I/S,e}&&\text{(\cite[\href{https://stacks.math.columbia.edu/tag/0B2D}{Tag 0B2D}]{SP})}
\\
&=T_e\mathrm{Aut}(Y)&&\text{(\cite[\href{https://stacks.math.columbia.edu/tag/0B2C}{Tag 0B2C}]{SP})}
\\
&=H^0(Y,T_Y)&&\text{(by \Cref{lem:tangent-aut})}
\\
&=0,
\end{align*}
which proves the claim.

Both $\cO_{S,s}$ and $\cO_{I,e}$ are noetherian, as $I$ is locally of finite type over the noetherian scheme $S$.
As $k(e)=k(s)=K$, we have
\[
\mathfrak{m}_e/(\mathfrak{m}_e^2+\mathfrak{m}_s\cO_{I,e})=\Omega^1_{I/S}\otimes k(e)\overset{\eqref{eq:unram-e}}{=}0
\]
(The first equality here is \cite[\href{https://stacks.math.columbia.edu/tag/0B2E}{Tag 0B2E}]{SP}, i.e., $\mathfrak{m}_e=\mathfrak{m}_s\cO_{I,e}+\mathfrak{m}_e^2$.)
As $\mathfrak{m}_e$ is a finite $\cO_{I,e}$-module, Nakayama's lemma \cite[\href{https://stacks.math.columbia.edu/tag/00DV}{Tag 00DV}(4)]{SP} gives
\[
\mathfrak{m}_s\cO_{I,e}=\mathfrak{m}_e,
\]
so $\dim_K\cO_{I,e}/\mathfrak{m}_s\cO_{I,e}=1$.
Hence $\hat{\cO}_{I,e}$ is the $\mathfrak{m}_s$-adic completion of $\cO_{I,e}$ and a finite $\hat{\cO}_{S,s}$-module \cite[\href{https://stacks.math.columbia.edu/tag/0394}{Tag 0394}(1)(2)]{SP}, and
\[
\hat{\cO}_{I,e}/\mathfrak{m}_s\hat{\cO}_{I,e}=\cO_{I,e}/\mathfrak{m}_s\cO_{I,e}=K
\]
\cite[\href{https://stacks.math.columbia.edu/tag/031C}{Tag 031C}]{SP}.
In particular, $\hat{\cO}_{S,s}\to\hat{\cO}_{I,e}$ is surjective.
With the section $\sigma$, it is an isomorphism, in particular \'{e}tale.
As $\cO_{I,e}$ is essentially of finite type over $\cO_{S,s}$, the map $\cO_{S,s}\to\cO_{I,e}$ is therefore an \'{e}tale homomorphism of local rings \cite[\href{https://stacks.math.columbia.edu/tag/039M}{Tag 039M}(1)(3)]{SP}.
Since $S$ is locally noetherian and $I\to S$ is locally of finite type, this means that $I\to S$ is \'{e}tale at $e$.
\end{proof}

\begin{para}
\label{para:prorep}
Let
\[
\mathcal{A}\mathit{rt}_K:=\bigl\{A'\text{ Artinian local $K$-algebra}\bigm| K\xrightarrow{\sim}A'/\mathfrak{m}_{A'}\bigr\},
\]
with local $K$-algebra homomorphisms as morphisms.
Let $S$ be a $K$-scheme, let $s\in S(K)$, and for $m\geq1$ put
\[
B_m:=\cO_{S,s}/\mathfrak{m}_s^m,
\]
so that $\hat{\cO}_{S,s}=\varprojlim_mB_m$.

Let $\pi\colon X\to S$ be a flat morphism, put $Y:=X_s$, let $\pi_s\colon Y\to\Spec K$ be the structure morphism.
We assume that
\begin{enumerate}[label=\textup{(\alph*)}]
\item $Y$ is proper over $K$;
\item\label{item:prorep-b} $H^0(Y,T_Y)=0$.
\end{enumerate}
Let
\begin{align*}
\overline{\mathrm{Def}}_Y\colon\mathcal{A}\mathit{rt}_K&\longrightarrow\mathrm{Sets}
\\
A'&\longmapsto\bigl\{(W,a)\bigm| W\text{ a flat }A'\text{-scheme},\ a\colon W\times_{\Spec A'}\Spec K\xrightarrow{\sim}Y\bigr\}/{\cong},
\end{align*}
where $(W,a)\cong(W',a')$ if there is an $A'$-isomorphism $W\xrightarrow{\sim}W'$ compatible with $a$ and $a'$.
On a morphism $\mu\colon A'\to A''$ of $\mathcal{A}\mathit{rt}_K$, it sends the class of $(W,a)$ to the class of $W\times_{\Spec A',\Spec\mu}\Spec A''$ with the isomorphism induced by $a$.
Since
\begin{equation}
\label{eq:der-vanish}
\begin{aligned}
\mathrm{Der}_K(\cO_Y,\cO_Y)&=\mathrm{Hom}_{\cO_Y}(\Omega^1_{Y/K},\cO_Y)
\\
&=H^0(Y,T_Y)
\\
&=0
\end{aligned}
\end{equation}
and $Y$ is proper over $K$, the functor $\overline{\mathrm{Def}}_Y$ is prorepresentable \cite[\href{https://stacks.math.columbia.edu/tag/0ET5}{Tag 0ET5}]{SP}.
By \cite[\href{https://stacks.math.columbia.edu/tag/06GX}{Tag 06GX}]{SP}, this means that there are a complete local noetherian $K$-algebra $A$ with residue field $K$ and bijections
\begin{equation}
\label{eq:prorep}
\mathrm{Hom}_{K\text{-loc}}(A,A')\xlongrightarrow{\simeq}\overline{\mathrm{Def}}_Y(A')
\end{equation}
natural in $A'\in\mathcal{A}\mathit{rt}_K$, where $\mathrm{Hom}_{K\text{-loc}}$ denotes the set of local $K$-algebra homomorphisms.
Let $q_A\colon A\to K$ be the residue map.

For $m\geq1$ put
\[
S_m:=\Spec B_m,\qquad X_m:=X\times_S S_m.
\]
We have $X_1=Y$.
As a base change of the flat morphism $\pi$, $X_m$ is flat over $B_m$, so $[X_m,\mathrm{id}_Y]\in\overline{\mathrm{Def}}_Y(B_m)$.
Let
\[
\chi_m\colon A\to B_m
\]
be the preimage of $[X_m,\mathrm{id}_Y]$ under \eqref{eq:prorep}.
Since $X_{m+1}\times_{S_{m+1}}S_m=X_m$, naturality of \eqref{eq:prorep} along $B_{m+1}\to B_m$ shows that the $\chi_m$ are compatible.
So they define a local $K$-algebra map
\begin{equation}
\label{eq:chi}
\chi:=\varprojlim_m\chi_m\colon A\longrightarrow\hat{\cO}_{S,s}.
\end{equation}

Let $\mathcal{D}\mathit{ef}$ be the category defined in \cite[\href{https://stacks.math.columbia.edu/tag/0DY7}{Tag 0DY7}]{SP}:
\begin{itemize}
\item an object is a pair $(A',W)$ of an object $A'$ of $\mathcal{A}\mathit{rt}_K$ and a scheme $W$ flat over $A'$;
\item a morphism $(\mu,h)\colon(B',W')\to(A',W)$ is a morphism $\mu\colon B'\to A'$ in $\mathcal{A}\mathit{rt}_K$ together with a cartesian square
\[
\begin{tikzcd}[column sep=3.6em, row sep=2.6em, ampersand replacement=\&]
W \arrow[r, "h"] \arrow[d] \arrow[dr, phantom, "\lrcorner", very near start] \& W' \arrow[d]
\\
\Spec A' \arrow[r, "\Spec\mu"'] \& \Spec B'.
\end{tikzcd}
\]
\end{itemize}
For $A'\in\mathcal{A}\mathit{rt}_K$, $\mathcal{D}\mathit{ef}(A')$ denotes the fiber category:
\begin{itemize}
\item its objects are the pairs $(A',W)$;
\item its morphisms are the morphisms $(\mathrm{id}_{A'},h)$.
\end{itemize}
\end{para}

\begin{lemma}
\label{lem:def-classifying}
In the setting of \ref{para:prorep}, assume that
\begin{enumerate}[label=\textup{(\alph*)}]
\item $Y$ is projective over $K$;
\item $\chi$ as defined in \eqref{eq:chi} agrees with the composition
\[
A\xrightarrow{q_A}K\subset\hat{\cO}_{S,s}.
\]
\end{enumerate}
Then for every $m\geq1$ there is a unique isomorphism
\[
\alpha_m\colon Y\times_K S_m\xlongrightarrow{\sim}X_m
\]
over $S_m$ restricting to $\mathrm{id}_Y$ modulo $\mathfrak{m}_s$, and $\alpha_{m+1}$ restricts to $\alpha_m$.
\end{lemma}

\begin{proof}
\textbf{Existence.}
We compare two descriptions of $\chi_m\colon A\to B_m$:
\begin{itemize}
\item by definition, $\chi_m$ is the preimage of $[X_m,\mathrm{id}_Y]$ under \eqref{eq:prorep};
\item by (b), $\chi_m=\bigl(A\xrightarrow{q_A}K\to B_m\bigr)$, which by naturality of \eqref{eq:prorep} along $K\to B_m$ corresponds to $[Y\times_K S_m,\mathrm{id}_Y]$, as $q_A$ corresponds to $[Y,\mathrm{id}_Y]$.
\end{itemize}
Since \eqref{eq:prorep} is bijective, $[X_m,\mathrm{id}_Y]=[Y\times_K S_m,\mathrm{id}_Y]$ in $\overline{\mathrm{Def}}_Y(B_m)$, i.e., for every $m\geq1$ there is an isomorphism $\alpha_m$ as in the statement.

\textbf{Uniqueness.}
If $\alpha'_m$ is another such isomorphism, then $\alpha_m^{-1}\circ\alpha'_m$ is an element of
\[
\mathrm{Ker}\bigl(\mathrm{Aut}_{S_m}(Y\times_K S_m)\to\mathrm{Aut}_K(Y)\bigr),
\]
where the map is restriction along $Y\hookrightarrow Y\times_K S_m$.
Moreover, $\mathcal{D}\mathit{ef}$ satisfies
\begin{itemize}
\item the Rim--Schlessinger condition (RS) \cite[\href{https://stacks.math.columbia.edu/tag/0DY8}{Tag 0DY8}]{SP}, that is \cite[\href{https://stacks.math.columbia.edu/tag/06J2}{Tag 06J2}]{SP}: \textit{for all morphisms $(A_1,W_1)\to(A_0,W_0)$ and $(A_2,W_2)\to(A_0,W_0)$ in $\mathcal{D}\mathit{ef}$ lying over morphisms $A_1\to A_0$ and $A_2\to A_0$ in $\mathcal{A}\mathit{rt}_K$ with $A_2\to A_0$ surjective, there is a fiber product $(A_1,W_1)\times_{(A_0,W_0)}(A_2,W_2)$ in $\mathcal{D}\mathit{ef}$ lying over $A_1\times_{A_0}A_2$};
\item $\mathrm{Ker}(\rho)=H^0(Y,T_Y)=0$ by \Cref{lem:tangent-aut} and assumption \hyperref[item:prorep-b]{(b)} of \Cref{para:prorep}.
\end{itemize}
Therefore, by \cite[\href{https://stacks.math.columbia.edu/tag/06JZ}{Tag 06JZ}]{SP}\footnote{$\epsilon$ here is squarezero as you imagine: it is recorded in \cite[\href{https://stacks.math.columbia.edu/tag/06HV}{Tag 06HV}]{SP}.}, $\mathrm{Ker}\bigl(\mathrm{Aut}_{S_m}(Y\times_K S_m)\to\mathrm{Aut}_K(Y)\bigr)=0$.
So $\alpha_m$ is unique.

\textbf{Compatibility.}
By uniqueness, $\alpha_{m+1}$ restricts to $\alpha_m$ via base change along $S_m\subset S_{m+1}$.
\end{proof}

\begin{para}
\label{para:hat-derivation}
We fix an affine open neighborhood $U$ of $s$ and $\partial\in\Gamma(U,T_{S/K})$.
Via $\Gamma(U,T_{S/K})\simeq\mathrm{Der}_K(\cO_U,\cO_U)$, the section $\partial$ is a $K$-derivation $\partial\colon\cO_U\to\cO_U$, and so is its stalk $\partial_s\colon\cO_{S,s}\to\cO_{S,s}$.
By the Leibniz rule, $\partial(\mathfrak{m}_s^{m+1})\subset\mathfrak{m}_s^m$, so $\partial$ induces $K$-linear maps $\partial^{(m)}\colon B_{m+1}\to B_m$, compatible in $m$, and a $K$-derivation
\[
\hat{\partial}:=\varprojlim_m\partial^{(m)}\colon\hat{\cO}_{S,s}\longrightarrow\hat{\cO}_{S,s}.
\]
\end{para}

\begin{lemma}
\label{lem:eta-sigma}
Assume the setting of \ref{para:prorep} and \ref{para:hat-derivation}.
Put $X_U:=\pi^{-1}(U)$, and let $\tilde{D}\in\Gamma(X_U,T_{X/K})$ be a section mapping to $\pi^*\partial$ under $T_{X/K}\to \pi^*T_{S/K}$.
For $m\geq1$, write $b\mapsto\bar{b}$ for the projection $B_{m+1}\to B_m$, and put
\begin{align*}
\eta_m\colon B_{m+1}&\longrightarrow B_m[\varepsilon]/\varepsilon^2, &\eta'_m\colon B_{m+1}&\longrightarrow B_m[\varepsilon]/\varepsilon^2,
\\
b&\longmapsto\bar{b}+\varepsilon\,\partial^{(m)}(b), &b&\longmapsto\bar{b}.
\end{align*}
Define
\begin{align*}
\gamma_m\colon\cO_{X_{m+1}}&\longrightarrow\cO_{X_m}[\varepsilon]/\varepsilon^2
\\
y&\longmapsto\bar{y}+\varepsilon\,\overline{\tilde{D}(y)}.
\end{align*}
Then the following hold.
\begin{enumerate}
\item
$\eta_m$ and $\eta'_m$ are morphisms of $\mathcal{A}\mathit{rt}_K$.
\item
$\tilde{D}$ is a lift of $\partial$, i.e., the following diagram is commutative:
\begin{equation}
\label{eq:ks-lift}
\begin{tikzcd}
\pi^{-1}\cO_U\arrow[r,"\pi^{-1}\partial"]\arrow[d] & \pi^{-1}\cO_U\arrow[d]
\\
\cO_{X_U}\arrow[r,"\tilde{D}"] & \cO_{X_U}\rlap{.}
\end{tikzcd}
\end{equation}
\item\label{item:etasigma-3}
$\gamma_m$ is a well-defined map of sheaves of rings, hence defines a morphism of schemes
\[
\Spec\gamma_m\colon X_m\times_{S_m}\Spec\bigl(B_m[\varepsilon]/\varepsilon^2\bigr)\longrightarrow X_{m+1}
\]
over $\Spec\eta_m\colon\Spec(B_m[\varepsilon]/\varepsilon^2)\to S_{m+1}$.
\end{enumerate}
\end{lemma}

\begin{proof}
\begin{enumerate}
\item
The map $\eta_m$ is a ring map: it is additive, $\eta_m(1)=1$ as $\partial^{(m)}(1)=0$, and for $a,b\in B_{m+1}$
\begin{align*}
\eta_m(a)\,\eta_m(b)&=\bigl(\bar a+\varepsilon\,\partial^{(m)}(a)\bigr)\bigl(\bar b+\varepsilon\,\partial^{(m)}(b)\bigr)
\\
&=\overline{ab}+\varepsilon\bigl(\bar a\,\partial^{(m)}(b)+\bar b\,\partial^{(m)}(a)\bigr)&&(\text{because }\varepsilon^2=0)
\\
&=\overline{ab}+\varepsilon\,\partial^{(m)}(ab)&&(\text{Leibniz rule})
\\
&=\eta_m(ab).
\end{align*}
The ring $B_m[\varepsilon]/\varepsilon^2$ is an Artinian local $K$-algebra with maximal ideal $\mathfrak{m}_{B_m}+\varepsilon B_m$ and residue field $K$; both $\eta_m$ and $\eta'_m$ are $K$-linear, for $\eta_m$ as $\partial^{(m)}(K)=0$ since $\partial$ is a $K$-derivation; and both map $\mathfrak{m}_{B_{m+1}}$ into $\mathfrak{m}_{B_m}+\varepsilon B_m$.
\item
By the universal property of K\"ahler differentials (and as $\Omega^1_{S/K}$ is locally free),
\begin{align*}
T_{X/K}|_{X_U}&=\mathcal{H}om_{\cO_{X_U}}(\Omega^1_{X_U/K},\cO_{X_U})=\mathrm{Der}_K(\cO_{X_U},\cO_{X_U}),
\\
\pi^*T_{S/K}|_{X_U}&=\mathcal{H}om_{\cO_{X_U}}(\pi^*\Omega^1_{U/K},\cO_{X_U})
\\
&=\mathcal{H}om_{\pi^{-1}\cO_U}(\pi^{-1}\Omega^1_{U/K},\cO_{X_U})
\\
&=\mathrm{Der}_K(\pi^{-1}\cO_U,\cO_{X_U}).
\end{align*}
Under these identifications, the map $T_{X/K}\to \pi^*T_{S/K}$ is precomposition with $\pi^*\Omega^1_{U/K}\to\Omega^1_{X_U/K}$, and as $\tilde{D}$ maps to $\pi^*\partial$, the commutative diagram
\[
\begin{tikzcd}
\pi^*\Omega^1_{U/K}\arrow[r,"\pi^*\partial"]\arrow[d] & \cO_{X_U}\arrow[d,equal]
\\
\Omega^1_{X_U/K}\arrow[r,"\tilde{D}"] & \cO_{X_U}\rlap{,}
\end{tikzcd}
\]
translates to \eqref{eq:ks-lift}.
\item
Let $\mathcal{J}\subset\cO_U$ be the ideal sheaf of $s$.
By the Leibniz rule and \eqref{eq:ks-lift},
\[
\tilde{D}(\mathcal{J}^{m+1}\cO_{X_U})\subset\mathcal{J}^m\cO_{X_U}.
\]
So $\gamma_m$ is well defined, and it is a map of sheaves of rings (as $\tilde{D}$ is a derivation and $\varepsilon^2=0$, by the same reasoning as for $\eta_m$ in (1)).
The map $\gamma_m$ restricts to $\eta_m$ on $B_{m+1}$: for $b\in B_{m+1}$,
\[
\gamma_m(b)=\bar{b}+\varepsilon\,\overline{\tilde{D}(b)}\overset{\eqref{eq:ks-lift}}{=}\bar{b}+\varepsilon\,\partial^{(m)}(b)=\eta_m(b),
\]
and modulo $\varepsilon$ the map $\gamma_m$ is the projection $\cO_{X_{m+1}}\to\cO_{X_m}$.
Hence $\gamma_m$ defines the morphism of the statement.
\end{enumerate}
\end{proof}

\begin{para}
\label{para:ks}
We leave the setting of \ref{para:prorep} for a moment.
Recall from \cite[(1.1.1)--(1.1.3)]{Kat72} that, for a ring $A$, a smooth $A$-scheme $S$, and a smooth morphism $\pi\colon X\to S$ (so $X$ is smooth over $A$), the relative tangent sequence
\begin{equation}\label{eq:ks-ses}
0\longrightarrow T_{X/S}\longrightarrow T_{X/A}\longrightarrow \pi^*T_{S/A}\longrightarrow 0
\end{equation}
is exact: it is the dual of the short exact sequence
\[
0\longrightarrow \pi^*\Omega^1_{S/A}\longrightarrow\Omega^1_{X/A}\longrightarrow\Omega^1_{X/S}\longrightarrow0
\]
of locally free sheaves, so the sequence is locally split and remains exact after dualizing, and that the \emph{Kodaira--Spencer class} of $\pi$ is the composite of the unit and the boundary map of $\pi_*$ applied to this sequence,
\[
\mathrm{KS}_{\pi}\colon T_{S/A}\longrightarrow \pi_*\pi^*T_{S/A}\xrightarrow{\ \delta\ }R^1\pi_*T_{X/S}.
\]
\end{para}

\begin{lemma}
\label{lem:ks-formal-derivative}
In the setting of \ref{para:prorep} and \ref{para:hat-derivation}, assume moreover that $S$ is smooth over $K$, that $\pi$ is smooth and qcqs, and that
\[
\mathrm{KS}_{\pi}(\partial)=0\quad\text{on }U.
\]
Then
\begin{equation}
\label{eq:dpsi-zero}
\hat{\partial}\circ\chi=0\colon A\longrightarrow\hat{\cO}_{S,s}.
\end{equation}
\end{lemma}

\begin{proof}
Applying $\pi_*$ to the short exact sequence \eqref{eq:ks-ses} over the base ring $K$, we get an exact sequence of $\cO_S$-modules
\[
0\longrightarrow \pi_*T_{X/S}\longrightarrow \pi_*T_{X/K}\longrightarrow \pi_*\pi^*T_{S/K}\xrightarrow{\ \delta\ }R^1\pi_*T_{X/S}.
\]
So the image of $\pi_*T_{X/K}\to \pi_*\pi^*T_{S/K}$ is $\mathrm{Ker}(\delta)$, and the kernel is $\pi_*T_{X/S}$ by left exactness of $\pi_*$:
\[
0\longrightarrow \pi_*T_{X/S}\longrightarrow \pi_*T_{X/K}\longrightarrow\mathrm{Ker}(\delta)\longrightarrow0.
\]
$\mathrm{KS}_{\pi}$ is the unit $\partial\mapsto \pi^*\partial$ followed by $\delta$ by definition, so
\[
\delta(\pi^*\partial)=\mathrm{KS}_{\pi}(\partial)=0\quad\text{on }U,
\]
i.e., $\pi^*\partial\in\mathrm{Ker}(\delta)(U)$.
As $U$ is affine and $\pi_*T_{X/S}$ is quasi-coherent (because $\pi$ is qcqs),
by Serre vanishing
\[
H^1(U,\pi_*T_{X/S})=0.
\]
So $H^0(U,\pi_*T_{X/K})\to H^0(U,\mathrm{Ker}(\delta))$ is surjective, and there is
\[
\tilde{D}\in H^0(U,\pi_*T_{X/K})=H^0(X_U,T_{X/K})
\]
mapping to $\pi^*\partial$.

Fix $m\geq1$.
Put
\[
S_\varepsilon:=\Spec\bigl(B_m[\varepsilon]/\varepsilon^2\bigr),\qquad X_{m,\varepsilon}:=X_m\mathop{\times}\limits_{S_m}S_\varepsilon,\qquad X_{m+1,\varepsilon}:=X_{m+1}\mathop{\times}\limits_{S_{m+1},\,\Spec\eta_m}S_\varepsilon.
\]
Let $\Spec\gamma_m$ be the morphism of \Cref{lem:eta-sigma}\hyperref[item:etasigma-3]{(3)} induced by $\tilde{D}$.
As it lies over $\Spec\eta_m$, together with the projection to $S_\varepsilon$ it induces by the universal property of the fiber product a unique morphism $h_m$ making the following diagram commute:
\[
\begin{tikzcd}[column sep=3.6em, row sep=2.6em, ampersand replacement=\&]
X_{m,\varepsilon}
\arrow[ddr, bend right=12]
\arrow[drr, bend left=12, "\Spec\gamma_m"]
\arrow[dr, dashed, "h_m"']
\& \&
\\
\&
X_{m+1,\varepsilon}
\arrow[r]
\arrow[d]
\arrow[dr, phantom, "\lrcorner", very near start]
\&
X_{m+1}
\arrow[d]
\\
\&
S_\varepsilon
\arrow[r, "\Spec\eta_m"']
\&
S_{m+1}\rlap{.}
\end{tikzcd}
\]
Then $h_m$ is a morphism over $S_\varepsilon$, and as $\eta_m$ reduces modulo $\varepsilon$ to the projection $B_{m+1}\to B_m$, its base change along $S_m\subset S_\varepsilon$ is $\mathrm{id}_{X_m}$.
So we have a commutative diagram of thickenings
\[
\begin{tikzcd}[column sep=1.2em, row sep=2.6em, ampersand replacement=\&]
(X_m\subset X_{m,\varepsilon})
\arrow[rr, "{(\mathrm{id}_{X_m},\,h_m)}"]
\arrow[dr]
\& \&
(X_m\subset X_{m+1,\varepsilon})
\arrow[dl]
\\
\&
(S_m\subset S_\varepsilon)
\&
\end{tikzcd}
\]
as in \cite[\href{https://stacks.math.columbia.edu/tag/06AG}{Tag 06AG}]{SP}, where:
\begin{itemize}
\item $S_m\subset S_\varepsilon$ is a first-order thickening;
\item $X_{m,\varepsilon}$ is flat over $S_\varepsilon$, as a base change of the flat morphism $X_m\to S_m$;
\item $X_m$ is the base change of both $X_{m,\varepsilon}$ and $X_{m+1,\varepsilon}$ along $S_m\subset S_\varepsilon$.
\end{itemize}
As $\mathrm{id}_{X_m}$ is an isomorphism, so is $h_m$ by part~(2) of \textit{loc.\ cit.}.
We claim that the following equality holds in $\overline{\mathrm{Def}}_Y(B_m[\varepsilon]/\varepsilon^2)$:
\begin{equation}
\label{eq:claim-eta}
\overline{\mathrm{Def}}_Y(\eta_m)\bigl([X_{m+1},\mathrm{id}_Y]\bigr)=\overline{\mathrm{Def}}_Y(\eta'_m)\bigl([X_{m+1},\mathrm{id}_Y]\bigr).
\end{equation}
Assuming \eqref{eq:claim-eta} for the moment, by naturality of \eqref{eq:prorep}, $\eta_m\circ\chi_{m+1}=\eta'_m\circ\chi_{m+1}$, i.e.,
\[
\partial^{(m)}\bigl(\chi_{m+1}(a)\bigr)=0\quad\text{in }B_m\qquad(a\in A).
\]
Since $\partial^{(m)}\circ\chi_{m+1}$ is the reduction of $\hat{\partial}\circ\chi$ modulo $\mathfrak{m}_s^m$ and this holds for every $m\geq1$,
\[
\hat{\partial}\circ\chi=0\colon A\longrightarrow\hat{\cO}_{S,s},
\]
which is \eqref{eq:dpsi-zero}.

It remains to prove \eqref{eq:claim-eta}.
By definition,
\begin{align*}
\overline{\mathrm{Def}}_Y(\eta_m)\bigl([X_{m+1},\mathrm{id}_Y]\bigr)&=[X_{m+1,\varepsilon},\mathrm{id}_Y],
\\
\overline{\mathrm{Def}}_Y(\eta'_m)\bigl([X_{m+1},\mathrm{id}_Y]\bigr)&=\bigl[X_{m+1}\mathop{\times}\limits_{S_{m+1},\,\Spec\eta'_m}S_\varepsilon,\mathrm{id}_Y\bigr],
\end{align*}
i.e., the leftmost and the rightmost objects in the following diagram:
\[
\begin{tikzcd}[column sep=3.2em, row sep=2.6em, ampersand replacement=\&]
(Y\subset X_{m+1,\varepsilon})
\arrow[dr, "{(\pi_s,\,\mathrm{pr}_2)}"', pos=0.3]
\&
(Y\subset X_{m,\varepsilon})
\arrow[l, "{(\mathrm{id}_Y,\,h_m)}"', "\sim"]
\arrow[r, "{(\mathrm{id}_Y,\,\mathrm{can})}", "\sim"']
\arrow[d, "{(\pi_s,\,\mathrm{pr}_2)}"']
\&
(Y\subset X_{m+1}\mathop{\times}\limits_{S_{m+1},\,\Spec\eta'_m}S_\varepsilon)
\arrow[dl, "{(\pi_s,\,\mathrm{pr}_2)}", pos=0.3]
\\
\&
(\Spec K\subset\Spec(B_m[\varepsilon]/\varepsilon^2))
\&
\end{tikzcd}
\]
As the two horizontal arrows are isomorphisms over $\Spec(B_m[\varepsilon]/\varepsilon^2)$ restricting to $\mathrm{id}_Y$, these two classes agree, which proves the claim.
\end{proof}

\begin{lemma}
\label{lem:const-formal}
In the setting of \ref{para:prorep}, let $C$ be a local $K$-algebra with residue field $K$ and residue map $q_C\colon C\to K$, and let $\chi'\colon C\to\hat{\cO}_{S,s}$ be a local $K$-algebra map.
Assume:
\begin{enumerate}[label=\textup{(\alph*)}]
\item $\operatorname{char}K=0$;
\item $S$ is smooth over $K$;
\item $\hat{\partial}\circ\chi'=0$ for every choice of $U$ and $\partial$ in \ref{para:hat-derivation}.
\end{enumerate}
Then
\[
\chi'=\bigl(C\xrightarrow{q_C}K\subset\hat{\cO}_{S,s}\bigr).
\]
\end{lemma}

\begin{proof}
As $S$ is smooth over $K$, $\Omega^1_{S/K}$ is finite locally free; let $r$ be its rank at $s$.
As $k(s)=K$,
\[
\mathfrak{m}_s/\mathfrak{m}_s^2=\Omega^1_{S/K,s}\otimes_{\cO_{S,s}}K,
\]
where the class of $x{}\in\mathfrak{m}_s$ corresponds to $dx\otimes1$.
Choose $x_1,\dots,x_r\in\mathfrak{m}_s$ whose classes form a basis of $\mathfrak{m}_s/\mathfrak{m}_s^2$.
Then $dx_1,\dots,dx_r$ generate the free $\cO_{S,s}$-module $\Omega^1_{S/K,s}$ of rank $r$, hence form a basis of it.
Let $\partial_1,\dots,\partial_r$ be the dual basis of
\[
T_{S/K,s}=\mathrm{Hom}_{\cO_{S,s}}(\Omega^1_{S/K,s},\cO_{S,s})\quad(\Omega^1_{S/K}\text{ is finite locally free}),
\]
and choose an affine open neighborhood $U$ of $s$ such that $x_1,\dots,x_r\in\cO_S(U)$ and $\partial_1,\dots,\partial_r$ are sections of $T_{S/K}$ over $U$.
By \ref{para:hat-derivation}, each $\partial_i$ is a $K$-derivation of $\cO_U$ and induces a $K$-derivation $\hat{\partial}_i$ of $\hat{\cO}_{S,s}$; moreover
\[
\partial_i(x_j)=\delta_{ij}\quad\text{in }\cO_{S,s}\qquad(1\leq i,j\leq r).
\]
By assumption,
\[
\hat{\partial}_i\circ\chi'=0\qquad(1\leq i\leq r).
\]
The local ring $\cO_{S,s}$ is regular, hence so is $\hat{\cO}_{S,s}$ \cite[\href{https://stacks.math.columbia.edu/tag/07NY}{Tag 07NY}]{SP}.
By the Cohen structure theorem \cite[\href{https://stacks.math.columbia.edu/tag/0C0S}{Tag 0C0S}(2)]{SP}, there is a $K$-algebra isomorphism
\[
\hat{\cO}_{S,s}\simeq K[[y_1,\dots,y_r]].
\]
Put $\mathfrak{n}:=(y_1,\dots,y_r)$.
Then $\mathfrak{n}^m=\mathfrak{m}_s^m\hat{\cO}_{S,s}=\mathrm{Ker}(\hat{\cO}_{S,s}\to B_m)$ \cite[\href{https://stacks.math.columbia.edu/tag/05GG}{Tag 05GG}]{SP}, so $\hat{\partial}_i(\mathfrak{n}^{m+1})\subset\mathfrak{n}^m$ for all $m\geq1$.

We claim that a $K$-derivation $\Theta$ of $K[[y_1,\dots,y_r]]$ with $\Theta(\mathfrak{n}^{m+1})\subset\mathfrak{n}^m$ for all $m\geq1$ satisfies
\begin{equation}\label{eq:theta-expansion}
\Theta=\sum_{k=1}^r\Theta(y_k)\frac{\partial}{\partial y_k}.
\end{equation}
In fact, the difference $\Theta'$ of the two sides vanishes on $K[y_1,\dots,y_r]$ and maps $\mathfrak{n}^{m+1}$ into $\mathfrak{n}^m$ for all $m\geq1$: for $g\in K[[y_1,\dots,y_r]]$, writing $g=g_{\leq m}+g_{>m}$ with $g_{\leq m}\in K[y_1,\dots,y_r]$ and $g_{>m}\in\mathfrak{n}^{m+1}$, we get $\Theta'(g)=\Theta'(g_{>m})\in\mathfrak{n}^m$ for all $m$.
By Krull's intersection theorem \cite[\href{https://stacks.math.columbia.edu/tag/00IP}{Tag 00IP}]{SP}, it follows that
\[
\Theta'(g)\in\bigcap_{m\geq1}\mathfrak{n}^m=0.
\]
Hence $\Theta'=0$, which proves the claim.

Put $a_{ik}:=\hat{\partial}_i(y_k)$.
For the image of $x_j$ in $\hat{\cO}_{S,s}$,
\[
\delta_{ij}=\hat{\partial}_i(x_j)\overset{\eqref{eq:theta-expansion}}{=}\sum_{k=1}^ra_{ik}\frac{\partial x_j}{\partial y_k}\qquad(1\leq i,j\leq r),
\]
i.e.,
\[
(a_{ik})_{i,k}\cdot\Bigl(\frac{\partial x_j}{\partial y_k}\Bigr)_{k,j}=(\delta_{ij})_{i,j}.
\]
So $\det(a_{ik})$ is a unit of $\hat{\cO}_{S,s}$, and the matrix $(a_{ik})$ is invertible.
For $b\in C$, by assumption,
\[
0=\hat{\partial}_i\bigl(\chi'(b)\bigr)\overset{\eqref{eq:theta-expansion}}{=}\sum_{k=1}^r a_{ik}\frac{\partial\,\chi'(b)}{\partial y_k}\qquad(1\leq i\leq r).
\]
As $(a_{ik})$ is invertible, $\partial\chi'(b)/\partial y_k=0$ for $1\leq k\leq r$.
As $\operatorname{char}K=0$,
\[
\bigcap_{k=1}^r\mathrm{Ker}\Bigl(\frac{\partial}{\partial y_k}\Bigr)=K\quad\text{in }K[[y_1,\dots,y_r]],
\]
so $\chi'(b)\in K$.
As $\chi'$ is local, $\chi'\bigl(b-q_C(b)\bigr)\in\mathfrak{n}\cap K=0$, so
\[
\chi'=\bigl(C\xrightarrow{q_C}K\subset\hat{\cO}_{S,s}\bigr).
\]
\end{proof}

\begin{theorem}[Vanishing KS $\Rightarrow$ $\pi$ isotrivial]
\label{thm:ks-isotrivial}
In the setting of \ref{para:prorep}, assume that
\begin{enumerate}[label=\textup{(\alph*)}]
\item\label{item:ksiso-a} $\operatorname{char}K=0$;
\item $S$ is smooth over $K$;
\item\label{item:ksiso-c} $\pi$ is smooth and projective;
\item the Kodaira--Spencer class $\mathrm{KS}_{\pi}\colon T_{S/K}\to R^1\pi_*T_{X/S}$ vanishes on an open neighborhood of $s$.
\end{enumerate}
Then there exist an \'{e}tale morphism $S'\to S$ with $S'$ connected and of finite type over $K$, and a $K$-point $s'\in S'(K)$ over $s$, such that, over $S'$,
\[
X\times_S S'\simeq Y\times_K S'.
\]
\end{theorem}

\begin{proof}
Replacing $S$ by an affine open neighborhood of $s$ contained in it, which is still smooth over $K$ and over which $\pi$ is still smooth and projective, we may assume $\mathrm{KS}_{\pi}=0$ on the whole $S$.
As a base change of the projective morphism $\pi$, $Y$ is projective over $K$.
As $S$ is of finite type over $K$, it is noetherian; let $Y_S$, $I$, and $e\in I(K)$ be as in \Cref{lem:tangent-aut-etale}, applied to $\pi$.

Since $\mathrm{KS}_{\pi}=0$ on $S$, \Cref{lem:ks-formal-derivative} implies that $\hat{\partial}\circ\chi=0$ as a map $A\to\hat{\cO}_{S,s}$.
As $\operatorname{char}K=0$, \Cref{lem:const-formal} applied with $C:=A$ and $\chi':=\chi$ gives
\[
\chi=\bigl(A\xrightarrow{q_A}K\subset\hat{\cO}_{S,s}\bigr).
\]
Let $\alpha_m\colon Y\times_K S_m\xlongrightarrow{\sim}X_m$ be as in \Cref{lem:def-classifying}.
Because $I$ represents $\mathfrak{Isom}_S(Y_S,X)$,
\[
\mathrm{Hom}_S(S_m,I)=\mathrm{Isom}_{S_m}\bigl(Y\times_K S_m,X_m\bigr).
\]
The compatible system $\alpha_m$ hence defines compatible $S$-morphisms $S_m\to I$ sending the closed point of $S_m$ to $e$.
As $S_m=\Spec B_m$ has a single point, each of them corresponds to a local $K$-algebra map
\[
\cO_{I,e}\to B_m,
\]
whose composition with $\cO_{S,s}\to\cO_{I,e}$ is the quotient map $\cO_{S,s}\to B_m$, as $S_m\to I$ is an $S$-morphism.
As this map sends $\mathfrak{m}_e$ into $\mathfrak{m}_{B_m}$ and $\mathfrak{m}_{B_m}^m=0$, it factors through
\[
\cO_{I,e}/\mathfrak{m}_e^m\to B_m.
\]
As $S_{m+1}\to I$ restricts to $S_m\to I$, these maps are compatible with the transition maps
\[
\cO_{I,e}/\mathfrak{m}_e^{m+1}\to\cO_{I,e}/\mathfrak{m}_e^m\quad\text{and}\quad B_{m+1}\to B_m.
\]
Taking limits gives a local $K$-algebra map
\[
\sigma\colon\hat{\cO}_{I,e}\to\hat{\cO}_{S,s}.
\]
For every $m$, the composition $B_m\to\cO_{I,e}/\mathfrak{m}_e^m\to B_m$ is $\mathrm{id}_{B_m}$, as the composition of $\cO_{I,e}\to B_m$ with $\cO_{S,s}\to\cO_{I,e}$ is the quotient map.
Taking limits, $\sigma\circ(\hat{\cO}_{S,s}\to\hat{\cO}_{I,e})=\mathrm{id}$, i.e., $\sigma$ is a formal section of $I\to S$ through $e$.
From this, together with $H^0(Y,T_Y)=0$ by assumption \hyperref[item:prorep-b]{(b)} of \Cref{para:prorep}, we deduce by \Cref{lem:tangent-aut-etale} that $I\to S$ is \'{e}tale at $e$.

Let $S'\subset I$ be a connected affine open neighborhood of $e$ that is \'{e}tale over $S$, and put $s':=e\in S'(K)$.
By the representability of $I$, $\mathrm{id}_I\in\mathrm{Hom}_S(I,I)$ corresponds to an isomorphism
\[
Y_S\times_S I\xrightarrow{\sim}X\times_S I
\]
over $I$.
Its base change along $S'\subset I$ is an isomorphism
\[
Y\times_K S'\xrightarrow{\sim}X\times_S S'
\]
over $S'$.
\end{proof}

\begin{corollary}
\label{cor:ks-iff}
In the setting of \ref{para:prorep}, assume \hyperref[item:ksiso-a]{(a)}--\hyperref[item:ksiso-c]{(c)} of \Cref{thm:ks-isotrivial}.
Then the following are equivalent.
\begin{enumerate}
\item The Kodaira--Spencer class $\mathrm{KS}_{\pi}$ vanishes on an open neighborhood of $s$.
\item There exist an \'{e}tale morphism $S'\to S$ with $S'$ connected and of finite type over $K$, and a $K$-point $s'\in S'(K)$ over $s$, such that, over $S'$,
\[
X\times_S S'\simeq Y\times_K S'.
\]
\end{enumerate}
\end{corollary}

\begin{proof}
(1)$\Rightarrow$(2): This is \Cref{thm:ks-isotrivial}.

(2)$\Rightarrow$(1): Let $g\colon S'\to S$ and $s'$ be as in (2), and let $\pi'\colon X\times_S S'\to S'$ be the base change of $\pi$.
\[
\begin{tikzcd}[column sep=3.6em, row sep=2.6em, ampersand replacement=\&]
X\times_SS' \arrow[r] \arrow[d, "\pi'"'] \& X \arrow[d, "\pi"] \\
S' \arrow[r, "g"] \& S
\end{tikzcd}
\]
As $g$ is \'{e}tale, $T_{S'/K}=g^*T_{S/K}$, and
\[
R^1\pi'_*T_{X\times_SS'/S'}=g^*R^1\pi_*T_{X/S}
\]
by flat base change \cite[\href{https://stacks.math.columbia.edu/tag/02KH}{Tag 02KH}]{SP}; under these identifications $\mathrm{KS}_{\pi'}\colon T_{S'/K}\to R^1\pi'_*T_{X\times_SS'/S'}$ equals $g^*\mathrm{KS}_{\pi}$, i.e., the diagram
\[
\begin{tikzcd}[column sep=4.4em, row sep=2.6em, ampersand replacement=\&]
T_{S'/K} \arrow[r, "\mathrm{KS}_{\pi'}"] \arrow[d, "\wr"'] \& R^1\pi'_*T_{X\times_SS'/S'} \arrow[d, "\wr"] \\
g^*T_{S/K} \arrow[r, "g^*\mathrm{KS}_{\pi}"] \& g^*R^1\pi_*T_{X/S}
\end{tikzcd}
\]
commutes.
By (2), the $S'$-isomorphism $X\times_SS'\simeq Y\times_KS'$ identifies the sequences \eqref{eq:ks-ses} for $\pi'$
\[
0\longrightarrow T_{X\times_SS'/S'}\longrightarrow T_{X\times_SS'/K}\longrightarrow \pi'^*T_{S'/K}\longrightarrow0
\]
and for the constant family $\mathrm{pr}_2\colon Y\times_KS'\to S'$
\begin{equation}\label{eq:ks-ses-const}
0\longrightarrow T_{Y\times_KS'/S'}\longrightarrow T_{Y\times_KS'/K}\longrightarrow \mathrm{pr}_2^*T_{S'/K}\longrightarrow0,
\end{equation}
so $\mathrm{KS}_{\pi'}$ is the Kodaira--Spencer class of the constant family.
The sequence \eqref{eq:ks-ses-const} is split, so $\mathrm{KS}_{\pi'}=0$, i.e., $g^*\mathrm{KS}_{\pi}=0$.
As $g$ is \'{e}tale, $g(S')$ is an open neighborhood of $s$ and $\{S'\to g(S')\}$ is an fpqc covering, so $\mathrm{KS}_{\pi}=0$ on $g(S')$ by fpqc descent \cite[\href{https://stacks.math.columbia.edu/tag/023T}{Tag 023T}]{SP}.
\end{proof}

\section{Integrality on the Hodge line}
\label{sec:dropping}

\begin{lemma}[Isotriviality of $\nabla_{\GM}$]
\label{lem:isotrivial-hodge-line}
Let $K$, $R$, $S$, $\pi$, $(\cE,\nabla)$, and $s$ be as in \Cref{thm:main}, and we assume
\begin{enumerate}[leftmargin=2em,labelwidth=1.6em,labelsep=0.4em,labelindent=0pt,align=left]
\item[\textup{(a$'$)}] $\theta\colon F_{\mathrm{Hdg}}^2\cE\to\bigl(F_{\mathrm{Hdg}}^1\cE/F_{\mathrm{Hdg}}^2\cE\bigr)\otimes_{\cO_S}\Omega^1_{S/R}$ is identically zero.
\item[\textup{(b)}] $s\in S(R)$ and $v\colon\Spec R\to E$ is a morphism over $s$ that factors through $F_{\mathrm{Hdg}}^2E$ and is not the zero section.
\end{enumerate}
Let $Y:=X_{s_K}$ be the fiber of $\pi$ over $s_K$, with structure morphism $\pi_{s_K}\colon Y\to\Spec K$.
Then there exist
\begin{enumerate}[label=\textup{(\roman*)},leftmargin=2em,labelwidth=1.6em,labelsep=0.4em,labelindent=0pt,align=left]
\item an \'{e}tale morphism $S'\to S_K$ with $S'$ connected and of finite type over $K$, and
\item a $K$-point $s'\in S'(K)$ over $s_K$,
\end{enumerate}
such that the following hold:
\begin{enumerate}[leftmargin=2em,labelwidth=1.6em,labelsep=0.4em,labelindent=0pt,align=left]
\item
The base change of $\pi$ to $S'$ is isomorphic to the constant family $Y\times_{K}S'\to S'$:
\[
X_K\times_{S_K}S'\simeq Y\times_{K}S'.
\]
\item\label{item:isohl-2}
The isomorphism in (1) induces an isomorphism
\[
(\cE,\nabla)_K|_{S'}\simeq\bigl(h^*R^2\pi_{s_K*}\Omega^\bullet_{Y/K},1\otimes d\bigr),
\]
where $h\colon S'\to\Spec K$ is the structure morphism.
\item\label{item:isohl-3}
This isomorphism maps $F_{\mathrm{Hdg}}^2\cE_K|_{S'}$ onto $h^*\pi_{s_K*}\Omega^2_{Y/K}$.
\end{enumerate}
\end{lemma}

\begin{proof}
\begin{figure}[ht]
\centering
\resizebox{0.85\textwidth}{!}{%
\begin{tikzpicture}[>=Stealth,line cap=round,line join=round]
\def\dx{1.6}\def\dy{1.8}
\draw[thick] (0,0) -- (8,0) node[right]{$\Spec R$};
\foreach \x/\l in {1/{\mathfrak p_1},2.6/{\mathfrak p_2},4.2/{\mathfrak p_3}}{\fill (\x,0) circle (1.5pt) node[below]{\small $\l$};}
\draw[thick,fill=white] (6.6,0) ellipse [x radius=5pt, y radius=2.6pt] node[below=2pt]{\small $\eta=\Spec K$};
\begin{scope}[yshift=1.1cm]
\fill[gray!10] (0,0) -- (8,0) -- (8+\dx,\dy) -- (\dx,\dy) -- cycle;
\begin{scope}
\clip (0,0) -- (8,0) -- (8+\dx,\dy) -- (\dx,\dy) -- cycle;
\foreach \w in {0.04,0.08,...,0.36}{\fill[red!80!black,opacity=0.09] (6.6-\w,0) -- (6.6+\w,0) -- (6.6+\dx+\w,\dy) -- (6.6+\dx-\w,\dy) -- cycle;}
\end{scope}
\draw[thick] (0,0) -- (8,0) -- (8+\dx,\dy) -- (\dx,\dy) -- cycle;
\foreach \x in {1,2.6,4.2}{\draw[gray] (\x,0) -- (\x+\dx,\dy);}
\node[gray] at (2.6+0.75*\dx,0.75*\dy) [right] {\small $S_{k(\mathfrak p_2)}$};
\node[red!70!black] at (6.6+\dx+0.1,\dy) [above right]{$S_K$};
\draw[blue!60!black,very thick] (0.2,0.35) .. controls (2.5,1.2) and (4.5,0.4) .. (6.6+0.45*\dx,0.45*\dy) .. controls (7.6,1.0) .. (8.4,0.95) node[right]{$s$};
\fill[blue!60!black] (6.6+0.45*\dx,0.45*\dy) circle (2.2pt);
\node[blue!60!black] at (6.6+0.45*\dx-0.15,0.45*\dy+0.2) [left] {\small $s_K$};
\node at (0.3,1.1) {$S$};
\end{scope}
\draw[->] (4,0.95) -- (4,0.2);
\begin{scope}[yshift=4.3cm]
\def\Xout{(0,0) .. controls (1.2,-0.5) and (2.2,0.35) .. (3.4,0.05) .. controls (4.6,-0.25) and (5.6,0.45) .. (6.6,0.1) .. controls (7.3,-0.15) and (7.7,-0.1) .. (8,0)
  .. controls (8.7,0.4) and (8.4,1.1) .. (8.9,1.3) .. controls (9.3,1.45) and (9.4,1.6) .. (8+\dx,\dy)
  .. controls (10.2,2.3) and (9.2,2.8) .. (9.8,3.3) .. controls (10.1,3.6) and (9.7,3.9) .. (8+\dx,\dy+2.4)
  .. controls (8.4,4.7) and (7.3,3.8) .. (5.8,4.25) .. controls (4.4,4.6) and (3.2,3.9) .. (\dx,\dy+2.4)
  .. controls (1.0,3.7) and (1.1,3.1) .. (0.5,2.9) .. controls (0.1,2.75) and (0.2,2.55) .. (0,2.4)
  .. controls (-0.5,1.8) and (0.35,1.4) .. (-0.1,0.9) .. controls (-0.4,0.5) and (0.2,0.3) .. (0,0) -- cycle}
\fill[gray!6] \Xout;
\draw[thin,dashed,gray!60!black] (\dx,\dy) .. controls (1.25,2.5) and (1.95,3.3) .. (\dx,\dy+2.4);
\draw[thin,dashed,gray!60!black] (0,0) .. controls (0.6,0.5) and (0.9,1.4) .. (\dx,\dy);
\draw[thin,dashed,gray!60!black] (\dx,\dy) .. controls (2.6,1.6) and (3.4,1.9) .. (4.2,1.75) .. controls (5.4,1.55) and (6.8,2.0) .. (8.2,1.8) .. controls (8.8,1.7) and (9.2,1.75) .. (8+\dx,\dy);
\begin{scope}
\clip \Xout;
\def\Sl{(2.6,0.105) .. controls (2.95,0.9) and (2.3,1.7) .. (2.6,2.502) .. controls (3.2,2.95) and (3.6,3.9) .. (4.2,4.3) .. controls (4.5,3.5) and (3.9,2.6) .. (4.2,1.75) .. controls (3.6,1.2) and (3.2,0.6) .. (2.6,0.105) -- cycle}
\fill[gray,opacity=0.3] \Sl;
\draw[gray!70!black] \Sl;
\def\Kl{(6.6,0.1) .. controls (6.95,0.9) and (6.3,1.7) .. (6.6,2.533) .. controls (7.2,2.95) and (7.6,3.9) .. (8.2,4.31) .. controls (8.5,3.5) and (7.9,2.6) .. (8.2,1.8) .. controls (7.6,1.2) and (7.2,0.6) .. (6.6,0.1) -- cycle}
\foreach \w in {-0.36,-0.32,...,0.37}{\fill[red!80!black,opacity=0.045,xshift=\w cm] \Kl;}
\end{scope}
\draw[blue!60!black,very thick,fill=blue!25] (6.6+0.45*\dx,0.45*\dy+1.2) ellipse (0.28 and 0.95);
\node[blue!60!black] at (6.6+0.45*\dx+0.05,0.45*\dy+1.2) {\small $Y$};
\draw[thick] \Xout;
\draw[thick] (0,2.4) .. controls (1.3,2.0) and (2.4,2.75) .. (3.6,2.45) .. controls (4.8,2.15) and (6.2,2.8) .. (8,2.4) .. controls (8.6,2.9) and (9.0,3.8) .. (8+\dx,\dy+2.4);
\draw[thick] (8,0) .. controls (8.4,0.8) and (7.7,1.5) .. (8,2.4);
\node[gray!60!black] at (2.6+\dx+0.1,\dy+2.4) [above] {$X_{k(\mathfrak p_2)}$};
\node[red!70!black] at (6.6+\dx+0.1,\dy+2.4) [above] {$X_K$};
\node at (0.4,2.0) {$X$};
\end{scope}
\draw[->] (4,4.15) -- node[right]{$\pi$} (4,3.05);
\begin{scope}[xshift=11.4cm,yshift=1.1cm]
\draw[green!40!black,very thick] (0,0.2) .. controls (0.6,0.9) and (1.0,1.0) .. (1.6,1.6) node[above right]{$S'$};
\fill[green!40!black] (0.72,0.81) circle (2.2pt) node[right]{\small $s'$};
\begin{scope}[yshift=3.2cm]
\fill[green!8] (0,0.2) .. controls (0.6,0.9) and (1.0,1.0) .. (1.6,1.6) -- (1.6,3.8) .. controls (1.0,3.2) and (0.6,3.1) .. (0,2.4) -- cycle;
\draw[green!40!black,thick] (0,0.2) .. controls (0.6,0.9) and (1.0,1.0) .. (1.6,1.6) -- (1.6,3.8) .. controls (1.0,3.2) and (0.6,3.1) .. (0,2.4) -- cycle;
\foreach \t/\px/\py in {0.15/0.25/0.48,0.75/1.2/1.25}{\draw[blue!50!black,fill=blue!10] (\px,\py+1.1) ellipse (0.22 and 0.75);}
\draw[blue!60!black,very thick,fill=blue!25] (0.72,0.81+1.1) ellipse (0.22 and 0.75);
\node[blue!60!black] at (0.72,1.91) {\scriptsize $Y$};
\node[green!40!black,align=left] at (2.0,4.3) {$X_{S'}\simeq Y\times_K S'$};
\end{scope}
\draw[->] (0.8,3.0) -- node[right]{\small $\pi_{S'}=\mathrm{pr}_2$} (0.8,2.0);
\end{scope}
\draw[->,green!40!black,thick] (11.2,2.0) .. controls (10.4,2.4) and (9.0,2.9) .. (7.75,2.45) node[pos=0.3,above,sloped]{$g$ \small (\'etale)};
\draw[->,green!40!black,thick] (11.2,6.0) .. controls (10.5,6.6) and (9.2,7.6) .. (7.85,7.05);
\draw[->,green!40!black,thick] (11.3,1.2) .. controls (10.2,0.55) and (8.2,0.45) .. (6.82,0.1) node[pos=0.35,above,sloped]{$h$};
\end{tikzpicture}%
}
\captionsetup{labelformat=simple,labelsep=period}
\caption{The schemes involved.}
\label{fig:schemes-61}
\end{figure}
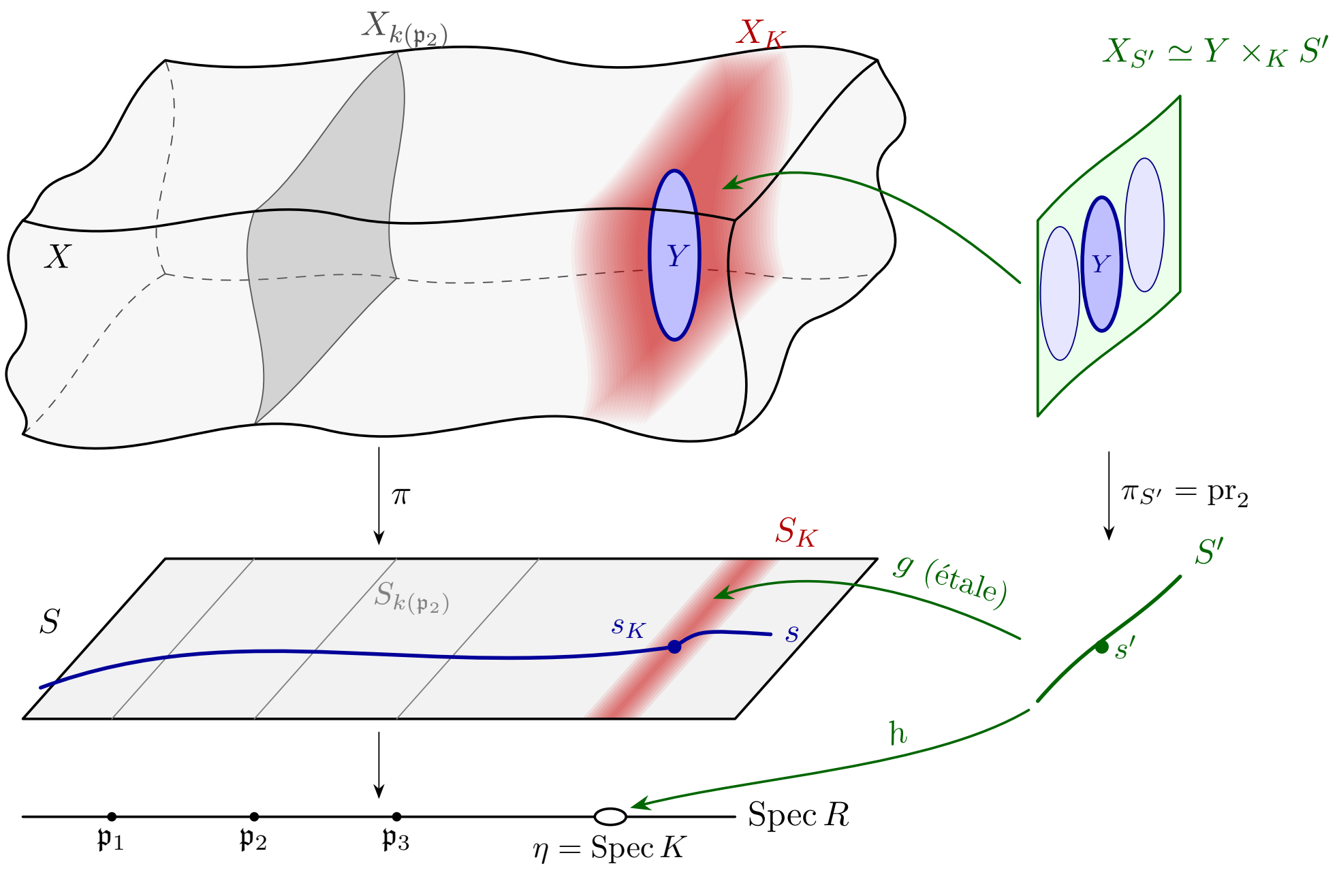
After replacing $R$ by $R[1/a]$ as in \Cref{rmk:local-freeness}, which does not change $S_K$ or $s_K$, we can assume
\begin{enumerate}[label=\textup{(\arabic*)},ref=\arabic*,labelindent=0pt,labelwidth=\widthof{\textup{(LFBS)}},labelsep=0.5em,leftmargin=!,align=left]
\item[\textup{\hyperref[item:lf-hodge]{(LFBS)}}]\phantomsection\label{item:61-lfbs} each $R^j\pi_*\Omega^i_{X/S}$ is locally free, and its formation commutes with base change;
\item[\textup{\hyperref[item:lf-degen]{(HdR)}}]\phantomsection\label{item:61-hdr} the Hodge--de Rham spectral sequence \eqref{eq:HdR-f} degenerates at $E_1$,
\end{enumerate}
by \Cref{rmk:local-freeness}(\ref{item:lf-221-lfbs}) and (\ref{item:lf-221-hdr}).
\begin{itemize}
\item On the one hand, \cite[(1.4.1.7.2)]{Kat72} for $(p,q)=(2,0)$ yields, for a local generator $\sigma$ of $\pi_*\Omega^2_{X/S}$,
\[
\theta(\partial)(\sigma)=\mathrm{KS}_{\pi}(\partial)\cup\sigma;
\]
\item On the other hand, the perfect pairing $\Omega^1_{X/S}\otimes\Omega^1_{X/S}\xrightarrow{\wedge}\Omega^2_{X/S}=\cO_X\cdot\sigma$ (because $X$ is a $K3$ family) gives an isomorphism $T_{X/S}\xrightarrow{\sim}\Omega^1_{X/S}$, $\xi\mapsto\sigma(\xi,-)$, which induces an isomorphism 
\[R^1\pi_*T_{X/S}\xrightarrow{\sim}R^1\pi_*\Omega^1_{X/S}\] 
given by $c\mapsto c\cup \sigma$.

\end{itemize}
Therefore $\theta=0$ implies $\mathrm{KS}_{\pi}(\partial)=0$ for all local sections $\partial$ of $T_{S/R}$, hence
\[\mathrm{KS}_{\pi}=0.\]
Write $\pi_K\colon X_K\to S_K$ for the base change of $\pi$.
By flat base change \cite[\href{https://stacks.math.columbia.edu/tag/02KH}{Tag 02KH}]{SP} along $\Spec K\to\Spec R$, $\mathrm{KS}_{\pi}=0$ implies
\[
\mathrm{KS}_{\pi_K}=0.
\]
As $s\in S(R)$, the induced point $s_K$ is a $K$-point of $S_K$.
The fiber $Y$ is a $K3$ surface over $K$, so $H^0(Y,T_Y)=0$ by \cite[Theorem 9.5.1]{Huy16}.
Since $H^0(Y,T_Y)=0$ and $\mathrm{KS}_{\pi_K}=0$, \Cref{thm:ks-isotrivial} applies to $X_K\to S_K$ and the point $s_K$, and gives
\begin{enumerate}[label=\textup{(\roman*)},leftmargin=2em,labelwidth=1.6em,labelsep=0.4em,labelindent=0pt,align=left]
\item an \'{e}tale morphism $S'\to S_K$ with $S'$ connected and of finite type over $K$, and
\item a point $s'\in S'(K)$ over $s_K$,
\end{enumerate}
such that
\begin{enumerate}[leftmargin=2em,labelwidth=1.6em,labelsep=0.4em,labelindent=0pt,align=left]
\item $X_K\times_{S_K}S'\simeq Y\times_K S'$ over $S'$.
\end{enumerate}
Write $g\colon S'\to S_K$ for the \'{e}tale morphism, $\pi_{S'}\colon X_{S'}:=X_K\times_{S_K}S'\to S'$ for the base change of $\pi_K$.
Since $g$ is flat, the base change map 
\[g^*\cE_K=g^*R^2\pi_{K*}\Omega^\bullet_{X_K/S_K}\longrightarrow R^2\pi_{S'*}\Omega^\bullet_{X_{S'}/S'}=:\cE_{S'}\]
is an isomorphism.
By functoriality of $\nabla_{\GM}$, so does the formation of $\nabla_{\GM}$, namely the diagram
\[
{\begin{tikzcd}[column sep=4.4em, row sep=2.6em, ampersand replacement=\&]
g^*\cE_K \arrow[r,"g^*\nabla_{\GM}"] \arrow[d,"\wr"'] \& g^*\cE_K\otimes_{\cO_{S'}}\Omega^1_{S'/K} \arrow[d,"\wr"] \\
\cE_{S'} \arrow[r,"\nabla_{\GM}"] \& \cE_{S'}\otimes_{\cO_{S'}}\Omega^1_{S'/K}
\end{tikzcd}}
\]
commutes.
On the constant family $Y\times_{K}S'$, \Cref{lem:dR-flat-bc}(\ref{item:dR-const}) applied to $\pi_{s_K}\colon Y\to\Spec K$ and $h\colon S'\to\Spec K$ gives
\[
(R^2\pi_{S'*}\Omega^\bullet_{Y\times_KS'/S'},\nabla_{\GM})\simeq(h^*R^2\pi_{s_K*}\Omega^\bullet_{Y/K},1\otimes d).
\]
This proves (2).
By \hyperref[item:61-hdr]{(HdR)} and \Cref{lem:HdR-degen}, the natural map $\pi_*\Omega^2_{X/S}=R^2\pi_*\Omega^{\geq2}_{X/S}\to\cE$ is injective with image $F_{\mathrm{Hdg}}^2\cE$.
By flat base change, $g^*\pi_{K*}\Omega^2_{X_K/S_K}\simeq\pi_{S'*}\Omega^2_{X_{S'}/S'}$ compatibly with the isomorphism $g^*\cE_K\xrightarrow{\sim}\cE_{S'}$.
As $g^*$ is exact, it preserves images, so this isomorphism maps $g^*F_{\mathrm{Hdg}}^2\cE_K$ onto $F_{\mathrm{Hdg}}^2\cE_{S'}$. On the constant family, $F_{\mathrm{Hdg}}^2$ is the trivial line subbundle $h^*\pi_{s_K*}\Omega^2_{Y/K}$.
This proves (3).
\end{proof}

\begin{theorem}
\label{thm:hodge-line-iff}
Let $K$, $R$, $S$, $\pi$, $s$, and $v$ be as in \Cref{thm:main}, but without assumption (a).
Write $\tilde{v}_K$ for the formal horizontal section of $(\cE,\nabla)_K$ through $v_K$.
The following are equivalent.
\begin{enumerate}
\item
The Higgs field $\theta\colon F_{\mathrm{Hdg}}^2\cE\to\bigl(F_{\mathrm{Hdg}}^1\cE/F_{\mathrm{Hdg}}^2\cE\bigr)\otimes_{\cO_S}\Omega^1_{S/R}$ is identically zero.
\item
The section $\tilde{v}_K$ is algebraic.
\item
The section $\tilde{v}_K$ is integral.
\item
The section $\tilde{v}_K$ is $\omega(p)$-integral.
\end{enumerate}
\end{theorem}

\begin{proof}
The implications $(2)\Rightarrow(3)$ and $(3)\Rightarrow(4)$ are \Cref{lem:alg-implies-int} and \Cref{lem:int-implies-omega}.
The implication $(4)\Rightarrow(1)$ is the contrapositive of \Cref{thm:main}.

It remains to show $(1)\Rightarrow(2)$.
If $\theta\equiv0$, then $\nabla_{\GM}$ preserves $F_{\mathrm{Hdg}}^2\cE$:
\[
\nabla_{\GM}\bigl(F_{\mathrm{Hdg}}^2\cE\bigr)\subset F_{\mathrm{Hdg}}^2\cE\otimes_{\cO_S}\Omega^1_{S/R}.
\]
So $\nabla_{\GM}$ restricts to a connection on $F_{\mathrm{Hdg}}^2\cE$, and $\tilde{v}_K$ is the formal horizontal section of the rank-one connection $(F_{\mathrm{Hdg}}^2\cE,\nabla_{\GM})_K$ through $v_K$.
{By \Cref{lem:isotrivial-hodge-line}\hyperref[item:isohl-2]{(2)}\hyperref[item:isohl-3]{(3)}, there exist an \'{e}tale morphism $h\colon$$S'\to S_K$ with $S'$ connected and of finite type over $K$, and a $K$-point $s'\in S'(K)$ over $s_K$, such that
\[
(F_{\mathrm{Hdg}}^2\cE,\nabla_{\GM})_K|_{S'}\simeq(\cO_{S'},d).
\]}%
Since $s'$ is a $K$-point, this isomorphism identifies the value $v_K$ at $s'$ with a scalar $c\in K$.
Let $u\in\Gamma\bigl(S',F_{\mathrm{Hdg}}^2\cE_K|_{S'}\bigr)$ be the section corresponding to the constant function $c$; then $\nabla_{\GM}u=0$ because $dc=0$.
Let $\hat h\colon\hat S'_{s'}\to\hat S_{s_K}$ be the morphism induced by $h$.
Then $u|_{\hat S'_{s'}}$ and $\tilde v_K\circ\hat h$ are formal horizontal sections of $(F_{\mathrm{Hdg}}^2\cE,\nabla_{\GM})_K$ over $\hat S'_{s'}$ with the same value $v_K$ at $s'$.
By uniqueness of formal horizontal sections \cite[Proposition 2.2.1]{LL25}, they coincide, i.e., the diagram \eqref{eq:cover-square} commutes for $u$.
So the algebraicity criterion \Cref{lem:cover-algebraic}, applied to $(F_{\mathrm{Hdg}}^2\cE,\nabla_{\GM})$, $s$, $v$, $h$, $s'$, and $u$, shows that $\tilde{v}_K$ is algebraic as a formal horizontal section of $(F_{\mathrm{Hdg}}^2\cE,\nabla_{\GM})_K$.
Since $F_{\mathrm{Hdg}}^2E_{\C}$ is a closed subscheme of $E_{\C}$, $\tilde{v}_K$ is algebraic as a formal horizontal section of $(\cE,\nabla_{\GM})_K$.
\end{proof}

\end{document}